\documentclass[a4paper,12pt,wlayout]{article}

\usepackage{array}
\usepackage{amssymb}
\usepackage{amsmath}
\usepackage{amsbsy}
\usepackage{epsfig}
\usepackage{mathrsfs}
\usepackage{graphics}
\usepackage{graphicx}
\usepackage[sort&compress,numbers]{natbib}
\usepackage{pstricks}
\usepackage{float}
\usepackage{amsthm}
\usepackage{slashed}

\newcommand{\eps}{\varepsilon}
\newcommand{\lb}{\label}

\newcommand{\go}{\rightarrow}

\newcommand{\ee}{\end{equation}}
\newcommand{\be}{\begin{equation}}
\newcommand{\bea}{\begin{eqnarray}}
\newcommand{\eea}{\end{eqnarray}}
\newcommand{\sbea}{\begin{subequations}\begin{eqnarray}}
\newcommand{\seea}{\end{eqnarray}\end{subequations}}
\newcommand{\ees}{\end{equation*}}
\newcommand{\bes}{\begin{equation*}}
\newcommand{\beas}{\begin{eqnarray*}}
\newcommand{\eeas}{\end{eqnarray*}}

\newcommand{\rf}[1]{(\ref{#1})}

\newcommand{\const}{\mathrm{const}}

\newcommand{\tb}{\color{blue} }
\providecommand{\keywords}[1]
{
	\small	
	\textbf{\textit{Keywords---}} #1
}

\begin{document}

\title{Composite solutions to a liquid bilayer model.}

\author{Georgy Kitavtsev\thanks{Middle East Technical University, Northern Cyprus Campus, Kalkanli, Güzelyurt, KKTC, 10 Mersin, Turkey. {\it email: georgy.kitavtsev@gmail.com}}}

\maketitle

\begin{abstract}
This article continues the research initiated in~\cite{JHKPW13}. We derive explicit formulae for the leading order profiles of eleven types of stationary solutions to a one-dimensional two-layer thin-film liquid model considered with an intermolecular potential depending on both layer heights. The found solutions are composed of the repeated elementary blocks (bulk, contact line and ultra-thin film ones) being consistently asymptotically matched together. We show that once considered on a finite interval  with Neumann boundary conditions these stationary solutions are either dynamically stable or weakly translationally unstable. Other composite solutions are found to be numerically unstable and rather exhibit complex coarsening dynamics.
\end{abstract}


%%%%%%%%%%%%%%%%%%%%%%%%%%%%%%%%%%%%%%%%%%%%%%%%%%%%%%%%%%%%%%%%%
\section{Introduction}
%%%%%%%%%%%%%%%%%%%%%%%%%%%%%%%%%%%%%%%%%%%%%%%%%%%%%%%%%%%%%%%%%

Evolution of polymer liquid films with micro- to nanoscale thickness governed by nonlinear competition of bulk, surface and intermolecular forces exhibits complex evolving patterns. Understanding and control of them is important question for numerous applications in printing and coating technologies, optoelectronics, cell dynamics and drug delivery problems in bio-microfluidics~\cite{CM09,BEIMR11}. In the last three decades, considerable amount of modeling and experimental work has been conducted on the related problem of evolution of two-layer immiscible thin liquid films deposed on a solid substrate. These so called bilayer liquid films demonstrate rich dynamics with complex evolving morphological structures  which  currently are partially understood even in the case when both liquid layers are Newtonian~\cite{PBMT04,PBMT05,BSTA21,Ne21}.

Mathematically seen {\it lubrication approximations} applied starting from the pioneering works~\cite{ODB97,Danov98,FG05,BGS05} showed that bilayer thin films can be effectively modeled by the coupled systems of fourth-order quasi-linear parabolic PDEs describing evolution of the liquid layer heights. Analytical and numerical properties of solutions to these system are complicated by the involved nonlinear terms degenerating when one of the layer heights shrinks to zero (a process called {\it rupture}). Here an important role is played by so called {\it disjoining pressure} potentials~\cite{Is92,SHJ01,GW03} which account for the action of intermolecular forces and prevent solutions of the bilayer systems from rupture: a small characteristic length scale parameter $\eps>0$ present in the intermolecular potentials typically sets a lower bound for the minima of the liquid layer heights.  Numerical simulations of the bilayer systems revealed existence of several stationary solutions with different non-constant profiles as well as complex dynamical patterns involving multiple coarsening stages of so called metastable (quasi-stationary or slow evolving) solutions~\cite{PBMT04,PBMT05,JPMW14,HJKP15}. 

Following seminal study of stationary solutions to single layer thin film equations~\cite{BGW01}, in~\cite{JHKPW13} authors considered a symmetric form of the intermolecular potential (see \rf{phi} below) depending on the height of the top liquid layer and showed that the non-constant stable stationary solutions (under broad types of boundary conditions) to the underlying one-dimensional bilayer system have a lens shape characterized by one liquid layer forming a sessile drop on top of the second one (cf.  Fig.2 {\bf (a)}), see also~\cite{KM03,CM06}. 

This article produces further steps in investigation of stationary solutions to the bilayer systems by considering a more realistic intermolecular potential \rf{phi_c} depending on the heights of both liquid layers. It is similar to the intermolecular potential considered in ~\cite{PBMT04,PBMT05} except that the repulsive terms in \rf{phi_c}--\rf{phi} have an algebraic rather than exponential decay for large heights. We demonstrate that on the one hand this potential modification results in a more reach and complex morphology of both stationary and dynamical solutions exhibiting multiple {\it triple phase contact lines}. The latter are tiny regions where three out of four participating phases (the two liquids, ambient atmosphere and solid substrate)  meet together in the asymptotic limit $\eps\go0$. On the other hand, classification and properties of the derived here solutions stay in good correspondence with the observations made by modeling and experimental studies of thin bilayer systems~\cite{NTGDC12,BSTA21,MAP01,PBMT04,PBMT05}. 
Following~\cite{JPMW14,JKT14,JHKPW13} we begin our study by considering a one-dimensional no-slip bilayer model:
\be
\lb{BS}
\left[\begin{array}{c}h_1\\ h\end{array}\right]_t=\frac{\partial}{\partial x}\left(\mathbf{Q}(h_1,h)\cdot\frac{\partial}{\partial x}\left[\begin{array}{c}-(\sigma+1)\frac{\partial^2 h_1}{\partial x^2}-\frac{\partial^2 h}{\partial x^2}+\frac{\partial\phi_\eps(h_1,\,h)}{\partial h_1}\\[2ex] -\frac{\partial^2 h_1}{\partial x^2}-\frac{\partial^2 h}{\partial x^2}+\frac{\partial\phi_\eps(h_1,\,h)}{\partial h}\end{array}\right]\right)
\ee
\noindent with $h_1(x,\,t)$ and $h(x,\,t)$ being the heights of the first and second liquid layers (cf. Fig.1), respectively, $\sigma$ the ratio of the surface tensions (of the lower liquid layer
%%%%%%%%%%%%%%%%%%%%%%%%%%%%%%%%%%%%%%%%%%%%%%%%%%%%%%%%%%%%%%%%%%%%%%%%%%%%%%%%%%%%%%%%%%%%%%%%%%%%%%%%%%%%%%%%
\begin{figure}[H] 
	\centering
	\vspace{-2.9cm}
	\hspace{-.4cm}\includegraphics[width=.3\textwidth]{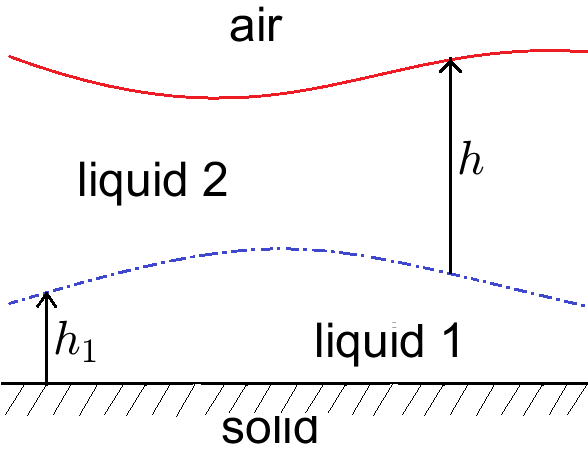}  
	\caption{\small Sketch of two liquid layer height profiles $h_1(x,\,t)$ and $h(x,\,t)$ in \rf{BS}.}	
\end{figure}
%%%%%%%%%%%%%%%%%%%%%%%%%%%%%%%%%%%%%%%%%%%%%%%%%%%%%%%%%%%%%%%%%%%%%%%%%%%%%%%%%%%%%%%%%%%%%%%%%%%%%%%%%%%%%%%%
\noindent over the upper one), and {\it the intermolecular potential}
\be
\lb{phi_c}
\phi_\eps(h_1,\,h)=\phi\left(\frac{h_1}{\eps}\right)+\phi\left(\frac{h}{\eps}\right)
\ee
with $\phi(h)$ denoting a Lenard-Jones type potential
\be
\lb{phi}
\phi(h)=\frac{1}{lh^l}-\frac{1}{nh^n}
\ee
and parameter $\eps>0$ scaling the location of its minimum. In \rf{phi}, $n$ and $l$ are the attractive and repulsive {\it positive integer exponents}, respectively, with $l>n\ge 2$. For the results of this article the concrete choice of $(n,\,l)$ is not important, the dependence on them will appear in {\it the leading order formulae} only through the absolute value of the potential minimum value (well's depth):
$$|\phi(1)|=\frac{1}{n}-\frac{1}{l}.$$
For polymeric dewetting liquid films typically used choices of $(n,\,l)$ are $(2,\,8)$ or $(2,\,3)$~\cite{Is92,SHJ01,BGW01}. In numerical simulations of this article the latter choice is used.

In \rf{BS},  $\mathbf{Q}(h_1,h)$ is the symmetric $2\times2$ {\it mobility matrix}:
\bes
\mathbf{Q}(h_1,h)=\frac{1}{\mu}\left[\begin{array}{cc}\frac{h_1^3}{3}&\frac{h_1^2h}{2}\\[2ex] \frac{h_1^2h}{2}& \quad\frac{\mu}{3}h^3+h_1h^2\end{array}\right]
\ees
with $\mu$ being the ratio of the liquid viscosities (of the lower layer over the upper one), which is set to $\mu=1$ in  numerical simulations of system \rf{BS}.

System \rf{BS} considered in interval $(-L,0)$ with Neumann boundary conditions:
\be
\lb{BC}
\frac{\partial h_1}{\partial x}=\frac{\partial h}{\partial x}=\frac{\partial^3 h_1}{\partial x^3}=\frac{\partial^3 h}{\partial x^3}\quad\text{at}\ x=-L,\,0, 
\ee
and initial conditions
\be
\lb{IC}
h_1(x,\,0)=h_{1,0}(x)\quad\text{and}\ h(x,\,0)=h_0(x).
\ee
has the unique positive classical smooth solution given sufficiently regular initial liquid height profiles $h_{1,0}(x)>0$ and $h_0(x)>0$~\cite{JKT14,EM14}.

Additionally, \rf{BS} equipped with \rf{BC} can be viewed as a gradient flow of the following energy functional~\cite{PBMT04,PBMT05,HJKP15}:
\be
\lb{En}
E(h_1,\,h)=\int_{-L}^0\left[\frac{\sigma}{2}\left|\frac{\partial h_1}{\partial x}\right|^2+\frac{1}{2}\left|\frac{\partial(h_1+h)}{\partial x}\right|^2+\phi_\eps(h_1,\,h)\right]\,dx.
\ee
The masses of two liquids are conserved in this evolution, i.e.
\sbea
\lb{MC1}
\int_{-L}^0h_1(x,\,t)\,dx&=&M_1,\\
\int_{-L}^0h(x,\,t)\,dx&=&M\quad\text{for all}\ t\ge 0
\lb{MC}
\seea
with positive constants $M_1$ and $M$ (for simplicity both liquid densities are assumed to be equal to $1$). {\it Positive stationary} solutions to system \rf{BS}--\rf{BC} are critical points of the energy functional \rf{En} and the corresponding Euler-Lagrange system can be written in the form:
\sbea
\lb{SSa}
\sigma h_1''(x)&=&\Pi_\eps(h_1(x))-\Pi_\eps(h(x))-\lambda_2+\lambda_1,\\[1ex]
\lb{SSb}
\sigma h''(x)&=&-\Pi_\eps(h_1(x))+(\sigma+1)\Pi_\eps(h(x))+\lambda_2-(\sigma+1)\lambda_1,\\[1ex]
h'(0)&=&h_1'(0)=h'(-L)=h_1'(-L)=0,
\lb{SSc}
\seea
with constants $\lambda_2$ and $\lambda_1$ being the Lagrange multipliers (or {\it constant hydrodynamic pressures} of the first and second liquid layers) associated with conservation of masses \rf{MC1} and
\rf{MC}, respectively, and  $\Pi_{\eps}(h)$ derivative of potential \rf{phi}:
\be
\Pi_\eps(h)=\frac{\partial \phi}{\partial h}\left(\frac{h}{\eps}\right)=\frac{1}{\eps}\left[\left(\frac{\eps}{h}\right)^{n+1}-\left(\frac{\eps}{h}\right)^{l+1}\right].
\lb{Pi}
\ee
Additionally, we note that \rf{SSa}--\rf{SSb} can be rewritten as
\sbea
\lb{Ms1}
(h+h_1)''&=&\Pi_\eps(h)-\lambda_1,\\[1ex]
(h+(\sigma+1)h_1)''&=&\Pi_\eps(h_1)-\lambda_2,
\lb{Ms2}
\seea
which after integration in $x\in(-L,\,0)$ using \rf{SSc} imply 
\be
\lambda_1=\tfrac{1}{L}\int_{-L}^0\Pi_\eps(h)\,dx\quad\text{and}\ \lambda_2=\tfrac{1}{L}\int_{-L}^0\Pi_\eps(h_1)\,dx.
\lb{lambda_int}
\ee
From \rf{Pi} and \rf{lambda_int} we see that both $\lambda_1$ and $\lambda_2$ are positive if $h_1(x)\ge\eps$ and $h(x)\ge\eps$ for all $x\in(-L,\,0)$. The latter conditions hold for most of the solutions to  system \rf{SSa}--\rf{SSc} derived in this article.

Furthermore, multiplying \rf{Ms1} and \rf{Ms2} by $h'(x)$ and $h_1'(x)$, respectively, and integrating them in $x\in(-L,\,0)$ by parts gives
\beas
-\int_{-L}^0h_1'h''\,dx&=&\phi(\tfrac{h(0)}{\eps})-\phi(\tfrac{h(-L)}{\eps})-\lambda_1(h(0)-h(-L)),\\[.5ex]
\int_{-L}^0h_1'h''\,dx&=&\phi(\tfrac{h_1(0)}{\eps})-\phi(\tfrac{h_1(-L)}{\eps})-\lambda_2(h_1(0)-h_1(-L)).
\eeas 
Adding these two equalities yields
\be
\hspace{-1cm}\lambda_1[h(0)-h(-L)]+\lambda_2[h_1(0)-h_1(-L)]=\phi(\tfrac{h(0)}{\eps})+\phi(\tfrac{h_1(0)}{\eps})-\phi(\tfrac{h(-L)}{\eps})-\phi(\tfrac{h_1(-L)}{\eps}).
\lb{lambda_rel}
\ee
When expanded in powers of $\eps$, \rf{lambda_rel} will provide us a relation between the leading orders of $\lambda_1$ and $\lambda_2$ for several solutions to system \rf{SSa}--\rf{SSc}.

In this article, we address the following two questions about solutions to the dynamical system \rf{BS}-\rf{BC} and its stationary counterpart \rf{SSa}--\rf{SSc}. Firstly, we are interested in understanding how the number and shape of solutions to \rf{SSa}--\rf{SSc} depends on the choice of physical and geometric parameters $L,\,\sigma,\,|\phi(1)|,\,M_1,\,M$. For a sake of convenience, instead of masses $M_1,\,M$ from \rf{MC1}--\rf{MC} we prefer to parameterize solution to \rf{SSa}--\rf{SSc} using parameters $h_1^m$ and $h^m$ related to $\displaystyle\max_{x\in(-L,\,0)}h_1(x)$ and $\displaystyle\max_{x\in(-L,\,0)}h(x)$, respectively. A similar height parameter was effectively employed by asymptotic matching of the lens stationary solution in section 4.1 of~\cite{JHKPW13}. Accordingly,  in sections $3-6$ for broad ranges of positive $L,\,\sigma,\,|\phi(1)|,\,h_1^m,\,h^m$ (called below the {\it model parameters}) we derive the leading order formulae (in the asymptotic limit $\eps\go 0$) for the corresponding pressures $(\lambda_1,\,\lambda_2)$ and spatial profiles $h(x)$ and $h_1(x)$ of several non-constant solutions to system \rf{SSa}-\rf{SSb}. Following~\cite{GW03,KW10,Ki14,JHKPW13} our approach relies on systematic asymptotic matching of the repeated elementary solution blocks (bulk, contact line and ultra-thin film ones) described in section $2$.

Correspondingly, we distinguish and describe eleven types of the so called {\it composite solutions} to \rf{SSa}--\rf{SSc} exhibiting from one up to four {\it triple phase contact lines}. Note that we do not address constant solutions to \rf{BS}-\rf{BC} or \rf{SSa}--\rf{SSc},  having $\lambda_1=\Pi_\eps(h),\,\lambda_2=\Pi_\eps(h_1)$, and their connection to those found in sections $3-6$  by analyzing possible bifurcation paths. The latter analysis could be similar to the one done in~\cite{BGW01,PBMT05} and is out of the scope this article. Still as a preparation step for it, in section $7$ we plot and analyze the combined diagrams showing for broad ranges of the {\it model parameters} typical shapes of the existence domains (EDs) for all found non-constant solutions to \rf{SSa}--\rf{SSc} together.

Secondly, our numerical simulations of bilayer system \rf{BS}-\rf{BC} show that the eleven stationary solutions (found in sections $3-6$ ) are either dynamically stable or weakly translationally unstable when tested as initial conditions \rf{IC}.  By numerical solving of system \rf{BS}-\rf{BC} we extended and used the fully implicit finite differences scheme with adaptive time step developed in~\cite{Pe08,KW10,KFE18} for other coupled thin film models.  Moreover, in section $8$ again using numerical simulations of system \rf{BS}-\rf{BC} we explain why all other composite stationary solutions to \rf{BS}-\rf{BC} are dynamically unstable and rather experience complex coarsening pathways (nevertheless converging in the long time to one of the stable states found in sections $3-6$). We conclude the article by discussing and comparing our results with the existing modeling and experimental ones and give an outlook in section 9.

%%%%%%%%%%%%%%%%%%%%%%%%%%%%%%%%%%%%%%%%%%%%%%%%%%%%%%%%%%%%%%%%%
\section{Bulk, contact line and UTF profiles}
%%%%%%%%%%%%%%%%%%%%%%%%%%%%%%%%%%%%%%%%%%%%%%%%%%%%%%%%%%%%%%%%%
Our aim is to construct different solutions to stationary system \rf{SSa}--\rf{SSc} by applying systematic asymptotic matching of the leading order expressions (as parameter $\eps\go0$ in \rf{phi}) of typical {\it bulk (B), contact line (CL) and ultra-thin film (UTF)} profiles described in this section. For that, we assume formal (uniform in $x$ variable) asymptotic expansions
\be
\hspace{-1.cm}h_1(x)=h_1^0(x)+\eps h_1^1(x)+\eps^2h_1^2(x)+O(\eps^3),\
h(x)=h^0(x)+\eps h^1(x)+\eps^2h^2(x)+O(\eps^3),
\lb{AE}
\ee  
for solutions to \rf{SSa}--\rf{SSc} and, similarly, expand
\bes
\lambda_i=\lambda_i^0+\eps\lambda_i^1+O(\eps^2)\quad\text{for}\ i=1,\,2.
\ees
The bulk regions have $O(1)$ length w.r.t $\eps$ and either $h_1^0\not=0$ or  $h^0\not=0$ holds there. In CL regions, either $h_1$ or $h$ experiences a transition on a small interval of length $O(\eps)$ around a point $x=-s$  from the magnitude of $O(1)$ to that one of UTF, i.e. to the leading order constant $\eps$ profile. Following the standard approach of~\cite{GW03,GW05,KW10,Ki14,JHKPW13} we use the change of variable
\be
z=\frac{x+s}{\eps},
\lb{CL_CV}
\ee
under which equations \rf{SSa}--\rf{SSb} in CL regions transform to equations
\sbea
\lb{SS_CLa}
\frac{\sigma}{\eps^2}h_1''(z)&=&\Pi_\eps(h_1(z))-\Pi_\eps(h(z))-\lambda_2+\lambda_1,\\[1ex]
\lb{SS_CLb}
\frac{\sigma}{\eps^2}h''(z)&=&-\Pi_\eps(h_1(z))+(\sigma+1)\Pi_\eps(h(z))+\lambda_2-(\sigma+1)\lambda_1,
\seea
considered for $z\in(-\infty,\,+\infty)$.

Below we distinguish three bulk, four CL and one UTF types of regions and describe the corresponding leading order profiles of $h_1$ and $h$ in them.

\vspace*{.5cm}
\underline{\bf Type I bulk profile:}
 both $h_1(x)$ and $h(x)$ are assumed to be $O(1)$ functions as $\eps\go 0$. Correspondingly, the leading order  equations \rf{SSa}--\rf{SSb} take the form:
\bes
\sigma (h_1^0)''=-\lambda_2^0+\lambda_1^0,\quad
\sigma (h^0)''=\lambda_2^0-(\sigma+1)\lambda_1^0,
\ees
Integrating them twice in $x$ without applying boundary conditions {\it the general type I bulk} profiles are obtained:
\be
\lb{BP1}
h_1^0(x)=\tfrac{\lambda_1^0-\lambda_2^0}{2\sigma}(x-x_{c1})^2+C_1,\quad
h^0(x)=\tfrac{\lambda_2^0-(\sigma+1)\lambda_1^0}{2\sigma}(x-x_c)^2+C.
\ee
In \rf{BP1}, $x_c,\,x_{c1}$ and $C,\,C_1$ are four integration constants, the former giving the centers and the latter the maximum (or minimum) values of parabolic profiles $h^0(x)$ and $h_1^0(x)$, respectively. 

\vspace*{.5cm}
\underline{\bf Type II bulk profile:}
 we assume $h_1(x)=O(1)$ and $h(x)=O(\eps)$ as $\eps\go 0$, i.e. in asymptotic expansions \rf{AE} we set $h^0(x)=0$. To balance the potential terms $\Pi_\eps(h)$ in \rf{SSa}--\rf{SSb} we additionally set $h^1(x)=1$. Next,
using the following expansion of \rf{Pi}:
\be
\Pi_\eps(\eps+\eps^2h)=(l-n)h+O(\eps),
\lb{Pi_ex}
\ee
the leading order of \rf{SSa}--\rf{SSb} becomes
\bes
\sigma (h_1^0)''=-(l-n)h^2-\lambda_2^0+\lambda_1^0,\quad
0=(\sigma+1)(l-n)h^2+\lambda_2^0-(\sigma+1)\lambda_1^0.
\ees
Expressing from the second equation
$$h^2=-\tfrac{\lambda_2^0-(\sigma+1)\lambda_1^0}{(\sigma+1)(l-n)},$$
substituting into the first one, and integrating the latter twice in $x$ yields {\it the general type II bulk} profiles:
\be
\lb{BP2}
h_1^0(x)=-\tfrac{\lambda_2^0}{2(\sigma+1)}(x-\tilde{x}_{c1})^2+\widetilde{C}_1,\quad
h^1(x)=1,
\ee
with $\tilde{x}_{c1}$ and $\widetilde{C}_1$ being two integration constants giving the center and the maximum values, respectively of concave parabolic profile $h_1^0(x)$. 

\vspace*{.5cm}
\underline{\bf Type III bulk profile:}
 we assume $h_1(x)=O(\eps)$ and $h(x)=O(1)$ as $\eps\go 0$, i.e. in \rf{AE} one has $h_1^0(x)=0$ and $h_1^1(x)=1$. Using again expansion \rf{Pi_ex} the leading order of \rf{SSa}--\rf{SSb} becomes
\bes
0=(l-n)h_1^2-\lambda_2^0+\lambda_1^0,\quad
\sigma (h^0)''=-(l-n)h_1^2+\lambda_2^0-(\sigma+1)\lambda_1^0.
\ees
Expressing from the first equation
$$h_1^2=-\tfrac{\lambda_1^0-\lambda_2^0}{l-n},$$
substituting into the second one, and integrating the latter twice in $x$ yields {\it the general type III bulk} profiles:
\be
\lb{BP3}
h_1^1(x)=1,\quad
h^0(x)=-\tfrac{\lambda_1^0}{2}(x-\tilde{x}_{c})^2+\widetilde{C},
\ee
with $\tilde{x}_{c}$ and $\widetilde{C}$ being the integration constants giving the center and the maximum values, respectively of concave parabolic profile $h^0(x)$. 

\vspace*{.5cm}
\underline{\bf UTF profile:}
 we assume $h_1^0(x)=h^0(x)=0$ and $h_1^1(x)=h^1(x)=1$ in \rf{AE}. Using again \rf{Pi_ex}, the leading order of \rf{SSa}--\rf{SSb} becomes
\beas
0&=&-(l-n)h^2+(l-n)h_1^2-\lambda_2^0+\lambda_1^0,\\
0&=&(\sigma+1)(l-n)h^2-(l-n)h_1^2+\lambda_2^0-(\sigma+1)\lambda_1^0,
\eeas
implying the second order constant profiles 
$$h_1^2(x)=\tfrac{\lambda_2^0}{l-n}\quad\text{and}\ h^2(x)=\tfrac{\lambda_1^0}{l-n}.$$
\underline{\bf Type I CL profile:}
 we assume $h_1^0=0$ and $h^0\not=0$  in \rf{AE}, i.e. $h_1(x)$ experiences a transition from $O(1)$ to the   leading order constant $\eps$ magnitude, while $h(x)$ stays $O(1)$. Applying the variable change \rf{CL_CV}, the leading order of the scaled system \rf{SS_CLa}--\rf{SS_CLb} becomes
\sbea
\lb{CL1a_lo}
\sigma(h_1^1)''(z)&=&\phi'(h_1^1(z)),\\
\sigma (h^0)''(z)&=&0
\lb{CL1b_lo}
\seea
with potential $\phi$ from \rf{phi}. In turn, the first order corrector equation for \rf{SS_CLb} takes the form $\sigma (h^1)''=-\phi'(h_1^1)=-\sigma(h_1^1)''$. Next, integrating \rf{CL1a_lo} once and \rf{CL1b_lo} twice  in $z$ one obtains
\sbea
\lb{CL1a}
|(h_1^1)'(z)|&=&\sqrt{\tfrac{2}{\sigma}[\phi(h_1^1(z))-\phi(1)]},\\[1ex]
\lb{CL1b}
h^0(z)&=&C_2,\quad h^1(z)=-h_1^1(z)+C_3z+C_*
\seea
with $C_2,\,C_3,\,C_*$ being the integration constants. The sign of  $(h_1^1)'(z)$ in \rf{CL1a} is positive (negative)
if $h_1^1(z)$ is increasing (decreasing).  

\vspace*{.5cm}
\underline{\bf Type II CL profile:}
 we assume $h^0=0$ and $h_1^0\not=0$  in \rf{AE}, i.e. $h(x)$ experiences a transition from $O(1)$ to the leading order constant $\eps$ magnitude, while $h_1(x)$ stays $O(1)$. Applying again \rf{CL_CV}, the leading order of  \rf{SS_CLa}--\rf{SS_CLb} becomes
\sbea
\lb{CL2a_lo}
\sigma(h_1^0)''(z)&=&0,\\
\sigma (h^1)''(z)&=&(\sigma+1)\phi'(h^1(z)).
\lb{CL2b_lo}
\seea
In turn, the first order corrector equation for \rf{SS_CLa} takes the form: $$\sigma (h_1^1)''=-\phi'(h^1)=-\tfrac{\sigma}{\sigma+1}(h^1)''.$$ Next, integrating once \rf{CL2b_lo} and \rf{CL2a_lo} twice  in $z$  one obtains
\sbea
\lb{CL2a}
h_1^0(z)&=&C_4,\quad h_1^1(z)=-\tfrac{h^1(z)}{\sigma+1}+C_5z+C_{**},\\[1ex]
\lb{CL2b}
|(h^1)'(z)|&=&\sqrt{\tfrac{2(\sigma+1)}{\sigma}[\phi(h^1(z))-\phi(1)]}
\seea
with $C_4,\,C_5,\,C_{**}$ being the integration constants.

\vspace*{.5cm}
\underline{\bf Type III CL profile:}
 we assume $h_1^0=h^0=0,\,h^1(x)=1$  in \rf{AE} with $h_1(x)$ experiencing a transition from $O(1)$ to the leading order  constant $\eps$ magnitude, while $h(x)=O(\eps)$. Applying \rf{CL_CV} the leading order of \rf{SS_CLa}--\rf{SS_CLb} becomes
\bes
\sigma(h_1^1)''(z)=\phi'(h_1^1(z))-\phi'(h^1(z)),\quad
\sigma (h^1)''(z)=-\phi'(h_1^1(z))+(\sigma+1)\phi'(h^1(z)).
\ees
This system can be rewritten as
\sbea
\lb{CL_s1}
(h_1^1+h^1)''(z)&=&-\phi'(h^1),\\
((\sigma+1)h_1^1+h^1)''(z)&=&\phi'(h_1^1).
\lb{CL_s2}
\seea
Multiplying \rf{CL_s1} and \rf{CL_s2} by $(h^1)'(z)$ and $(h_1^1)'(z)$, respectively, and integrating them for $z\in(-\infty,\,\infty)$ by parts gives
\bes
-\int_{-\infty}^{+\infty}(h_1^1)'(h^1)''\,dz=0,\quad
\tfrac{\sigma+1}{2}|(h_1^1)'|^2\Big|_{-\infty}^{+\infty}+\int_{-\infty}^{+\infty}(h_1^1)'(h^1)''\,dz=\phi(h_1^1(z))\Big|_{-\infty}^{+\infty},
\ees 
where we used $h^1(-\infty)=h^1(+\infty)=1$ and $(h^1)'(-\infty)=(h^1)'(+\infty)=0$. These  imply together
\bes
\tfrac{\sigma+1}{2}|(h_1^1)'|^2\Big|_{-\infty}^{+\infty}=\phi(h_1^1(z))\Big|_{-\infty}^{+\infty},
\ees
and, consequently, switching back to $x$ variable in \rf{CL_CV}, {\it the macroscopic contact angle}:
\be
|(h_1^1)'(-s)|=\sqrt{\tfrac{2}{\sigma+1}|\phi(1)|}
\lb{CA_h1}
\ee
is obtained. In the last step, we used $\phi(+\infty)=0$.

\vspace*{.5cm}
\underline{\bf Type IV CL profile:}
we assume $h_1^0=h^0=0,\,h_1^1(x)=1$ in \rf{AE} with $h(x)$ experiencing a transition from $O(1)$ to the leading order constant $\eps$ magnitude, while $h_1(x)=O(\eps)$. Applying \rf{CL_CV}, the leading order of \rf{SS_CLa}--\rf{SS_CLb} can be again written in the form \rf{CL_s1}--\rf{CL_s2}. Integrating it for $z\in(-\infty,\,\infty)$ leads now to
\bes
\tfrac{1}{2}|(h^1)'|^2\Big|_{-\infty}^{+\infty}-\int_{-\infty}^{+\infty}(h_1^1)'(h^1)''\,dz=-\phi(h^1(z))\Big|_{-\infty}^{+\infty},\quad
\int_{-\infty}^{+\infty}(h_1^1)'(h^1)''\,dz=0,
\ees 
where we used $h_1^1(-\infty)=h_1^1(+\infty)=1$ and $(h_1^1)'(-\infty)=(h_1^1)'(+\infty)=0$. These imply together
\bes
\tfrac{1}{2}|(h^1)'|^2\Big|_{-\infty}^{+\infty}=-\phi(h^1(z))\Big|_{-\infty}^{+\infty},
\ees
and, consequently, {\it the macroscopic contact angle}:  
\be
|(h^1)'(-s)|=\sqrt{2|\phi(1)|}.
\lb{CA_h}
\ee

%%%%%%%%%%%%%%%%%%%%%%%%%%%%%%%%%%%%%%%%%%%%%%%%%%%%%%%%%%%%%%%%%
\section{One-CL solutions}
%%%%%%%%%%%%%%%%%%%%%%%%%%%%%%%%%%%%%%%%%%%%%%%%%%%%%%%%%%%%%%%%%
In this section, we classify and derive the leading order  profiles of the solutions to system \rf{SSa}--\rf{SSc} having one contact line. Up to a possible inversion of $x$-variable we distinguish four types of such solutions: {\it lens}, {\it internal drop}, {\it$h_1$-drop}, and {\it$h$-drop} ones.\\[.5ex]

\underline{\bf Lens:}\hspace{.5cm}The solution is obtained by matching {\bf Type I} (with $x_c=x_{c1}=0$) and {\bf Type II} (with $\tilde{x}_{c1}=-L$) bulk solutions to {\bf Type II} CL ones centered around $x=-s$. For that we rewrite the bulk solutions \rf{BP1} and \rf{BP2} as functions of the inner CL variable $z$ using \rf{CL_CV}, expand them in powers of $\eps$, and subsequently match the leading orders of them and their first derivatives (w.r.t. $x$) to the corresponding expansions of CL solutions \rf{CL2a}--\rf{CL2b}. This procedure leads to the following set of matching conditions:
\be
\left\{\begin{array}{l}
	-\frac{\lambda_2^0}{2(\sigma+1)}(-s+L)^2+\widetilde{C}_1=C_4,\quad
	-\frac{\lambda_2^0}{\sigma+1}(-s+L)=C_5,\\[1.5ex]
	\frac{\lambda_1^0-\lambda_2^0}{2\sigma}s^2+C_1=C_4,\quad\hspace{2.cm}
	\frac{\lambda_2^0-(\sigma+1)\lambda_1^0}{2\sigma}s^2+C=0,\\[1.5ex]
	\frac{\lambda_1^0-\lambda_2^0}{\sigma}(-s)=-\sqrt{\frac{2|\phi(1)|}{\sigma(\sigma+1)}}+C_5,\quad
	\frac{\lambda_2^0-(\sigma+1)\lambda_1^0}{\sigma}(-s)=\sqrt{\frac{2|\phi(1)|(\sigma+1)}{\sigma}}.\end{array}\right.
\lb{0sol}
\ee
Additionally, we fix $C=h^m$ and $\widetilde{C}_1=h_1^m$ as the given maximum values of $h(x)$ and $h_1(x)$, respectively. System \rf{0sol} has $6$ equations with $6$ unknowns: $\lambda_1^0,\,\lambda_2^0,\,s$ and $C_1,\,C_4,\,C_5$. From its last two equations we find $\lambda_2^0s=(\sigma+1)C_5$. Substituting that into the $2^{\mathrm{nd}}$ equation in \rf{0sol} gives $\lambda_2^0=C_5=0$. Dividing the $4^{\mathrm{th}}$ over the last equations in \rf{0sol} provides expressions 
\bes
\lambda_1^0=\frac{|\phi(1)|}{h^m}\quad\text{and, subsequently,}\ s=\sqrt{\tfrac{2\sigma}{(\sigma+1)|\phi(1)|}}h^m.
\ees
The derived {\it lens} solution is well defined if the following two constraints are obeyed:
\be
L-s=L-\sqrt{\tfrac{2\sigma}{(\sigma+1)|\phi(1)|}}h^m>0,\quad C_1=h_1^m-\tfrac{h^m}{\sigma+1}>0.
\lb{lens_constr}
\ee 
Note that the second constraint sets a relation between the maximal heights $\bar{h}=h^m/h_1^m<\sigma+1$ for which both $h(x)$ and $h_1(x)$ stay positive for all $x\in(-L,0)$, while the first one controls the minimal interval length $L$ in which the solution can be allocated. Typical shape of  the {\it lens} solution is shown in Fig.2 {\bf(a)}.\\[.5ex]

\underline{\bf Internal drop:}\hspace{.5cm}The solution is obtained by matching {\bf Type I} (with $x_c=x_{c1}=0$) and {\bf Type III} (with $\tilde{x}_c=-L$) bulk solutions to {\bf Type I} CL ones centered around $x=-s$. Again we use \rf{CL_CV} and match \rf{BP1}, \rf{BP3} and their first derivatives (w.r.t. $x$) to the corresponding expansions of \rf{CL1a}--\rf{CL1b}. This procedure leads to the following set of matching conditions:
\be
\left\{\begin{array}{l}
	-\frac{\lambda_1^0}{2}(-s+L)^2+\widetilde{C}=C_2,\quad
	-\lambda_1^0(-s+L)=C_3,\\[1.5ex]
	\frac{\lambda_1^0-\lambda_2^0}{2\sigma}s^2+C_1=0,\hspace{1.8cm}
	\frac{\lambda_2^0-(\sigma+1)\lambda_1^0}{2\sigma}s^2+C=C_2,\\[1.5ex]
	\frac{\lambda_1^0-\lambda_2^0}{\sigma}(-s)=\sqrt{\frac{2|\phi(1)|}{\sigma}},\hspace{1.1cm}
	\frac{\lambda_2^0-(\sigma+1)\lambda_1^0}{\sigma}(-s)=-\sqrt{\frac{2|\phi(1)|}{\sigma}}+C_3.\end{array}\right.
\lb{1sol}
\ee
Additionally, we fix $\widetilde{C}=h^m$ and $C_1=h_1^m$. System \rf{1sol} has $6$ equations with $6$ unknowns: $\lambda_1^0,\,\lambda_2^0,\,s$ and $C,\,C_2,\,C_3$. Solving it proceeds similar to that one of system \rf{0sol} for the {\it lens} solution and yields the expressions
\bes
\lambda_1^0=0,\quad\lambda_2^0=\frac{|\phi(1)|}{h_1^m},\quad\text{and}\ s=\sqrt{\tfrac{2\sigma}{|\phi(1)|}}h_1^m,
\ees
%%%%%%%%%%%%%%%%%%%%%%%%%%%%%%%%%%%%%%%%%%%%%%%%%%%%%%%%%%%%%%%%%%%%%%%%%%%%%%%%%%%%%%%%%%%%%%%%%%%%%%%%%%%%%%%%
\begin{figure}[H] 
	\centering
	\vspace{-2.9cm}
	\hspace{-.4cm}\includegraphics[width=.5\textwidth]{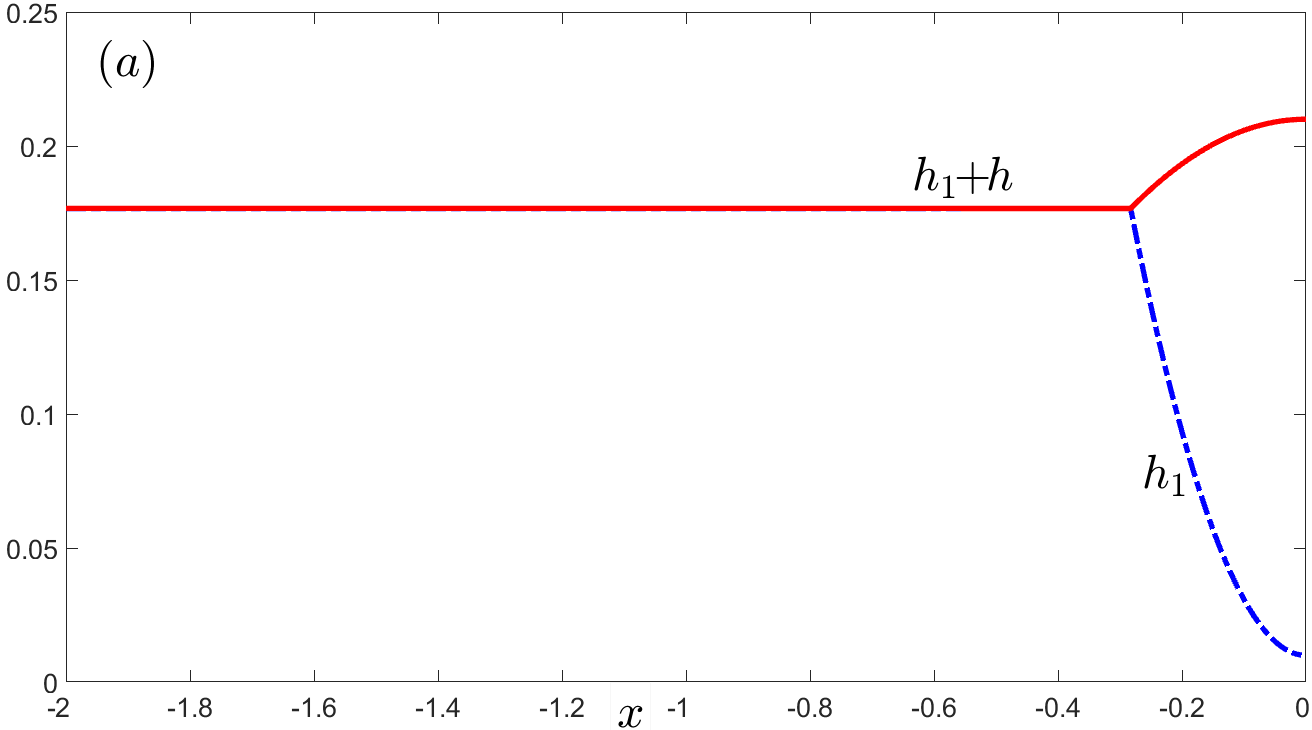}  
	\hspace{.2cm}\includegraphics[width=.5\textwidth]{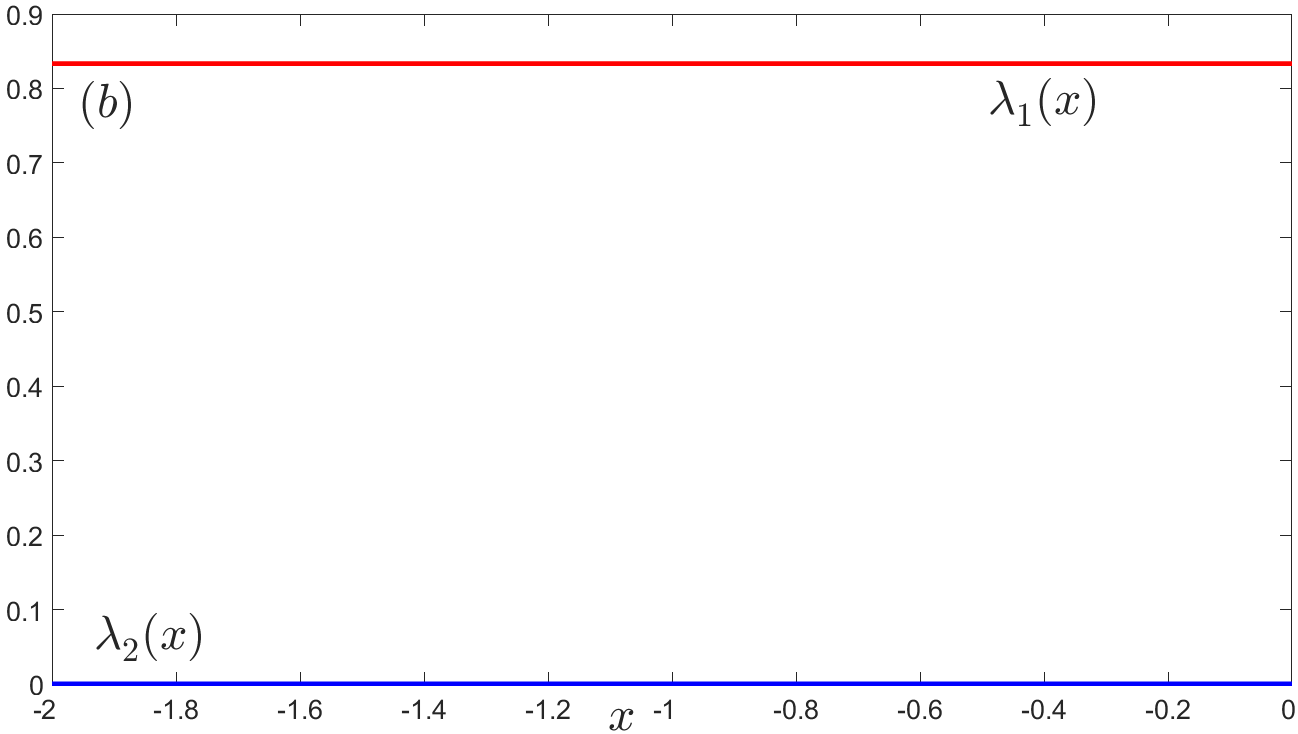}
	\caption{\small Numerical lens stationary solution {\bf(a)} to system \rf{BS}--\rf{BC} for\hspace{3.cm}$\eps\!=\!0.00005$, $\sigma\!=\!0.2,\,L\!=\!2.0$   with constant pressures {\bf(b)} $\lambda_1\!=\!0.833301,\, \lambda_2\!=\!0.0000519$.}	
\end{figure}
%%%%%%%%%%%%%%%%%%%%%%%%%%%%%%%%%%%%%%%%%%%%%%%%%%%%%%%%%%%%%%%%%%%%%%%%%%%%%%%%%%%%%%%%%%%%%%%%%%%%%%%%%%%%%%%%
\noindent as well as two constraints:
\be
L-s=L-\sqrt{\tfrac{2\sigma}{|\phi(1)|}}h_1^m>0\quad\text{and}\ C=h^m-h_1^m>0.
\lb{intdr_constr}
\ee 
Additionally, we note that $h_1^0(x)+h^0(x)=h^m$ holds for all $x\in(-L,\,0)$, i.e. the second fluid layer is {\it uniformly flat}, while the first one has the form of a parabolic drop. Typical shape of the  {\it internal drop} solution is shown in Fig.3  {\bf(a)}.
%%%%%%%%%%%%%%%%%%%%%%%%%%%%%%%%%%%%%%%%%%%%%%%%%%%%%%%%%%%%%%%%%%%%%%%%%%%%%%%%%%%%%%%%%%%%%%%%%%%%%%%%%%%%%%%%
\begin{figure}[H] 
	\centering
	\vspace{-.4cm}
	\hspace{-.4cm}\includegraphics[width=.5\textwidth]{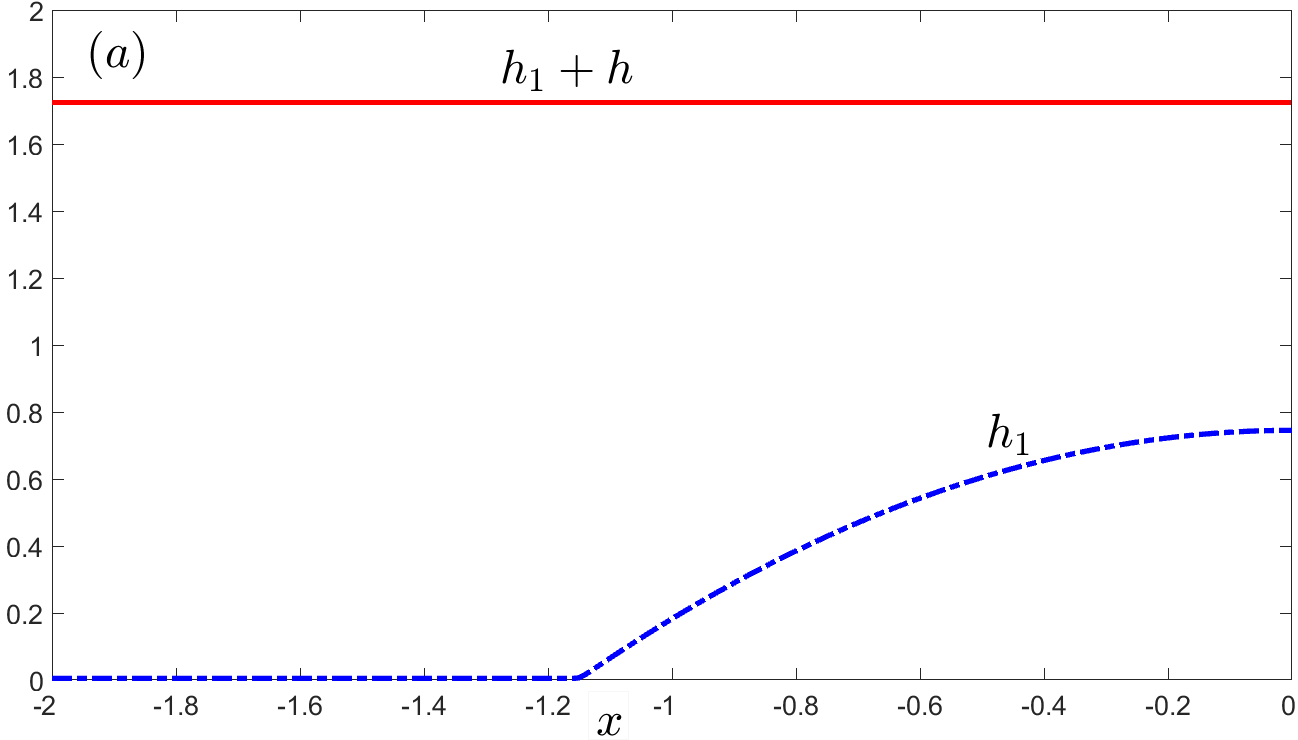} 
	\hspace{.2cm}\includegraphics[width=.5\textwidth]{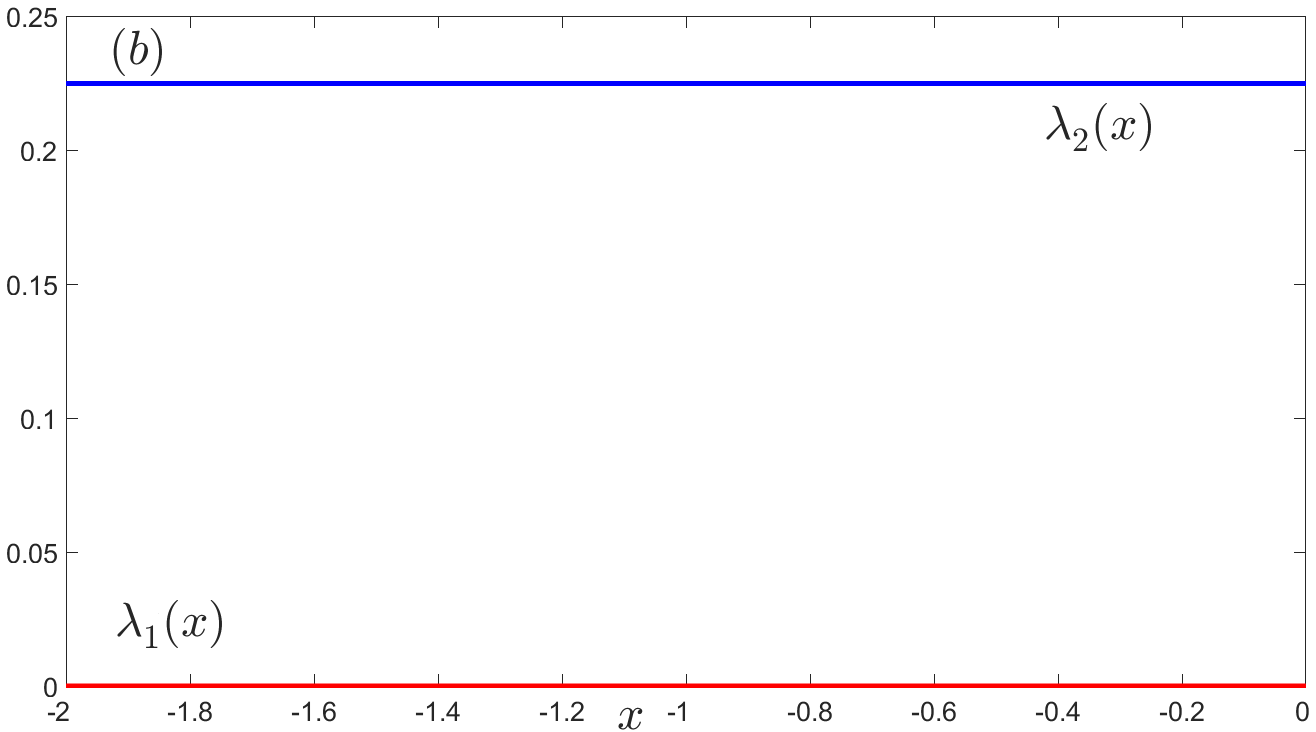}  
	\caption{\small Numerical internal drop stationary solution {\bf(a)} to system \rf{BS}--\rf{BC} for $\eps\!=\!0.005,\,\sigma\!=\!0.2,\,L\!=\!2.0$ with constant pressures {\bf(b)} $\lambda_1\!=\!0.0000113,\, \lambda_2\!=\!0.22512$.}	
\end{figure}
%%%%%%%%%%%%%%%%%%%%%%%%%%%%%%%%%%%%%%%%%%%%%%%%%%%%%%%%%%%%%%%%%%%%%%%%%%%%%%%%%%%%%%%%%%%%%%%%%%%%%%%%%%%%%%%%
\noindent\underline{\bf $h_1$-drop:}\hspace{.5cm}The solution is obtained by matching {\bf Type II} (with $\tilde{x}_{c1}=-L$) bulk and {\bf UTF} solutions to {\bf Type III} CL ones centered around $x=-s$. We match \rf{BP2} using contact angle condition \rf{CA_h1} to get
\be
-\tfrac{\lambda_2^0}{2(\sigma+1)}(-s+L)^2+h_1^m=0,\quad-\tfrac{\lambda_2^0}{\sigma+1}(-s+L)=-\sqrt{\tfrac{2|\phi(1)|}{\sigma+1}},
\lb{h1sol}
\ee
implying that
\bes
\lambda_2^0=\frac{|\phi(1)|}{h_1^m}\quad\text{and}\ s=L-\sqrt{\tfrac{2(\sigma+1)}{|\phi(1)|}}h_1^m.
\ees
The only constraint imposed on this solution is $s>0$, i.e. the one controlling the minimal length $L$:
\be
L>\sqrt{\tfrac{2(\sigma+1)}{|\phi(1)|}}h_1^m.
\lb{h1dr_constr}
\ee

Note, that system \rf{h1sol} does not determine the value of $\lambda_1^0$, which is observed to be nonzero and negative in the numerically found {\it $h_1$-drop} solutions, a typical shape of which  is shown in Fig.4  {\bf(a)}.\\[.5ex]
%%%%%%%%%%%%%%%%%%%%%%%%%%%%%%%%%%%%%%%%%%%%%%%%%%%%%%%%%%%%%%%%%%%%%%%%%%%%%%%%%%%%%%%%%%%%%%%%%%%%%%%%%%%%%%%%
\begin{figure}[H] 
	\centering
	\vspace{-2.9cm}
	\hspace{-.4cm}\includegraphics[width=.5\textwidth]{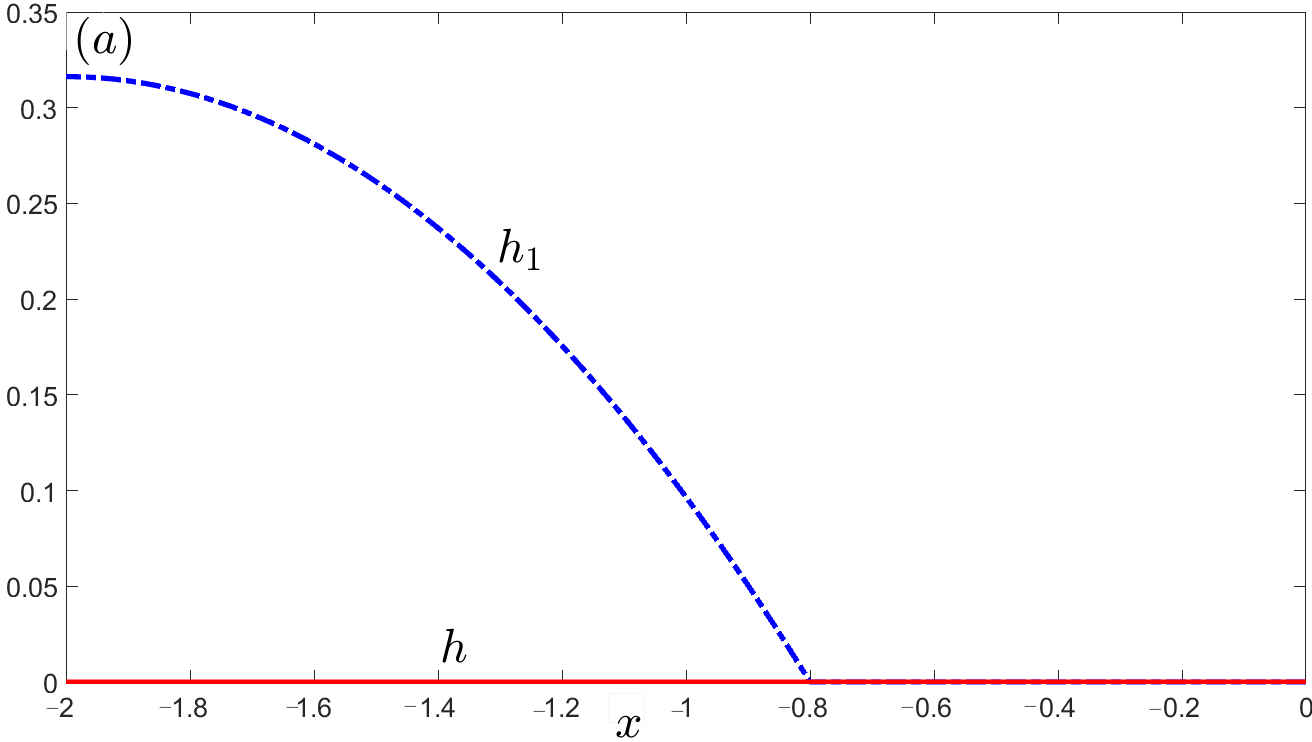} 
	\hspace{.2cm}\includegraphics[width=.5\textwidth]{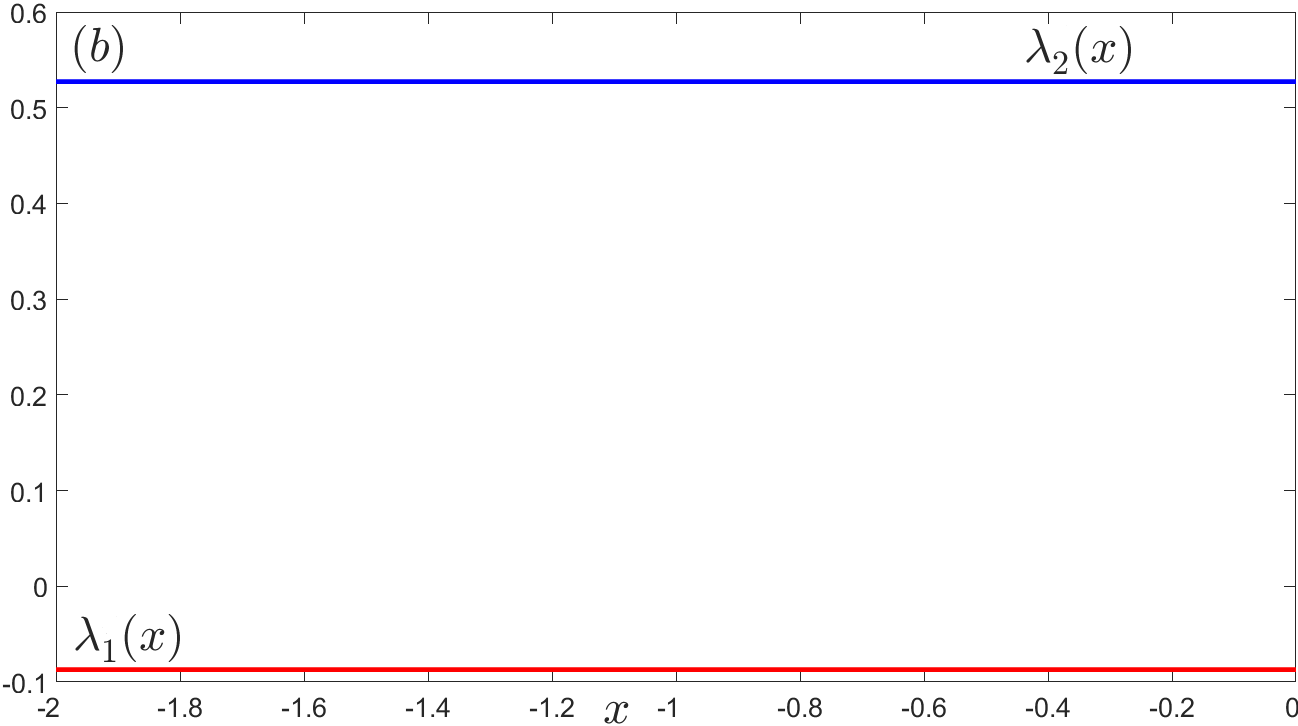} 
	\caption{\small Numerical $h_1$-drop stationary solution {\bf(a)} to system \rf{BS}--\rf{BC} for $\eps\!=\!0.00005,\,\sigma\!=\!0.2,\,L\!=\!2.0$ with constant pressures {\bf(b)} $\lambda_1\!=\!-0.087101$, $\lambda_2\!=\!0.527111$.}	
\end{figure}
%%%%%%%%%%%%%%%%%%%%%%%%%%%%%%%%%%%%%%%%%%%%%%%%%%%%%%%%%%%%%%%%%%%%%%%%%%%%%%%%%%%%%%%%%%%%%%%%%%%%%%%%%%%%%%%%
\underline{\bf $h$-drop:}\hspace{.5cm}The solution is obtained by matching {\bf Type III} (with $\tilde{x}_{c}=-L$) bulk and {\bf UTF} solutions to {\bf Type IV} CL ones centered around $x=-s$. We match \rf{BP3} using contact angle condition \rf{CA_h} to get
\be
-\tfrac{\lambda_1^0}{2}(-s+L)^2+h^m=0,\quad-\lambda_1^0(-s+L)=-\sqrt{2|\phi(1)|},
\lb{hsol}
\ee
implying that
\bes
\lambda_1^0=\frac{|\phi(1)|}{h^m}\quad\text{and}\ s=L-\sqrt{\tfrac{2}{|\phi(1)|}}h^m.
\ees
%%%%%%%%%%%%%%%%%%%%%%%%%%%%%%%%%%%%%%%%%%%%%%%%%%%%%%%%%%%%%%%%%%%%%%%%%%%%%%%%%%%%%%%%%%%%%%%%%%%%%%%%%%%%%%%%
\begin{figure}[H] 
	\centering
	\vspace{-.45cm}
	\hspace{-.4cm}\includegraphics[width=.5\textwidth]{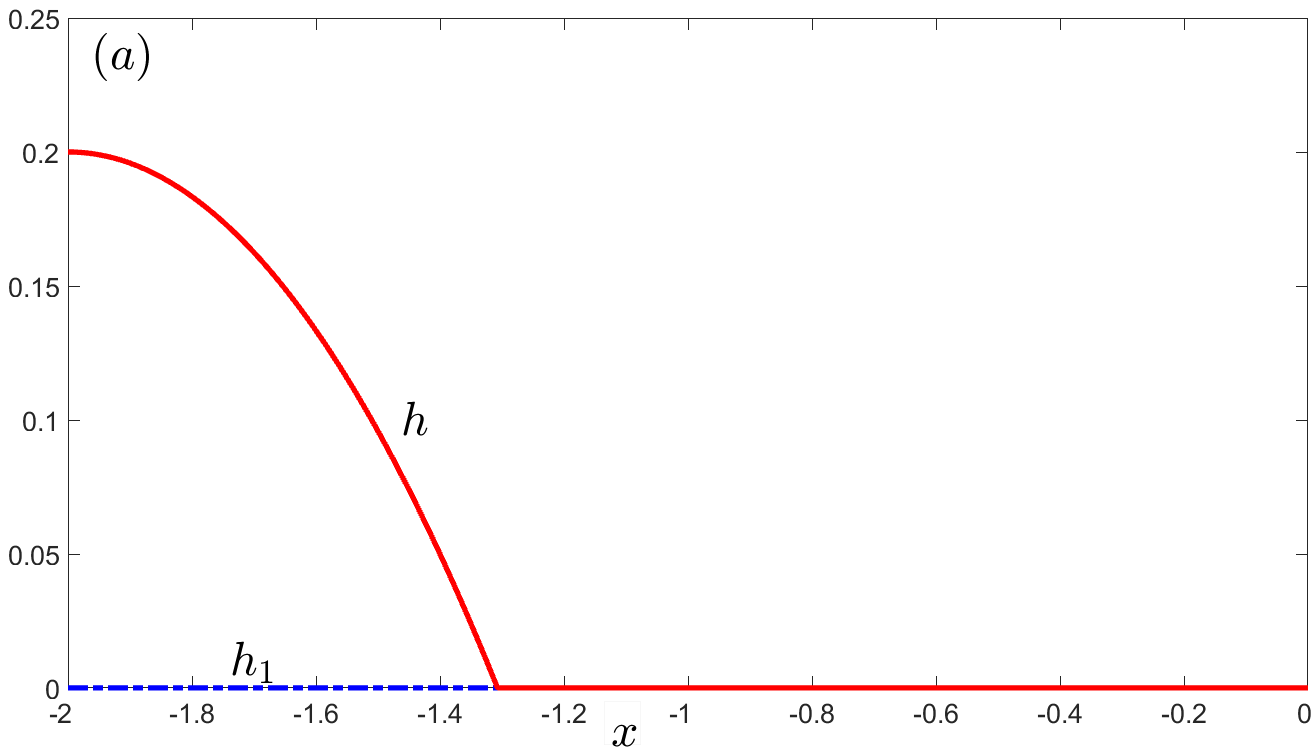}  
	\hspace{.2cm}\includegraphics[width=.5\textwidth]{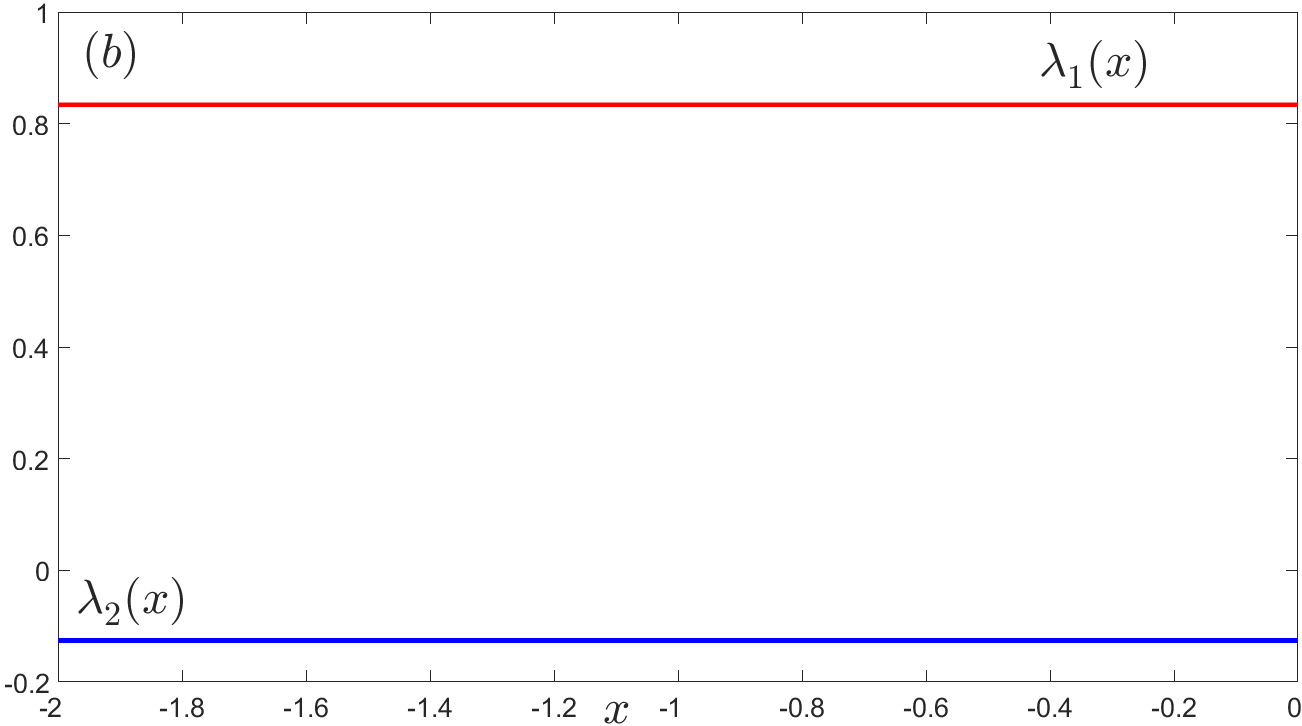} 
	\caption{\small Numerical $h$-drop stationary solution  {\bf(a)} to system \rf{BS}--\rf{BC} for $\eps\!=\!0.00005,\,\sigma\!=\!0.2,\,L\!=\!2.0$ with constant pressures {\bf(b)} $\lambda_1\!=\!0.833511$, $\lambda_2\!=\!-0.125582$.}	
\end{figure}
%%%%%%%%%%%%%%%%%%%%%%%%%%%%%%%%%%%%%%%%%%%%%%%%%%%%%%%%%%%%%%%%%%%%%%%%%%%%%%%%%%%%%%%%%%%%%%%%%%%%%%%%%%%%%%%%
\noindent The only constraint imposed on this solution is again $s>0$, i.e.
\be
L>\sqrt{\tfrac{2}{|\phi(1)|}}h^m.
\lb{hdr_constr}
\ee

Again system \rf{hsol} does not determine the value of $\lambda_2^0$, which is observed to be nonzero and negative in the  numerically found {\it $h$-drop} solutions, a typical shape of which is shown in Fig.5 {\bf(a)}.

%%%%%%%%%%%%%%%%%%%%%%%%%%%%%%%%%%%%%%%%%%%%%%%%%%%%%%%%%%%%%%%%%
\section{Two-CL solutions}
%%%%%%%%%%%%%%%%%%%%%%%%%%%%%%%%%%%%%%%%%%%%%%%%%%%%%%%%%%%%%%%%%
In this section, we classify and derive the leading order profiles of the solutions to system \rf{SSa}--\rf{SSc} having two contact lines. Up to a possible inversion of $x$-variable we distinguish four types of them: {\it zig-zag}, {\it sessile lens}, {\it sessile internal drop}, and {\it $2$-drops} ones.\\[.5ex]

\underline{\bf Zig-zag:}\hspace{.5cm}The solution is obtained from the combined matching of general {\bf Type I} and {\bf Type III} (with $\tilde{x}_{c}=-L$) bulk solutions to {\bf Type I} CL ones centered around $x=-s_{1}$ ($1^{\mathrm{st}}$ CL) as well as of general {\bf Type I} and {\bf Type II} (with $\tilde{x}_{c1}=0$) bulk solutions to {\bf Type II} CL ones centered around $x=-s$ ($2^{\mathrm{nd}}$ CL). Schematically this matching chain can be encrypted as $\bf{IIIB-ICL-IB-IICL-IIB}$ once moving in $x$ from $-L$ to $0$.

At the contact line points $x=-s_1$ and $x=-s$ we apply analogous procedures to those used in section 3 for the one-CL solutions involving the matching of them and their first derivatives. For shortness, we omit the details and just state the resulting system of the leading order matching conditions:
\be
\left\{\begin{array}{l}
	-\frac{\lambda_1^0}{2}(-s_1+L)^2+\widetilde{C}=C_2,\hspace{1.8cm}-\lambda_1^0(-s_1+L)=C_3,\\[1.5ex]
	\frac{\lambda_1^0-\lambda_2^0}{2\sigma}(s_1+x_{c1})^2+C_1=0,\hspace{1.8cm}
	\frac{\lambda_2^0-(\sigma+1)\lambda_1^0}{2\sigma}(s_1+x_c)^2+C=C_2,\\[1.5ex]
	\frac{\lambda_1^0-\lambda_2^0}{\sigma}(-s_1-x_{c1})=\sqrt{\frac{2|\phi(1)|}{\sigma}},\hspace{1.5cm}
	\frac{\lambda_2^0-(\sigma+1)\lambda_1^0}{\sigma}(-s_1-x_c)=-\sqrt{\frac{2|\phi(1)|}{\sigma}}+C_3,\\[1.5ex]
	\frac{\lambda_1^0-\lambda_2^0}{2\sigma}(s+x_{c1})^2+C_1=C_4,\hspace{1.6cm}
	\frac{\lambda_2^0-(\sigma+1)\lambda_1^0}{2\sigma}(s+x_c)^2+C=0,\\[1.5ex]
	\frac{\lambda_1^0-\lambda_2^0}{\sigma}(-s-x_{c1})=\sqrt{\frac{2|\phi(1)|}{\sigma(\sigma+1)}}+C_5,\quad
	\frac{\lambda_2^0-(\sigma+1)\lambda_1^0}{\sigma}(-s-x_c)=-\sqrt{\frac{2(\sigma+1)|\phi(1)|}{\sigma}},\\[1.5ex]
	-\frac{\lambda_2^0}{2(\sigma+1)}s^2+\widetilde{C}_1=C_4,\hspace{2.6cm}
	-\frac{\lambda_2^0}{\sigma+1}(-s)=C_5.\end{array}\right.
\lb{1m0sol}
\ee
Additionally, we fix $\widetilde{C}=h^m$ and $\widetilde{C}_1=h_1^m$. System \rf{1m0sol} has $12$ equations with $12$ unknowns: $\lambda_1^0,\,\lambda_2^0,\,s,\,s_1,\,x_c,\,x_{c1}$ and $C,\,C_1,\,C_2,\,C_3,\,C_4,\,C_5$.
From the linear part of this system, i.e. by performing joint manipulations with the 2$^{\mathrm{nd}}$, 5-6$^{\mathrm{th}}$, 9-10$^{\mathrm{th}}$, and 12$^{\mathrm{th}}$ equations in \rf{1m0sol}, the following expressions for positions are obtained:
\bea
s&=&\tfrac{\sqrt{2\sigma(\sigma+1)|\phi(1)|}}{\lambda_2^0-(\sigma+1)\lambda_1^0}-x_c,\quad s_1=\tfrac{\sqrt{2\sigma|\phi(1)|}}{\lambda_2^0-\lambda_1^0}-x_{c1},\nonumber\\[1.ex]
x_c&=&\tfrac{\lambda_1^0L(\sigma+1)}{\lambda_2^0-(\sigma+1)\lambda_1^0},\quad\quad\quad\quad\quad x_{c1}=\tfrac{\lambda_1^0L}{\lambda_2^0-\lambda_1^0}.
\lb{sx_expr}
\eea
In turn, the rest of the equations of system \rf{1m0sol} imply two quadratic relations
\sbea
\lb{QP1}
\hspace{-1cm}\tfrac{\lambda_2^0-(\sigma+1)\lambda_1^0}{2\sigma}\left[(s_1+x_c)^2-(s+x_c)^2\right]&=&-\tfrac{\lambda_1^0}{2}(-s_1+L)^2+h^m,\\[1.ex]
\tfrac{\lambda_1^0-\lambda_2^0}{2\sigma}\left[(s+x_{c1})^2-(s_1+x_{c1})^2\right]&=&-\tfrac{\lambda_2^0}{2(\sigma+1)}s^2+h_1^m,
\lb{QP2}
\seea
from which the expressions for $\lambda_1^0$ and $\lambda_2^0$ in terms of the {\it model parameters}\\ $L,\,\sigma,\,|\phi(1)|,\,h_1^m,\,h^m$ can be derived as follows.

First, we substitute expressions \rf{sx_expr} into \rf{QP1} and divide both sides of the latter by $(\lambda_2^0-(\sigma+1)\lambda_1^0)/(2\sigma)$ to get

\bes
\hspace{-2.cm}\left[\tfrac{\sqrt{2\sigma|\phi(1)|}}{\lambda_2^0-\lambda_1^0}+\tfrac{\sigma L\lambda_1^0\lambda_2^0}{(\lambda_2^0-(\sigma+1)\lambda_1^0)(\lambda_2^0-\lambda_1^0)}\right]^2-
\left[\tfrac{L\lambda_2^0-\sqrt{2\sigma|\phi(1)|}}{\lambda_2^0-\lambda_1^0}\right]^2
=\tfrac{2\sigma h^m}{\lambda_2^0-(\sigma+1)\lambda_1^0}
+\left[\tfrac{\sqrt{2\sigma(\sigma+1)|\phi(1)|}}{\lambda_2^0-(\sigma+1)\lambda_1^0}\right]^2-\tfrac{(L\lambda_2^0-\sqrt{2\sigma|\phi(1)|})^2}{(\lambda_2^0-(\sigma+1)\lambda_1^0)(\lambda_2^0-\lambda_1^0)}.
\ees

Next, expanding the squares on both sides of this equality and, subsequently, multiplying it by $(\lambda_2^0-(\sigma+1)\lambda_1^0)^2(\lambda_2^0-\lambda_1^0)$ yields after few term cancellations
\be
L^2|\lambda_1^0|^2|\lambda_2^0|^2-2\sigma|\phi(1)|\lambda_1^0\lambda_2^0=2h^m\lambda_1^0(\lambda_2^0-(\sigma+1)\lambda_1^0)(\lambda_2^0-\lambda_1^0).
\lb{lambda_12_1}
\ee
Similarly, substituting \rf{sx_expr} into \rf{QP2} one derives another relation
\bes
-L^2|\lambda_1^0|^2|\lambda_2^0|^2+2\sigma|\phi(1)|\lambda_1^0\lambda_2^0=-2h_1^m\lambda_2^0(\lambda_2^0-(\sigma+1)\lambda_1^0)(\lambda_2^0-\lambda_1^0).
\ees
Adding the last two relations gives
\bes
2(h^m\lambda_1^0-h_1^m\lambda_2^0)(\lambda_2^0-(\sigma+1)\lambda_1^0)(\lambda_2^0-\lambda_1^0)=0,
\ees
implying that
\be
\frac{\lambda_2^0}{\lambda_1^0}=\frac{h^m}{h_1^m}=\overline{h}.
\lb{lambda_r}
\ee
Finally, dividing \rf{lambda_12_1} by $|\lambda_1^0|^2\not=0$ and using \rf{lambda_r} yields a quadratic equation for $\lambda_2^0$:
\bes
L^2|\lambda_2^0|^2-2h_1^m(\overline{h}-1)(\overline{h}-\sigma-1)\lambda_2^0-2\sigma|\phi(1)|\overline{h}=0,
\ees 
having single positive root:
\be
\lambda_2^0=\frac{h_1^m(\overline{h}-1)(\overline{h}-\sigma-1)}{L^2}\left[1+\sqrt{1+\frac{2\sigma|\phi(1)|L^2\overline{h}}{|h_1^m|^2(\overline{h}-1)^2(\overline{h}-\sigma-1)^2}}\right].
\lb{lambda2_zz}
\ee
Note that taking $\lambda_2^0<0$, instead, would imply $\lambda_1^0<0$ by \rf{lambda_r} and, subsequently, contradict to conditions $s_1>s>0$ (this can be checked using explicit expressions \rf{sx_expr}). 

In summary, formulae \rf{sx_expr}, \rf{lambda_r} and \rf{lambda2_zz} provide a complete information about the leading order (as $\eps\go0$) profile of the derived {\it zig-zag} solution with typical shapes of which being shown in Fig.6. 
They also cover the two limiting cases $\lambda_2^0=\lambda_1^0$ and $\lambda_2^0=(1+\sigma)\lambda_1^0$ in system \rf{1m0sol}, which correspond to 
$\overline{h}=1$ or $\overline{h}=\sigma+1$ values in \rf{lambda_r}, respectively. In these cases, profile of $h_1^0(x)$ or $h^0(x)$, in the bulk region $x\in(-s_1,\,-s)$ is given by a linear segment and 
\bes
\lambda_1^0=\tfrac{\sqrt{2\sigma|\phi(1)|}}{L}\quad\text{or}\ \lambda_1^0=\tfrac{\sqrt{2\sigma|\phi(1)|/(\sigma+1)}}{L}\quad\text{holds, respectively}.
\ees
The corresponding limiting formulae for positions $s,\,s_1,\,x_c,\,x_{c1}$ can be derived from \rf{sx_expr} then.

The constraints on the {\it zig-zag} solution are deduced from the length ones 
\be
L>s_1>s>0
\lb{LC}
\ee
%%%%%%%%%%%%%%%%%%%%%%%%%%%%%%%%%%%%%%%%%%%%%%%%%%%%%%%%%%%%%%%%%%%%%%%%%%%%%%%%%%%%%%%%%%%%%%%%%%%%%%%%%%%%%%%%
\begin{figure}[H] 
	\centering
	\vspace{-2.9cm}
	\hspace{-.4cm}\includegraphics[width=.5\textwidth]{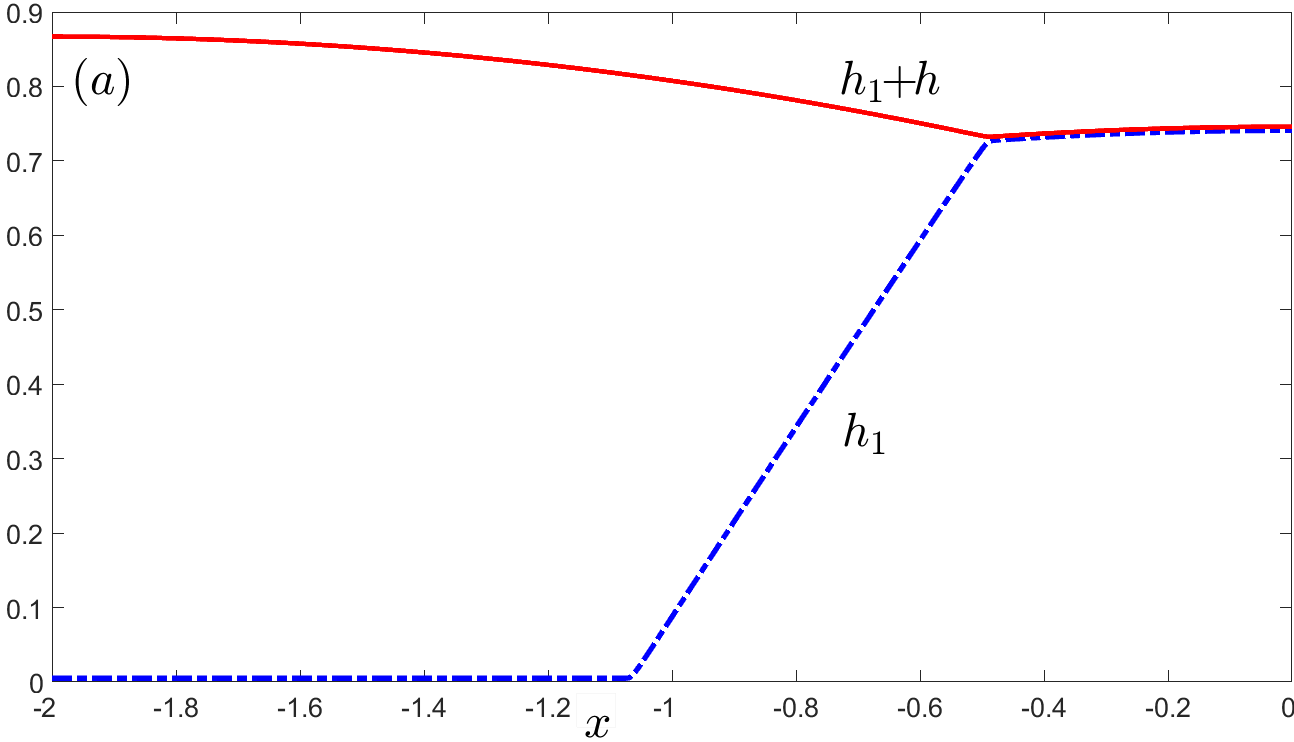}  
	\hspace{.2cm}\includegraphics[width=.5\textwidth]{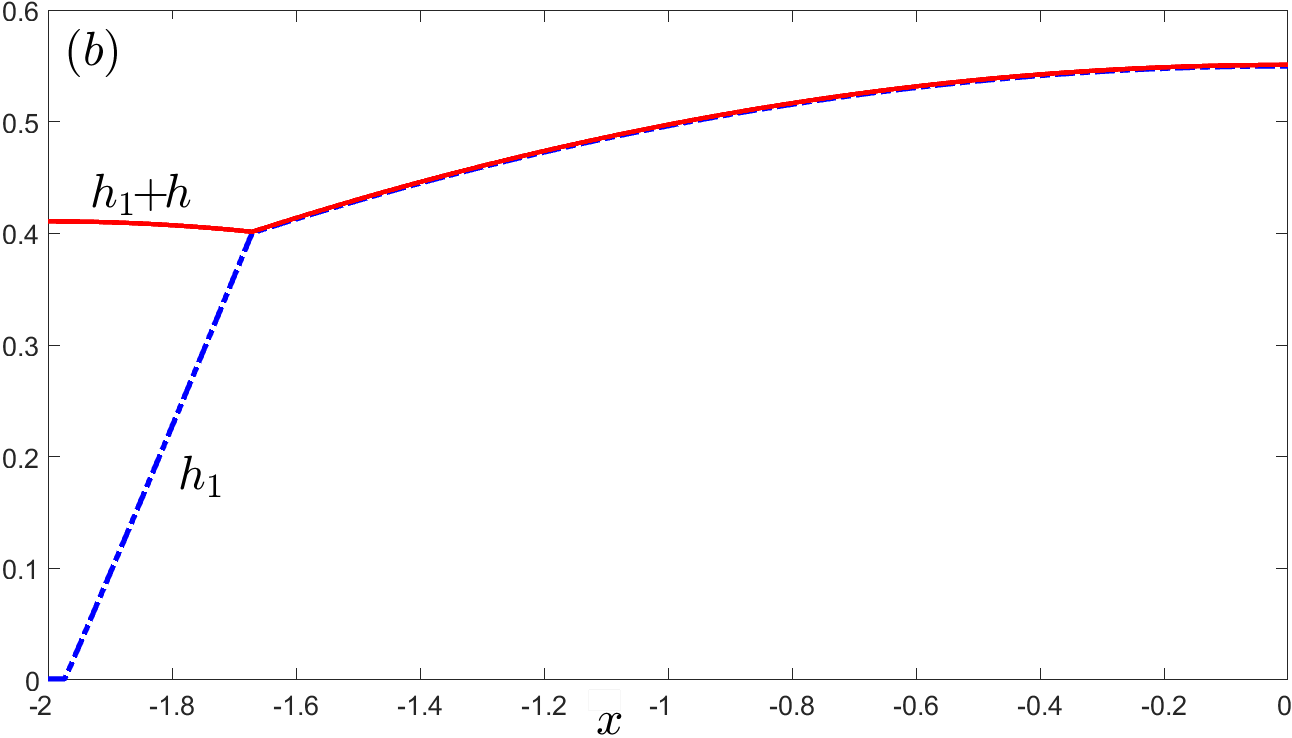}  
	
	\vspace{.2cm}
	\hspace{-.4cm}\includegraphics[width=.5\textwidth]{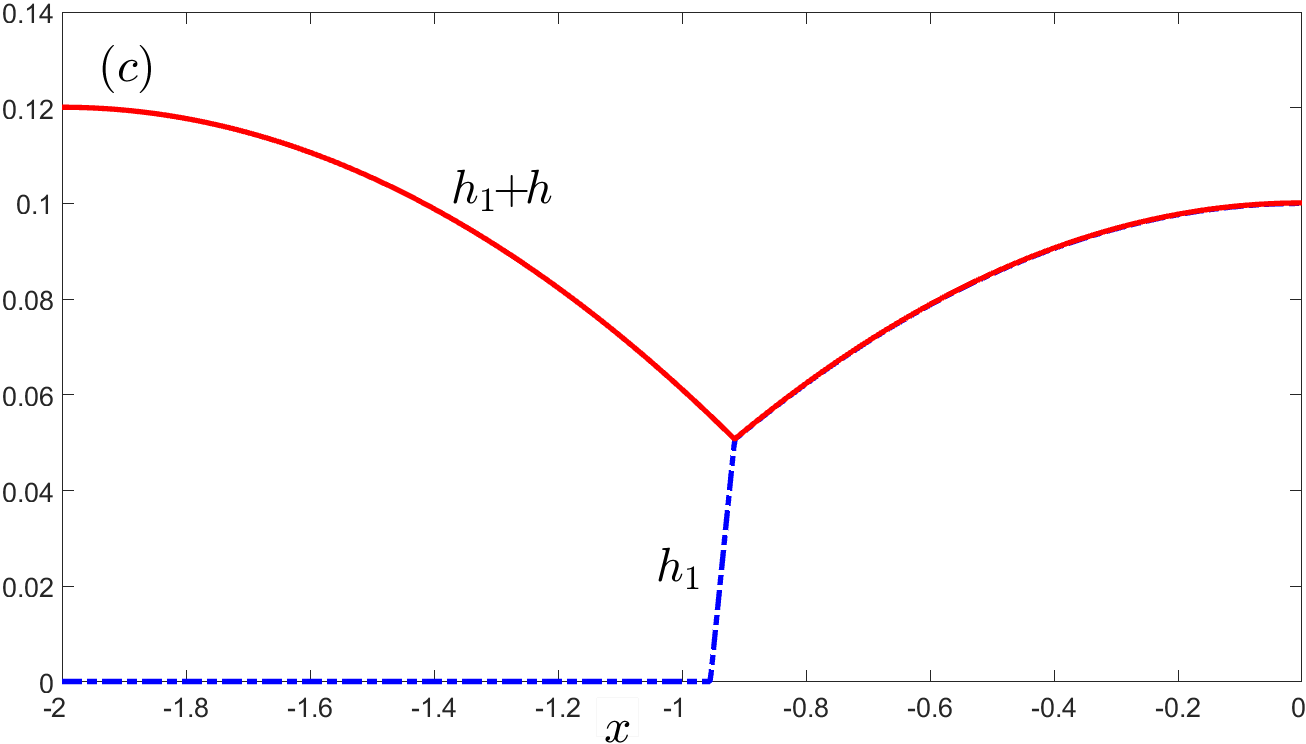}  
	\hspace{.2cm}\includegraphics[width=.5\textwidth]{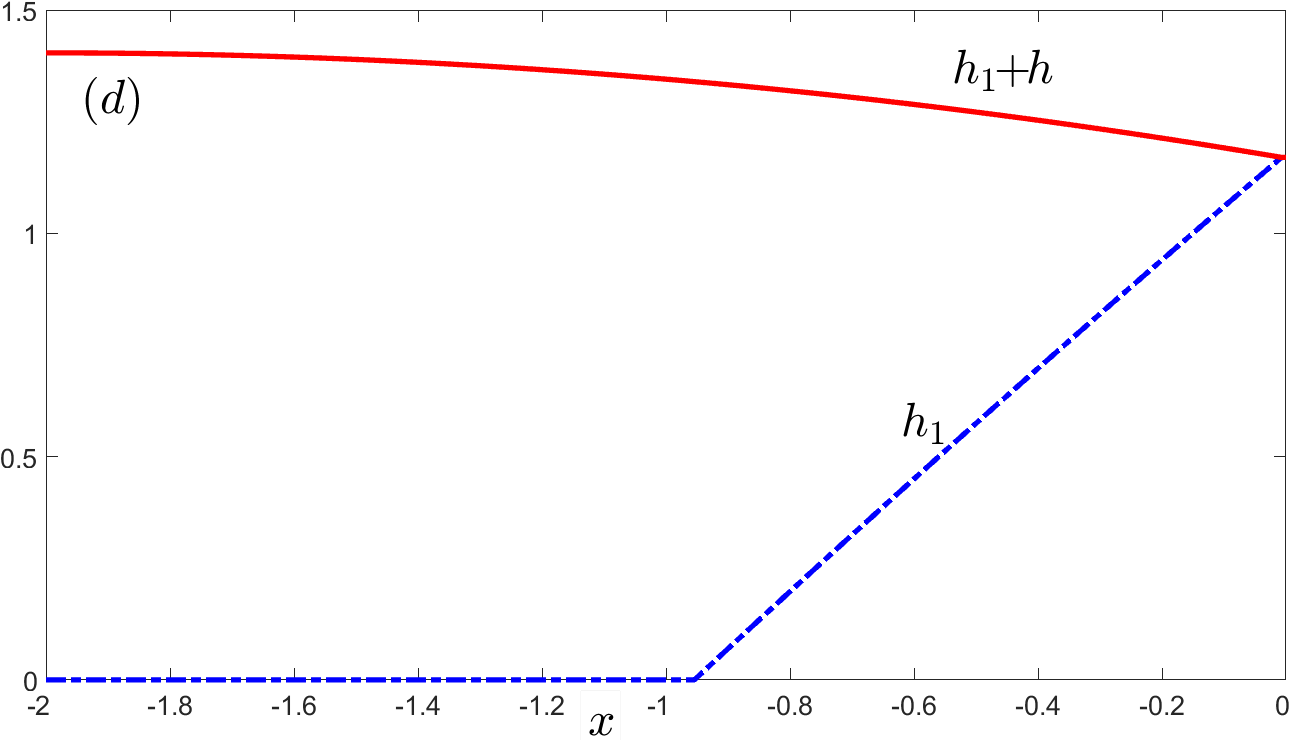}  
	
	\caption{\small\hspace{-.15cm}Numerical zig-zag stationary solutions to system \rf{BS}--\rf{BC} for $\eps\!=\!0.001$, $\sigma\!=\!0.2,\,L\!=\!2.0$ with observed: {\bf (a)} $\lambda_1\!=\!0.1188,\, \lambda_2\!=\!0.1383$; {\bf (b)} $\lambda_1\!=\!0.1726,\, \lambda_2\!=\!0.1285$; {\bf (c)} $\lambda_1\!=\!0.11786,\, \lambda_2\!=\!(\sigma+1)\lambda_1$; {\bf (d)} $\lambda_1\!=\!0.11785,\, \lambda_2\!=\!(\sigma+1)\lambda_1$.}	
\end{figure}
%%%%%%%%%%%%%%%%%%%%%%%%%%%%%%%%%%%%%%%%%%%%%%%%%%%%%%%%%%%%%%%%%%%%%%%%%%%%%%%%%%%%%%%%%%%%%%%%%%%%%%%%%%%%%%%%
\noindent as follows. First, note that from the $2^{\mathrm{nd}}$, 6$^{\mathrm{th}}$ and 10$^{\mathrm{th}}$ equation of system \rf{1m0sol} it follows that $h^0(x)$ is monotonically decreasing function. Similarly, from the 12$^{\mathrm{th}}$, 9$^{\mathrm{th}}$ and 5$^{\mathrm{th}}$ equations in \rf{1m0sol} it follows that $h_1^0(x)$ is monotonically increasing. Therefore, if the length conditions \rf{LC} are assured then the leading order profile of {\it zig-zag} solution described above is well defined.

For instance, let $\overline{h}>\sigma+1$. Using formulae \rf{sx_expr} conditions \rf{LC} can be written then as
\bes
\tfrac{L\lambda_2^0}{\sqrt{2\sigma|\phi(1)|}}>1,\quad
\tfrac{L\lambda_2^0}{\sqrt{2\sigma|\phi(1)|}}>\tfrac{\overline{h}+\sqrt{\sigma+1}}{\sqrt{\sigma+1}+1},\quad
\tfrac{L\lambda_2^0}{\sqrt{2\sigma|\phi(1)|}}<\tfrac{\overline{h}}{\sqrt{\sigma+1}}.
\ees 
As $\overline{h}+\sqrt{\sigma+1}>1+\sqrt{\sigma+1}$, these inequalities reduce to constraint
\be
\tfrac{L\lambda_2^0}{\sqrt{2\sigma|\phi(1)|}}\in\left(\tfrac{\overline{h}+\sqrt{\sigma+1}}{\sqrt{\sigma+1}+1},\,\tfrac{\overline{h}}{\sqrt{\sigma+1}}\right).
\lb{Ms}
\ee
Next, substituting the expression \rf{lambda2_zz} for $\lambda_2^0$ in \rf{Ms} and using $\overline{h}>1$ allow us to deduce
\be
\hspace{-2.cm}L\sqrt{2\sigma|\phi(1)|}\in\left(2\sqrt{\sigma+1}h_1^m(\overline{h}-1),\,\frac{\sigma h_1^m(\overline{h}-1)(\overline{h}-\sigma-1)(\overline{h}+\sqrt{\sigma+1})(\sqrt{\sigma+1}-1)}{(\overline{h}+\sqrt{\sigma+1})^2(\sqrt{\sigma+1}-1)^2-\sigma^2\overline{h}}\right).
\lb{Ms3}
\ee
Finally, using decomposition $$(\overline{h}+\sqrt{\sigma+1})^2(\sqrt{\sigma+1}-1)^2-\sigma^2\overline{h}=(\sqrt{\sigma+1}-1)^2(\overline{h}-1)(\overline{h}-\sigma-1)$$
simplifies \rf{Ms3} further and implies the following bounds on the interval length $L$ from above and below:
\bes
L\in \tfrac{\sqrt{2}h_1^m}{\sqrt{\sigma|\phi(1)|}}\left(\sqrt{\sigma+1}(\overline{h}-1),\,(\sqrt{\sigma+1}+1)(\overline{h}+\sqrt{\sigma+1})\right).
\ees

Similar constraints on $L$ can be deduced using formulae \rf{sx_expr} and \rf{lambda2_zz} for the other two ranges $\overline{h}\in(1,\,\sigma+1)$ and $\overline{h}\in(0,\,1)$. We omit the details and just state their combined final form being valid for all positive $\overline{h}$:
\be
L\in\left\{\begin{array}{ll}\tfrac{\sqrt{2}h_1^m}{\sqrt{\sigma|\phi(1)|}}\left(\sigma+1-\overline{h},\,(\sqrt{\sigma+1}+1)(\overline{h}+\sqrt{\sigma+1})\right),&\text{if}\ \overline{h}\in(0,\,\sqrt{\sigma+1}],\\[3ex]
\tfrac{\sqrt{2}h_1^m}{\sqrt{\sigma|\phi(1)|}}\left(\sqrt{\sigma+1}(\overline{h}-1),\,(\sqrt{\sigma+1}+1)(\overline{h}+\sqrt{\sigma+1})\right),&\text{if}\ \overline{h}\in[\sqrt{\sigma+1},\,+\infty).\end{array}\right.
\lb{zz_constr}
\ee
Constraints \rf{zz_constr} cover also the two limiting cases $\lambda_2^0\!=\!\lambda_1^0$ and $\lambda_2^0\!=\!(1+\sigma)\lambda_1^0$ occurring for $\overline{h}\!=\!1$ and $\overline{h}\!=\!\sigma+1$ values, respectively. Note that the interval boundaries $\sigma+1-\overline{h}$, $\sqrt{\sigma+1}(\overline{h}-1)$ and $(\sqrt{\sigma+1}+1)(\overline{h}+\sqrt{\sigma+1})$ in \rf{zz_constr} correspond to the critical merges $s_1\go L$, $s\go0$ and $s_1\go s$, respectively, occurring when one of the bulk regions in {\it zig-zag} solution shrinks to zero (cf. Fig.6 {\bf (b)-(d)}). Accordingly, special value $\overline{h}=\sqrt{\sigma+1}$ in \rf{zz_constr} sets a threshold between the cases when either $s_1\!\go\!L$ or $s\!\go\!0$ merge occurs first by decreasing $L$. Additionally, note that merge  $s_1\go s$ occurs simultaneously with $h_1^0(-s)\go0$ and $h^0(-s_1)\go0$, i.e. the solution splits into a {\it 2-drops} one (described below in this section) then (cf. Fig.6 {\bf(c)}).\\[.5ex]

\underline{\bf Sessile lens:}\hspace{.5cm}The solution is obtained by a combined matching of {\bf Type I} (with $x_{c1}=x_c=-L$) and general {\bf Type II}  bulk solutions to {\bf Type II} CL ones centered around $x=-s_{1}$ ($1^{\mathrm{st}}$ CL) as well as of general {\bf Type II}  bulk  and UTF solutions to {\bf Type III} CL ones centered around $x=-s$ ($2^{\mathrm{nd}}$ CL). Schematically this matching chain can be encrypted as $\bf{IB-IICL-IIB-IIICL-UTF}$ once moving in $x$ from $-L$ to $0$.

After the matching procedures at contact line points $x=-s_1$ and $x=-s$ the resulting system of the leading order conditions takes the form:
\be
\left\{\begin{array}{l}
	\frac{\lambda_1^0-\lambda_2^0}{2\sigma}(-s_1+L)^2+C_1=C_4,\hspace{1.4cm}
	\frac{\lambda_2^0-(\sigma+1)\lambda_1^0}{2\sigma}(-s_1+L)^2+C=0,\\[1.5ex]
	\frac{\lambda_1^0-\lambda_2^0}{\sigma}(-s_1+L)=\sqrt{\frac{2|\phi(1)|}{\sigma(\sigma+1)}}+C_5,\hspace{.5cm}
	\frac{\lambda_2^0-(\sigma+1)\lambda_1^0}{\sigma}(-s_1+L)=-\sqrt{\frac{2(\sigma+1)|\phi(1)|}{\sigma}},\\[1.5ex]
	-\frac{\lambda_2^0}{2(\sigma+1)}(s_1+\tilde{x}_{c1})^2+\widetilde{C}_1=C_4,\hspace{1.cm}
	-\frac{\lambda_2^0}{\sigma+1}(-s_1-\tilde{x}_{c1})=C_5\\[1.5ex]
	-\frac{\lambda_2^0}{2(\sigma+1)}(s+\tilde{x}_{c1})^2+\widetilde{C}_1=0,\hspace{1.5cm}
	-\frac{\lambda_2^0}{\sigma+1}(-s-\tilde{x}_{c1})=-\sqrt{\frac{2|\phi(1)|}{\sigma+1}}.\end{array}\right.
\lb{02sol}
\ee  
Additionally, we fix $C=h^m$ and $C_1=h_1^m$. System \rf{02sol} has $8$ equations with $8$ unknowns: $\lambda_1^0,\,\lambda_2^0,\,s,\,s_1,\,\tilde{x}_{c1}$ and $C_4,\,C_5,\,\widetilde{C}_1$. From the linear part of system \rf{02sol}, i.e. by performing joint manipulations with its $3^{\mathrm{th}}-4^{\mathrm{th}}$ and 6$^{\mathrm{th}}$ equations, one deduces that $\tilde{x}_{c1}=-L$. Next, from the last two equations in \rf{02sol} one obtains $\widetilde{C}_1=|\phi(1)|/\lambda_2^0$. Substituting that into the $5^{\mathrm{th}}$ equation in  \rf{02sol} and, subsequently, subtracting from it the $1^{\mathrm{st}}$ and $2^{\mathrm{nd}}$ equations yields a relation
\be
-h_1^m-\tfrac{h^m}{\sigma+1}+\tfrac{|\phi(1)|}{\lambda_2^0}=0,\quad\text{implying}\ 
\lambda_2^0=\tfrac{(\sigma+1)|\phi(1)|}{(\sigma+1)h_1^m+h^m}.
\lb{lambda2_gl}
\ee
Next, from the $2^{\mathrm{nd}}$ and $4^{\mathrm{th}}$ equations in \rf{02sol} one deduces a relation
\be
\lambda_2^0-(\sigma+1)\lambda_1^0=-\tfrac{2|\phi(1)|}{(\sigma+1)h^m},
\lb{Ms4}
\ee
which using \rf{lambda2_gl} gives 
\be
\lambda_1^0=\tfrac{(\sigma+1)h_1^m+2h^m}{(\sigma+1)h_1^m+h^m}\cdot\tfrac{|\phi(1)|}{h^m}.
\lb{lambda1_gl}
\ee
Finally, substituting \rf{lambda2_gl}--\rf{lambda1_gl} into the $4^{\mathrm{th}}$ and  $8^{\mathrm{th}}$ equations in \rf{02sol} yields the expressions for the contact line positions:
\be
\lb{s_gl}
s_1=L-h^m\sqrt{\tfrac{2\sigma}{(\sigma+1)|\phi(1)|}},\quad
s=L-((\sigma+1)h_1^m+h^m)\sqrt{\tfrac{2}{(\sigma+1)|\phi(1)|}}.
\ee

In summary, formulae $\tilde{x}_{c1}=-L$ and \rf{lambda2_gl}--\rf{s_gl} provide a complete information about the derived {\it sessile lens} solution, typical shapes of which are shown in Fig.7.
In particular, from the $1^{\mathrm{st}}$ equation in \rf{02sol} one gets
\be
C_4=\tfrac{h^m}{\sigma+1}\cdot\tfrac{(1-\sigma)h^m+(\sigma+1)h_1^m}{(\sigma+1)h_1^m+h^m}+h_1^m,
\lb{C4}
\ee
implying the following combined formula for the maximum value of $h_1^0$:
\be
\displaystyle\max_{x\in(-L,\,0)}h_1^0(x)=\left\{\begin{array}{ll}
	h_1^0(0)=h_1^m,& \text{if}\ \sigma>1\,\text{and}\, \frac{h^m}{h_1^m}\ge\tfrac{\sigma+1}{\sigma-1},\\[1.5ex]
	h_1^0(-s_1)=C_4,& \text{otherwise}.\end{array}\right.
\lb{h1m_ls}
\ee  
Additionally, for $\sigma>1$ \rf{lambda2_gl} and \rf{lambda1_gl} imply the presence of critical case $\lambda_2^0=\lambda_1^0$ for which $\bar{h}=\frac{h^m}{h_1^m}=\tfrac{\sigma+1}{\sigma-1}$ and  $h_1^0(x)=h_1^m=\const$ holds for all $x\in(-L,\,-s_1)$, i.e. the leading order profile of $h_1$ in the left bulk region is {\it flat} (cf. Fig.7 {\bf (c)}). From \rf{Ms4} we observe that the second critical case $\lambda_2^0=(\sigma+1)\lambda_1^0$ is not possible here.

The constraints on the {\it sessile lens} solution are deduced again from conditions \rf{LC}
using \rf{s_gl} and take the form:
\be
\left\{\begin{array}{c}L>((\sigma+1)h_1^m+h^m)\sqrt{\tfrac{2}{(\sigma+1)|\phi(1)|}},\\[2ex]
		\tfrac{h^m}{h_1^m}<\tfrac{\sigma+1}{\sqrt{\sigma}-1}\quad\text{if}\ \sigma>1.
\end{array}\right.
\lb{gl_constr}
\ee
The first condition in \rf{gl_constr} provides the bound on the minimal interval length $L$, while the second one ensures $s_1>s$ and is required only if $\sigma>1$. The threshold 
$\bar{h}=h^m/h_1^m=\frac{\sigma+1}{\sqrt{\sigma}-1}$ corresponds to the contact line merge $s_1\go s$ occurring simultaneously with $h_1^0(-s_1)\go0$ (cf. \rf{C4}--\rf{h1m_ls} and Fig.7 {\bf (d)}).\\[.5ex] 

\underline{\bf Sessile internal drop:}\hspace{.5cm}The solution is obtained by a combined matching of {\bf Type I} (with $x_{c1}=x_c=-L$) and general {\bf Type III}  bulk solutions to {\bf Type I} CL ones centered around $x=-s_{1}$ ($1^{\mathrm{st}}$ CL) as well as of general {\bf Type III}  bulk  and {\bf UTF} solutions to {\bf Type IV} CL ones centered around $x=-s$ ($2^{\mathrm{nd}}$ CL). 
%%%%%%%%%%%%%%%%%%%%%%%%%%%%%%%%%%%%%%%%%%%%%%%%%%%%%%%%%%%%%%%%%%%%%%%%%%%%%%%%%%%%%%%%%%%%%%%%%%%%%%%%%%%%%%%%
\begin{figure}[H] 
	\centering
	\vspace{-2.9cm}
	\hspace{-.4cm}\includegraphics[width=.5\textwidth]{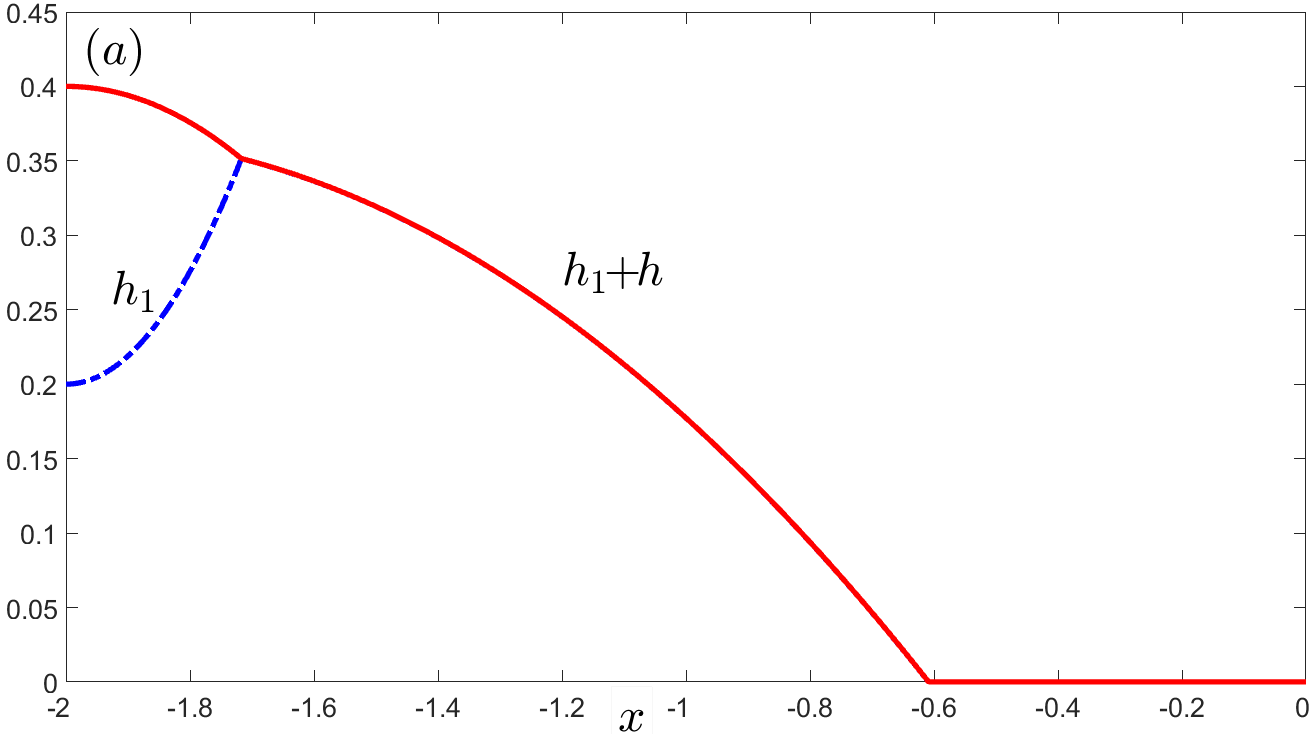}  
	\hspace{.2cm}\includegraphics[width=.5\textwidth]{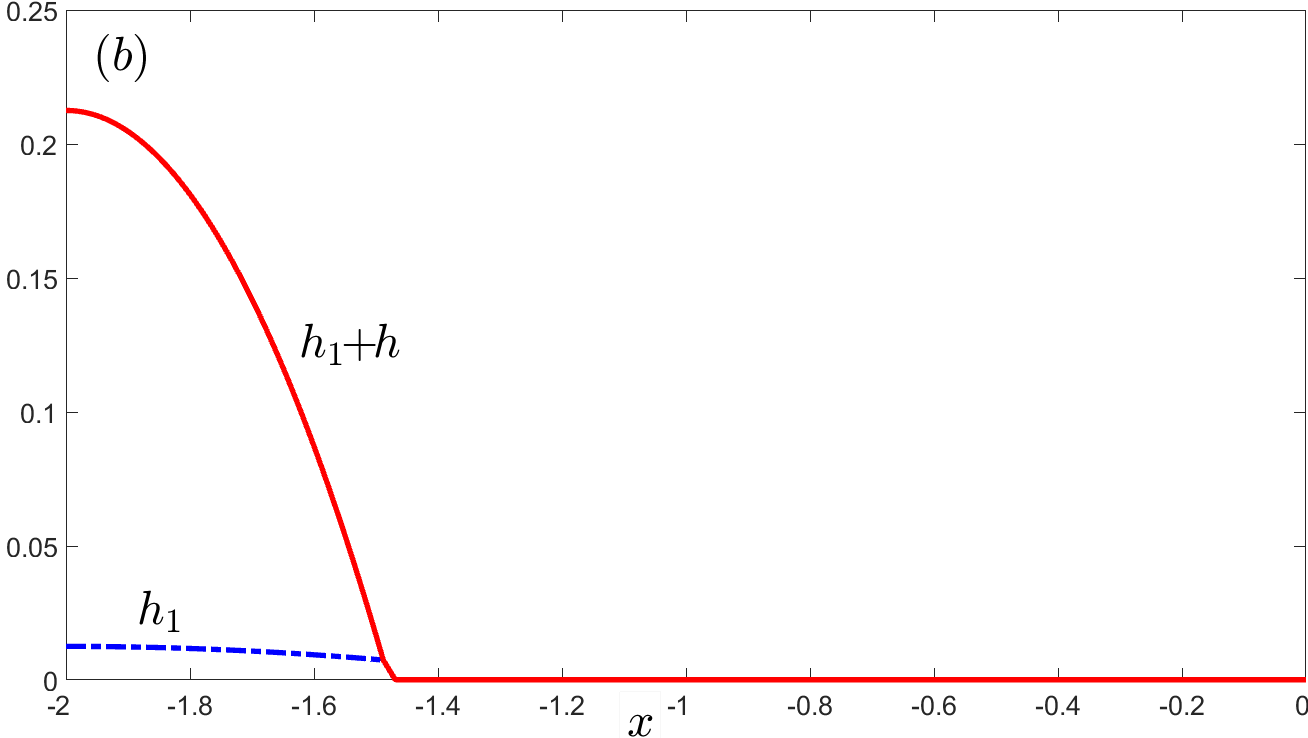}  
	
	\vspace{.2cm}
	\hspace{-.4cm}\includegraphics[width=.5\textwidth]{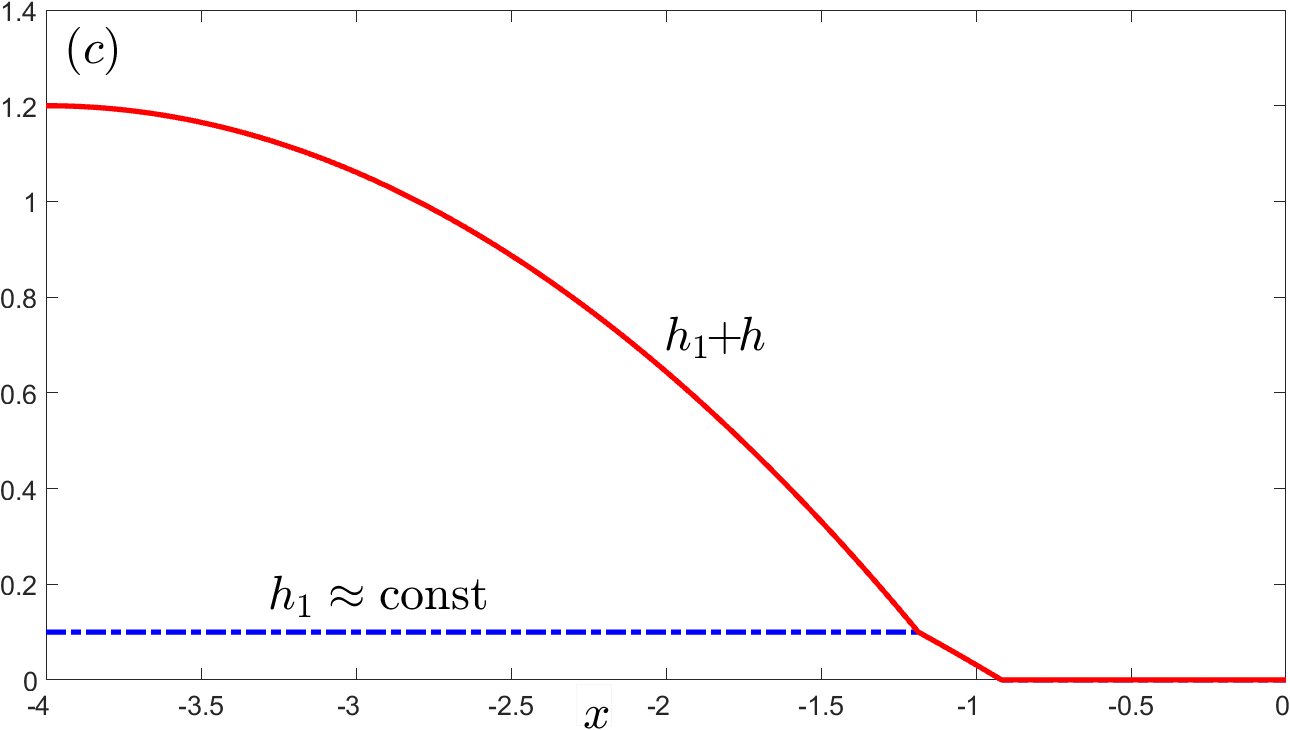}  
	\hspace{.2cm}\includegraphics[width=.5\textwidth]{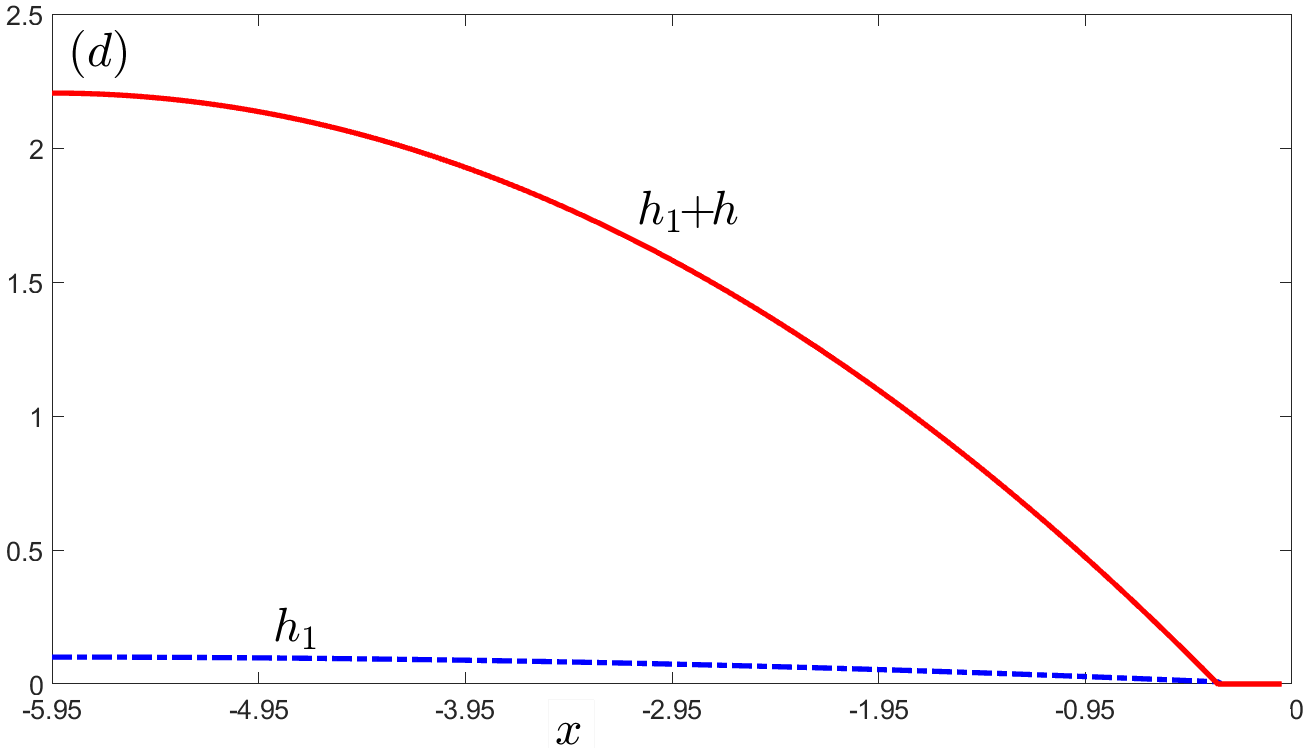}  
	
	\caption{\small Numerical sessile lens stationary solution to system \rf{BS}--\rf{BC} for\hspace{2.cm}$\eps\!=\!0.00005$: {\bf (a)} $\sigma\!=\!0.2,\,L\!=\!2$ and $\lambda_1\!=\!1.212163,\, \lambda_2\!=\!0.454598$; {\bf (b)} $\sigma\!=\!1.2,\,L\!=\!2$ and $\lambda_1\!=\!1.566342,\, \lambda_2\!=\!1.614762$; {\bf (c)} $\sigma\!=\!1.2,\,L\!=\!4$ and $\lambda_1\!=\!0.277782$, $\lambda_2\!=\!0.277799$; {\bf (d)} $\sigma\!=\!1.2,\,L\!=\!5.95$ and $\lambda_1\!=\!0.144377,\, \lambda_2\!=\!0.151355$.}	
\end{figure}
%%%%%%%%%%%%%%%%%%%%%%%%%%%%%%%%%%%%%%%%%%%%%%%%%%%%%%%%%%%%%%%%%%%%%%%%%%%%%%%%%%%%%%%%%%%%%%%%%%%%%%%%%%%%%%%%
\noindent The matching chain can be encrypted as $\bf{IB-ICL-IIIB-IVCL-UTF}$ once moving in $x$ from $-L$ to $0$.
After the matching procedures at contact lines $x=-s_1$ and $x=-s$ the resulting system of the leading order conditions takes the form:
\be
\left\{\begin{array}{l}
	\frac{\lambda_1^0-\lambda_2^0}{2\sigma}(-s_1+L)^2+C_1=0,\hspace{1.3cm}
	\frac{\lambda_2^0-(\sigma+1)\lambda_1^0}{2\sigma}(-s_1+L)^2+C=C_2,\\[1.5ex]
	\frac{\lambda_1^0-\lambda_2^0}{\sigma}(-s_1+L)=-\sqrt{\frac{2|\phi(1)|}{\sigma}},\hspace{.8cm}
	\frac{\lambda_2^0-(\sigma+1)\lambda_1^0}{\sigma}(-s_1+L)=\sqrt{\frac{2|\phi(1)|}{\sigma}}+C_3,\\[1.5ex]
	-\frac{\lambda_1^0}{2}(s_1+\tilde{x}_c)^2+\widetilde{C}_0=C_2,\hspace{1.4cm}
	-\lambda_1^0(-s_1-\tilde{x}_c)=C_3\\[1.5ex]
	-\frac{\lambda_1^0}{2}(s+\tilde{x}_c)^2+\widetilde{C}_0=0,\hspace{1.8cm}
	-\lambda_1^0(-s-\tilde{x}_c)=-\sqrt{2|\phi(1)|}.\end{array}\right.
\lb{13sol}
\ee  
Additionally, we fix $C=h^m$ and $C_1=h_1^m$. System \rf{13sol} has $8$ equations with $8$ unknowns: $\lambda_1^0,\,\lambda_2^0,\,s,\,s_1,\,\tilde{x}_c$ and $C_2,\,C_3,\,\widetilde{C}_0$. Its solution proceeds
similarly to the one of system \rf{02sol} for the {\it sessile lens} solution. From the linear part of system \rf{13sol} one deduces that $\tilde{x}_c=-L$. Next, from the last two equations in \rf{02sol} one obtains $\widetilde{C}_0=|\phi(1)|/\lambda_1^0$. Substituting that into the $5^{\mathrm{th}}$ equation in  \rf{02sol} and, subsequently, subtracting from it the $1^{\mathrm{st}}$ and $2^{\mathrm{nd}}$ equations one arrives at relation
\be
-h_1^m-h^m+\tfrac{|\phi(1)|}{\lambda_1^0}=0\quad\text{implying}\ \lambda_1^0=\tfrac{|\phi(1)|}{h_1^m+h^m}.
\lb{lambda1_gdr}
\ee
Next, from the $1^{\mathrm{st}}$ and $3^{\mathrm{rd}}$ equations in \rf{13sol} one deduces
\be
\lambda_1^0-\lambda_2^0=-\tfrac{|\phi(1)|}{h_1^m}\quad\text{impliyng}\ 
\lambda_2^0=\tfrac{h^m+2h_1^m}{h_1^m+h^m}\cdot\tfrac{|\phi(1)|}{h_1^m}.
\lb{lambda2_gdr} 
\ee
Finally, substituting \rf{lambda1_gdr}--\rf{lambda2_gdr} into the $3^{\mathrm{rd}}$ and  $8^{\mathrm{th}}$ equations in \rf{13sol} yields:
\be
\lb{s_gdr}
s_1=L-h_1^m\sqrt{\tfrac{2\sigma}{|\phi(1)|}},\quad
s=L-(h_1^m+h^m)\sqrt{\tfrac{2}{|\phi(1)|}}.
\ee
%%%%%%%%%%%%%%%%%%%%%%%%%%%%%%%%%%%%%%%%%%%%%%%%%%%%%%%%%%%%%%%%%%%%%%%%%%%%%%%%%%%%%%%%%%%%%%%%%%%%%%%%%%%%%%%%
\begin{figure}[H] 
	\centering
	\vspace{-2.9cm}
	\hspace{-.4cm}\includegraphics[width=.5\textwidth]{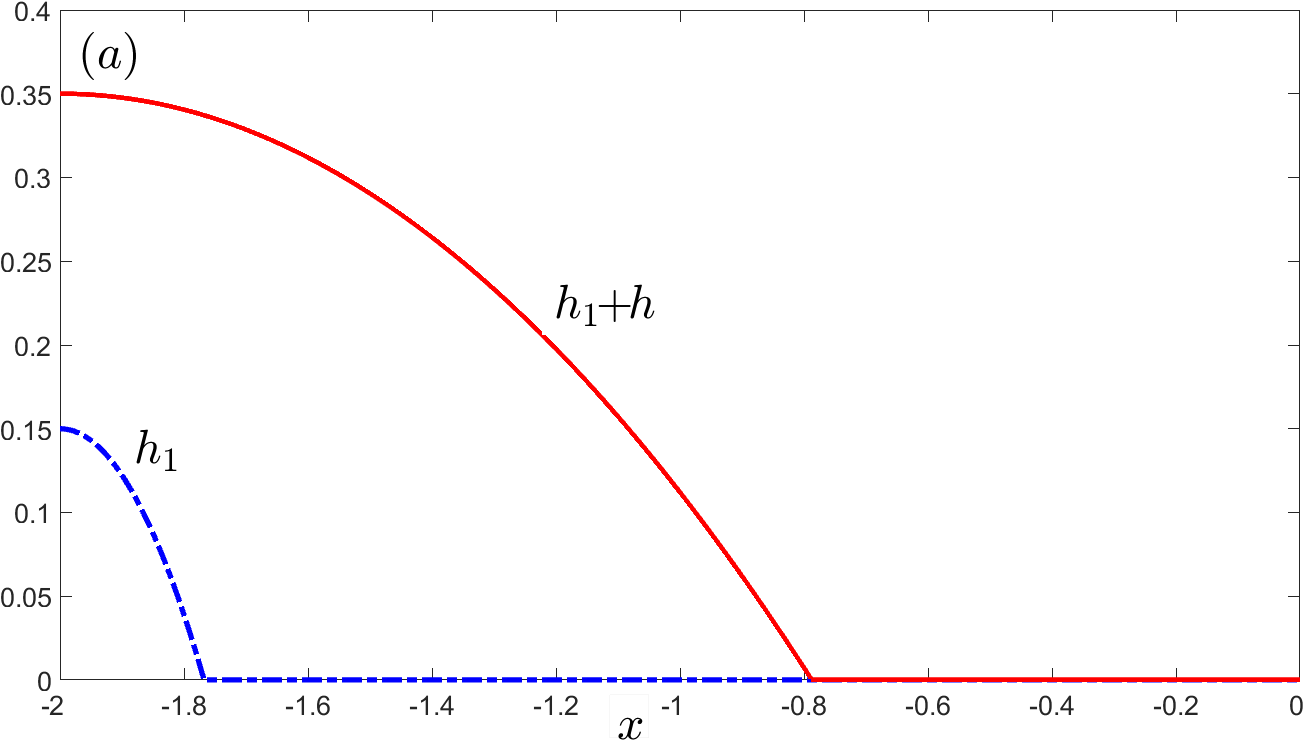}  
	\hspace{.2cm}\includegraphics[width=.5\textwidth]{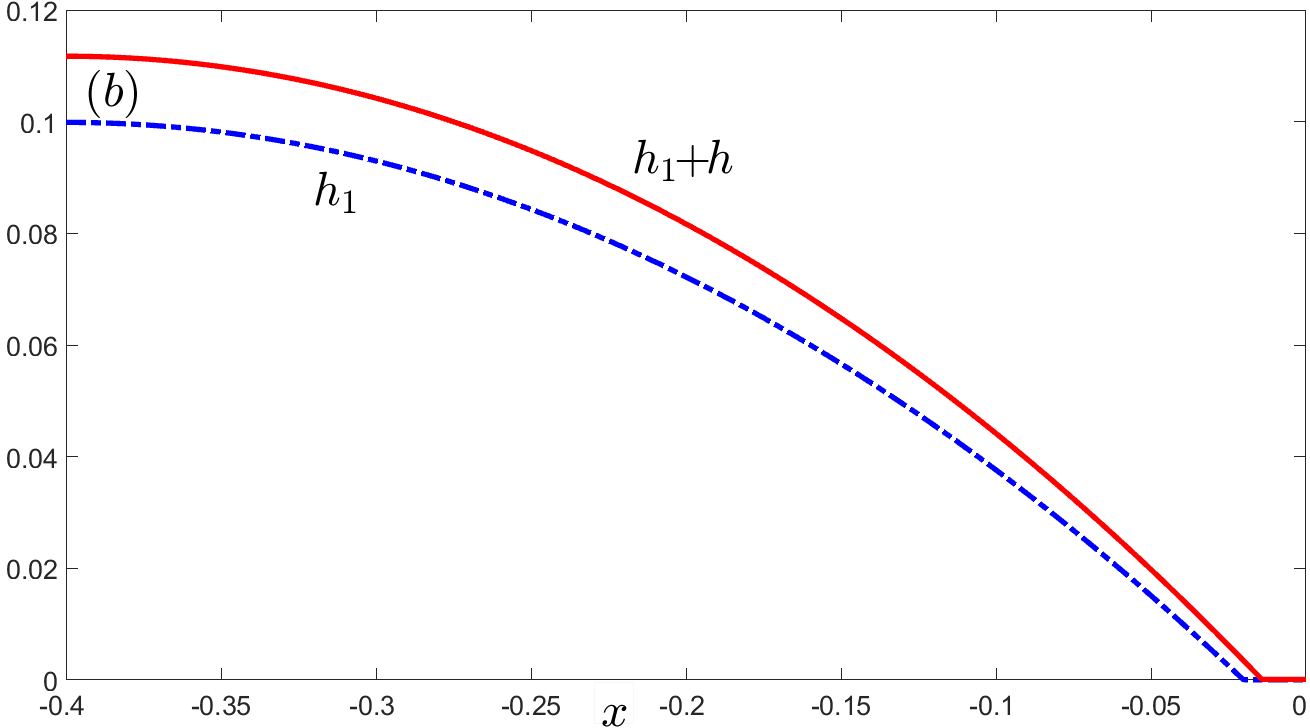}  
	
	\caption{\small Numerical sessile internal drop stationary solution to system \rf{BS}--\rf{BC} for $\eps\!=\!0.00005$: {\bf (a)} $\sigma\!=\!0.2,\,L\!=\!2.0$ and $\lambda_1\!=\!0.476158,\, \lambda_2\!=\!1.587123$; {\bf (b)} $\sigma\!=\!1.2,\,L\!=\!0.4$ and $\lambda_1\!=\!1.50337,\, \lambda_2\!=\!3.16151$.}	
\end{figure}
%%%%%%%%%%%%%%%%%%%%%%%%%%%%%%%%%%%%%%%%%%%%%%%%%%%%%%%%%%%%%%%%%%%%%%%%%%%%%%%%%%%%%%%%%%%%%%%%%%%%%%%%%%%%%%%%
In summary, formulae $\tilde{x}_c=-L$ and \rf{lambda1_gdr}--\rf{s_gdr} provide a complete information about the derived {\it sessile internal drop} solution, typical shapes of which are shown in Fig.8. In particular, expressions \rf{lambda1_gdr}--\rf{lambda2_gdr} imply for $\sigma>1$ the presence of critical case $\lambda_2^0=(\sigma+1)\lambda_1^0$ for which $\bar{h}=h^m/h_1^m=\sigma-1$ and  $h^0(x)=h^m=\const$ holds for all $x\in(-L,\,-s_1)$, i.e. the leading order profile of $h$ in the left bulk region is {\it flat}. Accordingly, one has
\be
\displaystyle\max_{x\in(-L,\,0)}h^0(x)=\left\{\begin{array}{ll}
	h^0(-s_1)=C_2,& \text{if}\ \frac{h^m}{h_1^m}\ge \sigma-1,\\[1.5ex]
	h^0(0)=h^m,& \text{otherwise}.\end{array}\right.
\lb{hm_ids}
\ee
From \rf{lambda2_gdr} we note that the second critical case $\lambda_2^0=\lambda_1^0$ is not possible here.

The constraints on the {\it sessile internal drop} are deduced again from conditions \rf{LC}
using \rf{s_gdr} and take the form:
\be
\left\{\begin{array}{c}L>(h_1^m+h^m)\sqrt{\tfrac{2}{|\phi(1)|}},\\[2ex]
	\tfrac{h^m}{h_1^m}>\sqrt{\sigma}-1\quad\text{if}\ \sigma>1.
\end{array}\right.
\lb{gdr_constr}
\ee
The second condition in \rf{gdr_constr} ensures $s_1>s$ and is required only if $\sigma>1$. The threshold 
$\bar{h}=h^m/h_1^m=\sqrt{\sigma}-1$ corresponds to the contact line merge $s_1\go s$ with $h^0(-s_1)\go0$ (cf. Fig.8 {\bf (b)}).\\[.5ex]

\underline{\bf 2-drops:}\hspace{.5cm}The solution is obtained from simple concatenation of {\it $h$-drop} and {\it $h_1$-drop} solutions. The matching chain can be encrypted as  $\bf{IIIB-IVCL-UTF}$--$\bf{-IIICL-IIB}$ once moving in $x$ from $-L$ to $0$. The corresponding system of matching conditions is simply the union of \rf{h1sol} and \rf{hsol} ones with the former ones being shifted to $x=0$ instead of $x=-L$ and, subsequently, reflected around the origin, while in the latter ones $s$ replaced by $s_1$. The resulting system of four equations implies the following expressions:
\bes
\lambda_1^0=\frac{|\phi(1)|}{h^m},\quad \lambda_2^0=\frac{|\phi(1)|}{h_1^m},\quad
s_1=L-\sqrt{\tfrac{2}{|\phi(1)|}}h^m,\quad s=\sqrt{\tfrac{2(\sigma+1)}{|\phi(1)|}}h_1^m.
\ees
A typical shape of {\it 2-drops} solution is shown in Fig.9. The only constraint imposed on it is
\be
L>\sqrt{\tfrac{2}{|\phi(1)|}}(h^m+\sqrt{\sigma+1}h_1^m).
\lb{2dr_constr}
\ee
%%%%%%%%%%%%%%%%%%%%%%%%%%%%%%%%%%%%%%%%%%%%%%%%%%%%%%%%%%%%%%%%%%%%%%%%%%%%%%%%%%%%%%%%%%%%%%%%%%%%%%%%%%%%%%%%
\begin{figure}[H] 
	\centering
	\vspace{-2.9cm}
	\hspace{-.4cm}\includegraphics[width=.5\textwidth]{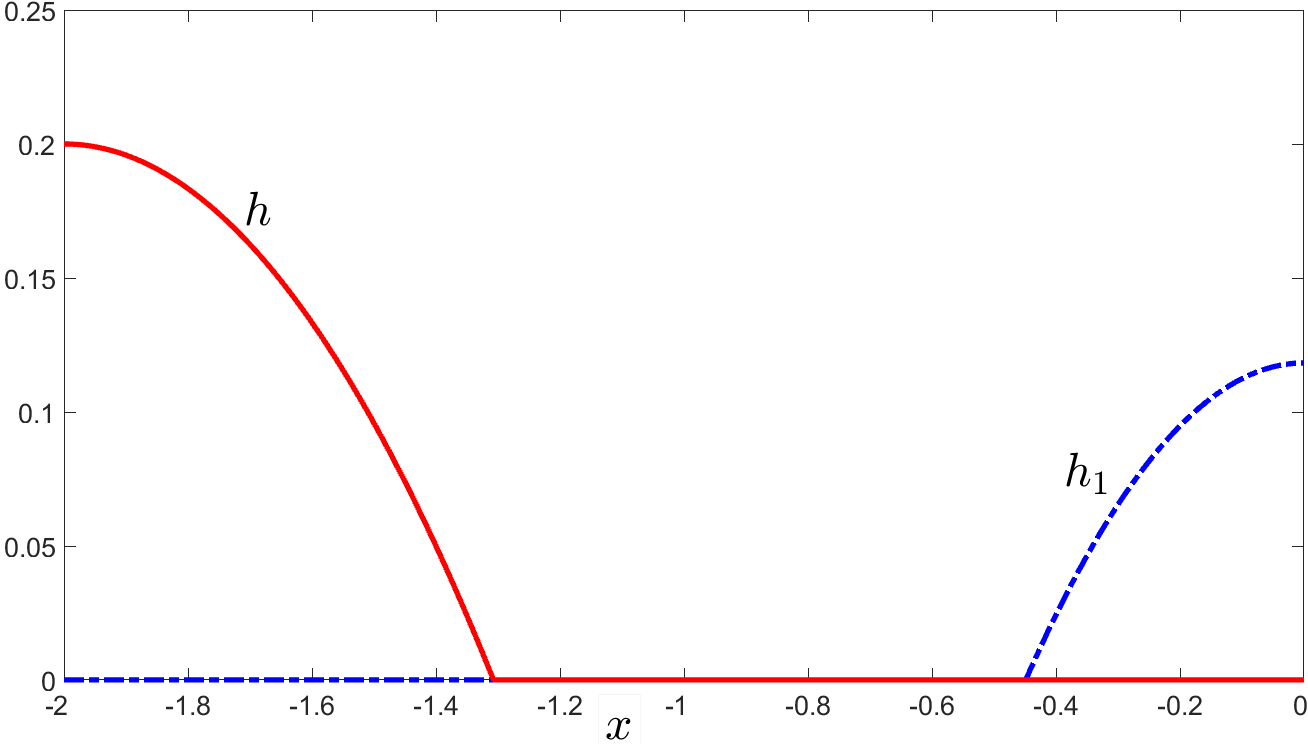}  
	\caption{\small Numerical 2-drops stationary solution to system \rf{BS}--\rf{BC} for\hspace{2.cm}$\eps\!=\!0.00005,\,\sigma\!=\!0.2,\,L\!=\!2.0$ with $\lambda_1\!=\!0.83349,\,\lambda_2\!=\!1.40943$.}
	
\end{figure}
%%%%%%%%%%%%%%%%%%%%%%%%%%%%%%%%%%%%%%%%%%%%%%%%%%%%%%%%%%%%%%%%%%%%%%%%%%%%%%%%%%%%%%%%%%%%%%%%%%%%%%%%%%%%%%%%

%%%%%%%%%%%%%%%%%%%%%%%%%%%%%%%%%%%%%%%%%%%%%%%%%%%%%%%%%%%%%%%%%
\section{Four-CL solutions}
%%%%%%%%%%%%%%%%%%%%%%%%%%%%%%%%%%%%%%%%%%%%%%%%%%%%%%%%%%%%%%%%%
For certain symmetry reasons we classify and describe the three contact line solutions in the next section. Before that in this section, we derive the leading order profiles of the solutions to system \rf{SSa}--\rf{SSc} having four contact lines. 
Up to possible inversion of $x$-variable we find a single solution type denoted below as {\it 2-side sessile zig-zag}.\\[.5ex]

\underline{\bf 2-side sessile zig-zag:}\hspace{.5cm}The solution is obtained by a combined matching of the following leading order profiles: {\bf UTF} and general {\bf Type II} bulk solutions to {\bf Type III} CL centered around $x=-s_3$ ($1^{\mathrm{st}}$ CL); general {\bf Type I} and {\bf Type II} bulk to {\bf Type II} CL ones centered around $x=-s_2$ ($2^{\mathrm{nd}}$ CL); general {\bf Type I} and {\bf Type III} bulk solutions to {\bf Type I} CL ones centered around $x=-s_1$ ($3^{\mathrm{rd}}$ CL); as well as {\bf UTF} and general {\bf Type III} bulk solutions to {\bf Type IV} CL centered around $x=-s$ ($4^{\mathrm{th}}$ CL). Schematically this matching chain can be encrypted as $\bf{UTF-IIICL-IIB-IICL-IB-ICL}$--$\bf{-IIIB-IVCL-UTF}$ once moving in $x$ from $-L$ to $0$.

At the four contact line points we apply analogous procedures to those used in sections 3-4 involving matching of the solutions and of their first derivatives. We omit the details and just state the resulting system of the matching conditions:
\be
\left\{\begin{array}{l}
	-\frac{\lambda_2^0}{2(\sigma+1)}(s_3+\widetilde{x}_{c1})^2+\widetilde{C}_1=0,\hspace{1.8cm}
	-\frac{\lambda_2^0}{\sigma+1}(-s_3-\widetilde{x}_{c1})=\sqrt{\frac{2|\phi(1)|}{\sigma+1}},\\[1.5ex]
	-\frac{\lambda_2^0}{2(\sigma+1)}(s_2+\widetilde{x}_{c1})^2+\widetilde{C}_1=C_4,\hspace{1.6cm}
	-\frac{\lambda_2^0}{\sigma+1}(-s_2-\widetilde{x}_{c1})=C_5,\\[1.5ex]
	 \frac{\lambda_1^0-\lambda_2^0}{2\sigma}(s_2+x_{c1})^2+C_1=C_4,\hspace{2.cm}
	 \frac{\lambda_2^0-(\sigma+1)\lambda_1^0}{2\sigma}(s_2+x_c)^2+C=0,\\[1.5ex]
	 \frac{\lambda_1^0-\lambda_2^0}{\sigma}(-s_2-x_{c1})=-\sqrt{\frac{2|\phi(1)|}{\sigma(\sigma+1)}}+C_5,\quad
	 \frac{\lambda_2^0-(\sigma+1)\lambda_1^0}{\sigma}(-s_2-x_c)=\sqrt{\frac{2(\sigma+1)|\phi(1)|}{\sigma}},\\[1.5ex]
	 \frac{\lambda_1^0-\lambda_2^0}{2\sigma}(s_1+x_{c1})^2+C_1=0,\hspace{2.2cm}
	 \frac{\lambda_2^0-(\sigma+1)\lambda_1^0}{2\sigma}(s_1+x_c)^2+C=C_2,\\[1.5ex]
	 \frac{\lambda_1^0-\lambda_2^0}{\sigma}(-s_1-x_{c1})=-\sqrt{\frac{2|\phi(1)|}{\sigma}},\hspace{1.6cm}
	 \frac{\lambda_2^0-(\sigma+1)\lambda_1^0}{\sigma}(-s_1-x_c)=\sqrt{\frac{2|\phi(1)|}{\sigma}}+C_3,\\[1.5ex]	 
	 -\frac{\lambda_1^0}{2}(s_1+\widetilde{x}_c)^2+\widetilde{C}=C_2,\hspace{2.6cm}
	 -\lambda_1^0(-s_1-\widetilde{x}_c)=C_3,\\[1.5ex]
	-\frac{\lambda_1^0}{2}(s+\widetilde{x}_c)^2+\widetilde{C}=0,\hspace{3.cm}
	-\lambda_1^0(-s-\widetilde{x}_c)=-\sqrt{2|\phi(1)|}.    
\end{array}\right.
\lb{2m0m13sol}
\ee
Additionally, we fix $\widetilde{C}=h^m$ and $\widetilde{C}_1=h_1^m$. System \rf{2m0m13sol} has $16$ equations with $16$ unknowns: $\lambda_1^0,\,\lambda_2^0,\,s,\,s_1,\,s_2,\,s_3,\,x_c,\,x_{c1},\,\widetilde{x}_c,\,\widetilde{x}_{c1}$ and $C,\,C_1,\,C_4,\,C_5,\,C_2,\,C_3$. In fact, due to a translation invariance (within interval $(-L,\,0)$) possessed by this solution one of the four parameters $x_c,\,x_{c1},\,\widetilde{x}_c,\,\widetilde{x}_{c1}$ can be chosen as a free {\it (shift)} parameter, which makes \rf{2m0m13sol} an overdetermined algebraic system.

From the first and the last two equations in \rf{2m0m13sol} one easily obtains expressions
\be
\lambda_1^0=\tfrac{|\phi(1)|}{h^m}\quad\text{and}\quad\lambda_2^0=\tfrac{|\phi(1)|}{h_1^m}.
\lb{lambda_grzz}
\ee
Similarly, dividing the $4^{\mathrm{th}}$ over the $3^{\mathrm{rd}}$ equation in \rf{2m0m13sol} and using \rf{lambda_grzz} yields 
\be
C_5^2=\tfrac{2}{\sigma+1}(|\phi(1)|-C_4\lambda_2^0),
\lb{Ms6}
\ee
while dividing  the $9^{\mathrm{th}}$ over the $11^{\mathrm{th}}$ equation in \rf{2m0m13sol} gives $C_1=\tfrac{|\phi(1)|}{\lambda_2^0-\lambda_1^0}$.
In turn, dividing the $7^{\mathrm{th}}$ over the $5^{\mathrm{th}}$ equation in \rf{2m0m13sol} and using the last two expressions yields
\be
C_5=\tfrac{C_4(\lambda_2^0-(\sigma+1)\lambda_1^0)}{\sqrt{2\sigma(\sigma+1)|\phi(1)|}}.
\lb{C3}
\ee
Comparing expressions \rf{Ms6}--\rf{C3} one deduces a quadratic equation for $C_4$:
\be
\tfrac{|C_4|^2(\lambda_2^0-(\sigma+1)\lambda_1^0)^2}{2\sigma|\phi(1)|}=2(|\phi(1)|-C_4\lambda_2^0).
\lb{C2}
\ee
Further, substituting \rf{lambda_grzz} into \rf{C2} yields a unique positive solution to the latter:
\be
C_4=\frac{-2\sigma\overline{h}+\sqrt{4\sigma^2\overline{h}^2+4\sigma(\overline{h}-\sigma-1)^2}}{(\overline{h}-\sigma-1)^2}\cdot h^m\quad\text{with}\quad\overline{h}=\tfrac{h^m}{h_1^m}.
\lb{C2_expr}
\ee
Proceeding similarly, after deducing expression $C=\tfrac{(\sigma+1)|\phi(1)|}{\lambda_2^0-(\sigma+1)\lambda_1^0}$
from the $8^{\mathrm{th}}$ and $6^{\mathrm{th}}$ equations in \rf{2m0m13sol} and manipulating with $10^{\mathrm{th}}-14^{\mathrm{th}}$ ones yields relations
\sbea
\lb{C5}
C_3&=&\tfrac{C_2(\lambda_2^0-\lambda_1^0)}{\sqrt{2\sigma|\phi(1)|}},\\
\tfrac{|C_2|^2(\lambda_2^0-\lambda_1^0)^2}{2\sigma|\phi(1)|}&=&2(|\phi(1)|-C_2\lambda_1^0).
\lb{C4_1}
\seea
Next, substituting \rf{lambda_grzz} into \rf{C4_1} one finds a unique positive solution to the latter:
\be
C_2=\frac{-2\sigma+\sqrt{4\sigma^2+4\sigma(\overline{h}-1)^2}}{(\overline{h}-1)^2}\cdot h^m\quad\text{with}\quad\overline{h}=\tfrac{h^m}{h_1^m}.
\lb{C4_expr}
\ee

Finally, letting $\widetilde{x}_c$ be a free {\it (shift)} parameter from the linear part of system \rf{2m0m13sol}, i.e. by manipulating with the $2^{\mathrm{nd}}$,  $4^{\mathrm{th}}$, $7^{\mathrm{th}}-8^{\mathrm{th}}$, $11^{\mathrm{th}}-12^{\mathrm{th}}$, $14^{\mathrm{th}}$, and $16^{\mathrm{th}}$ equations in \rf{2m0m13sol}, the following expressions for the positions are obtained inductively:
\bea
s&=&-\tfrac{\sqrt{2|\phi(1)|}}{\lambda_1^0}-\widetilde{x}_c,\quad\quad\quad\quad s_1=\tfrac{C_3}{\lambda_1^0}-\widetilde{x}_c,\nonumber\\[1.ex]
x_c&=&-\tfrac{\sqrt{2\sigma|\phi(1)|}+\sigma C_3}{\lambda_2^0-(\sigma+1)\lambda_1^0}-s_1,\quad\quad x_{c1}=\tfrac{\sqrt{2\sigma|\phi(1)|}}{\lambda_1^0-\lambda_2^0}-s_1,\nonumber\\[1.ex]
s_2&=&-\tfrac{\sqrt{2(\sigma+1)\sigma|\phi(1)|}}{\lambda_2^0-(\sigma+1)\lambda_1^0}-x_c=\tfrac{\sqrt{2\sigma|\phi(1)|/(\sigma+1)}-\sigma C_5}{\lambda_1^0-\lambda_2^0}-x_{c1},\nonumber\\[1.ex]
\widetilde{x}_{c1}&=&\tfrac{(\sigma+1)C_5}{\lambda_2^0}-s_2,\quad\quad\quad\quad\quad s_3=\tfrac{\sqrt{2(\sigma+1)|\phi(1)|}}{\lambda_2^0}-\widetilde{x}_{c1}.
\lb{sx_expr_2}
\eea
Note that two different expressions for $s_2$ in \rf{sx_expr_2} are consistent with each other. They can be rewritten as
\bes
\tfrac{\sqrt{2\sigma|\phi(1)|}(1-1/\sqrt{\sigma+1})+\sigma C_5}{\lambda_1^0-\lambda_2^0}=
\tfrac{\sqrt{2\sigma|\phi(1)|}(\sqrt{\sigma+1}-1)-\sigma C_3}{\lambda_2^0-(\sigma+1)\lambda_1^0}.
\ees
One can show validity of the last relation using expressions \rf{C3}--\rf{C2} and \rf{C5}--\rf{C4_1} arguing similarly as in the next following argument that crosschecks satisfaction of the quadratic equations in system \rf{2m0m13sol}.

By manipulating with the 3$^{\mathrm{rd}}$, 5$^{\mathrm{th}}$, and 9$^{\mathrm{th}}$ equations in \rf{2m0m13sol} one deduces a combined quadratic relation
\be
\tfrac{\lambda_1^0-\lambda_2^0}{2\sigma}\left[(s_2+x_{c1})^2-(s_1+x_{c1})^2\right]=h_1^m-\tfrac{\lambda_2^0}{2(\sigma+1)}(s_2+\widetilde{x}_{c1})^2.
\lb{qr1}
\ee
Using \rf{sx_expr_2} it can be rewritten as
\bes
\tfrac{1}{2\sigma(\lambda_1^0-\lambda_2^0)}\left[(\sqrt{2\sigma|\phi(1)|/(\sigma+1)}-\sigma C_5)^2-2\sigma|\phi(1)|\right]=h_1^m-\tfrac{\sigma+1}{2\lambda_2^0}C_5^2.
\ees
Next, multiplying both sides by $2\sigma(\lambda_1^0-\lambda_2^0)\lambda_2^0$ along with bracket opening yields
\bes
C_5^2(\lambda_2^0-(\sigma+1)\lambda_1^0)-2\sqrt{\tfrac{2\sigma|\phi(1)|}{\sigma+1}}\lambda_2^0C_5-2\lambda_2^0\left(\tfrac{\sigma|\phi(1)|}{\sigma+1}+(\lambda_1^0-\lambda_2^0)h_1^m\right)=0.
\ees
Finally, substituting expression \rf{C3} for $C_5$ and subsequently using \rf{C2} reduces the last relation simply to the expression of $\lambda_2^0$  in \rf{lambda_grzz}.

Similarly, using \rf{C5}--\rf{C4_1} and  \rf{sx_expr_2} one verifies another combined quadratic relation
\be
\tfrac{\lambda_2^0-(\sigma+1)\lambda_1^0}{2\sigma}\left[(s_1+x_c)^2-(s_2+x_c)^2\right]=h^m-\tfrac{\lambda_1^0}{2}(s_1+\widetilde{x}_c)^2,
\lb{qr2}
\ee
which is deduced from the 13$^{\mathrm{rd}}$, 10$^{\mathrm{th}}$, and 6$^{\mathrm{th}}$ equations in \rf{2m0m13sol}.

In summary, formulae \rf{lambda_grzz},   \rf{C3}, \rf{C2_expr}, \rf{C5} and \rf{C4_expr}--\rf{sx_expr_2} provide a complete information about the leading order (as $\eps\go0$) profile of the derived {\it 2-side sessile zig-zag} solution, typical shapes of which are shown in Fig.10. They also cover the two limiting cases $\lambda_2^0=\lambda_1^0$ and $\lambda_2^0=(1+\sigma)\lambda_1^0$ in system \rf{2m0m13sol} corresponding to 
%%%%%%%%%%%%%%%%%%%%%%%%%%%%%%%%%%%%%%%%%%%%%%%%%%%%%%%%%%%%%%%%%%%%%%%%%%%%%%%%%%%%%%%%%%%%%%%%%%%%%%%%%%%%
\begin{figure}[H] 
	\centering
	\vspace{-2.9cm}
	\hspace{-.4cm}\includegraphics[width=.5\textwidth]{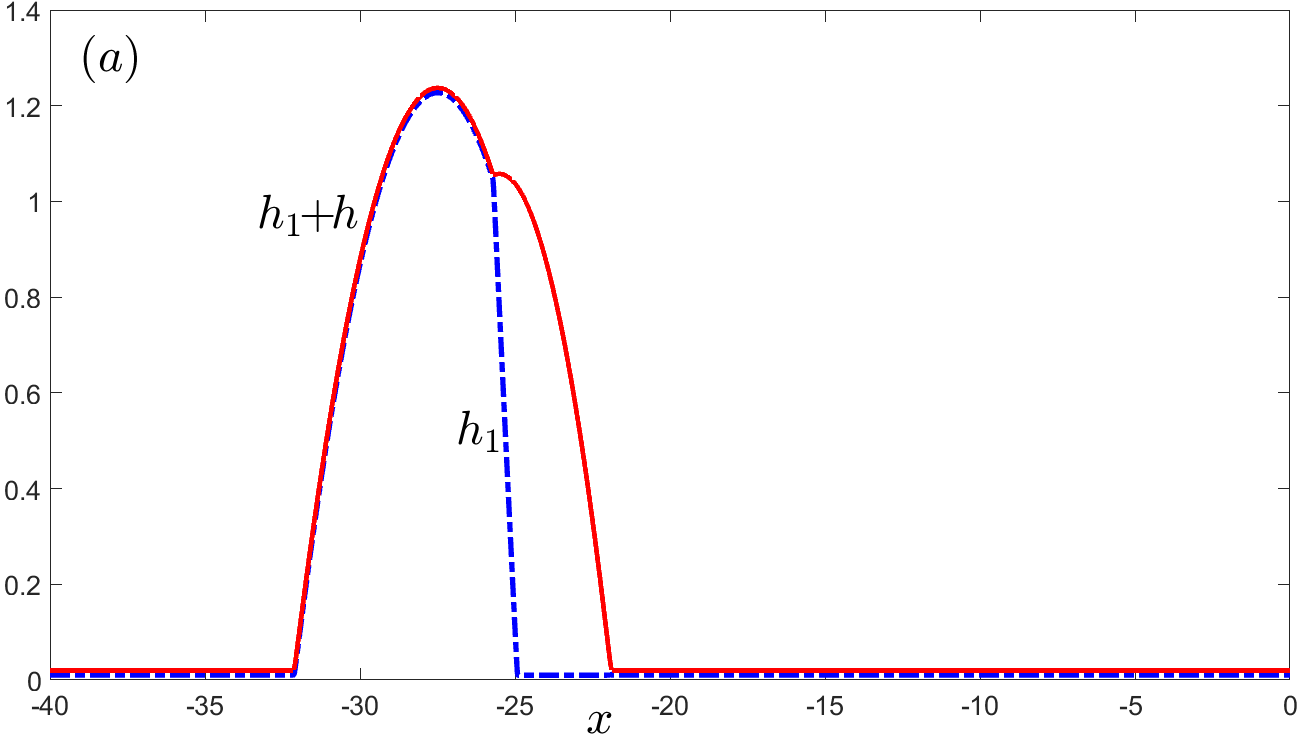}  
	\hspace{.2cm}\includegraphics[width=.5\textwidth]{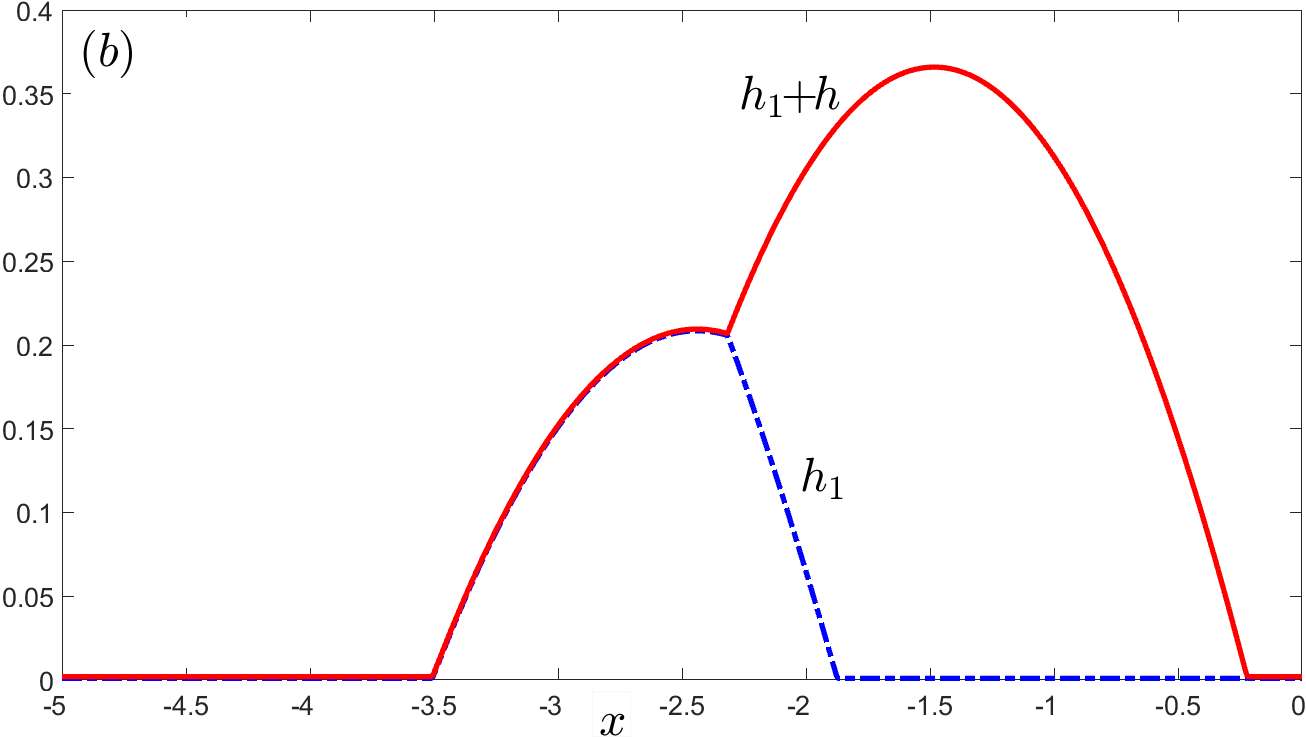}  
	
	\vspace{.2cm}
	\hspace{-.4cm}\includegraphics[width=.5\textwidth]{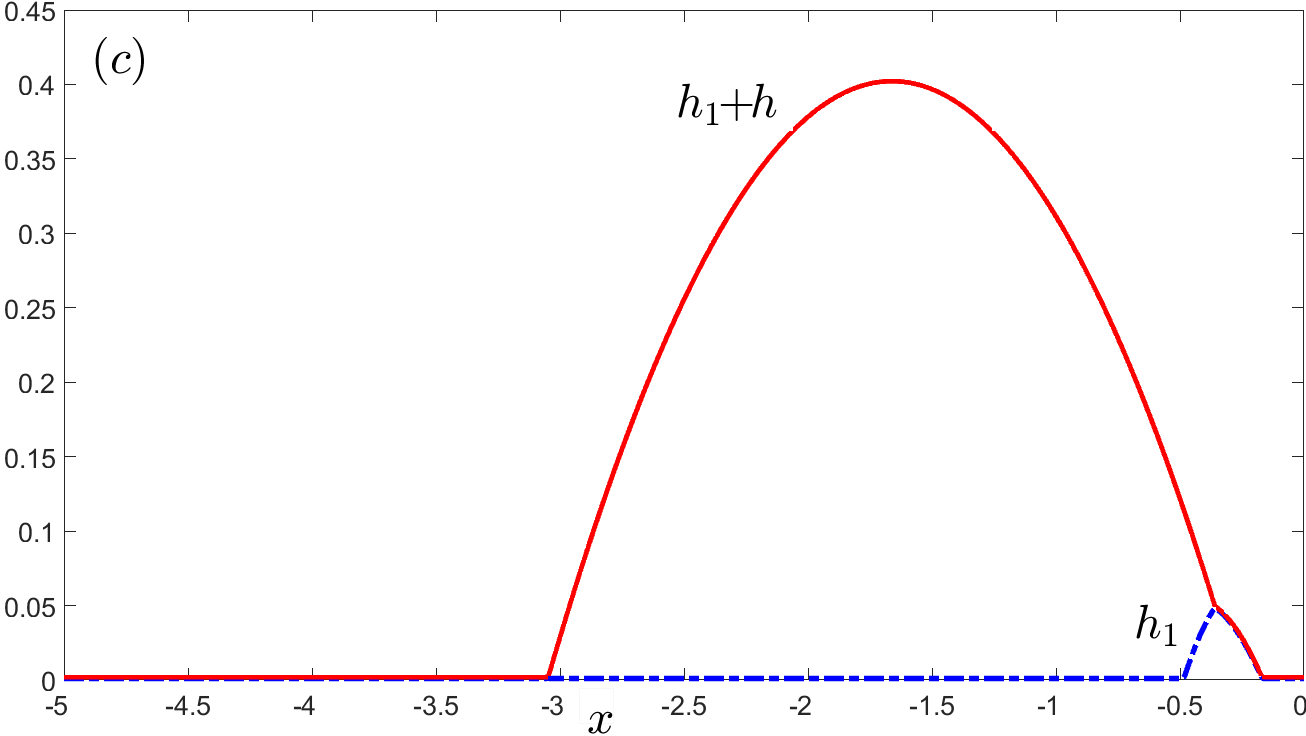}  
	\hspace{.2cm}\includegraphics[width=.5\textwidth]{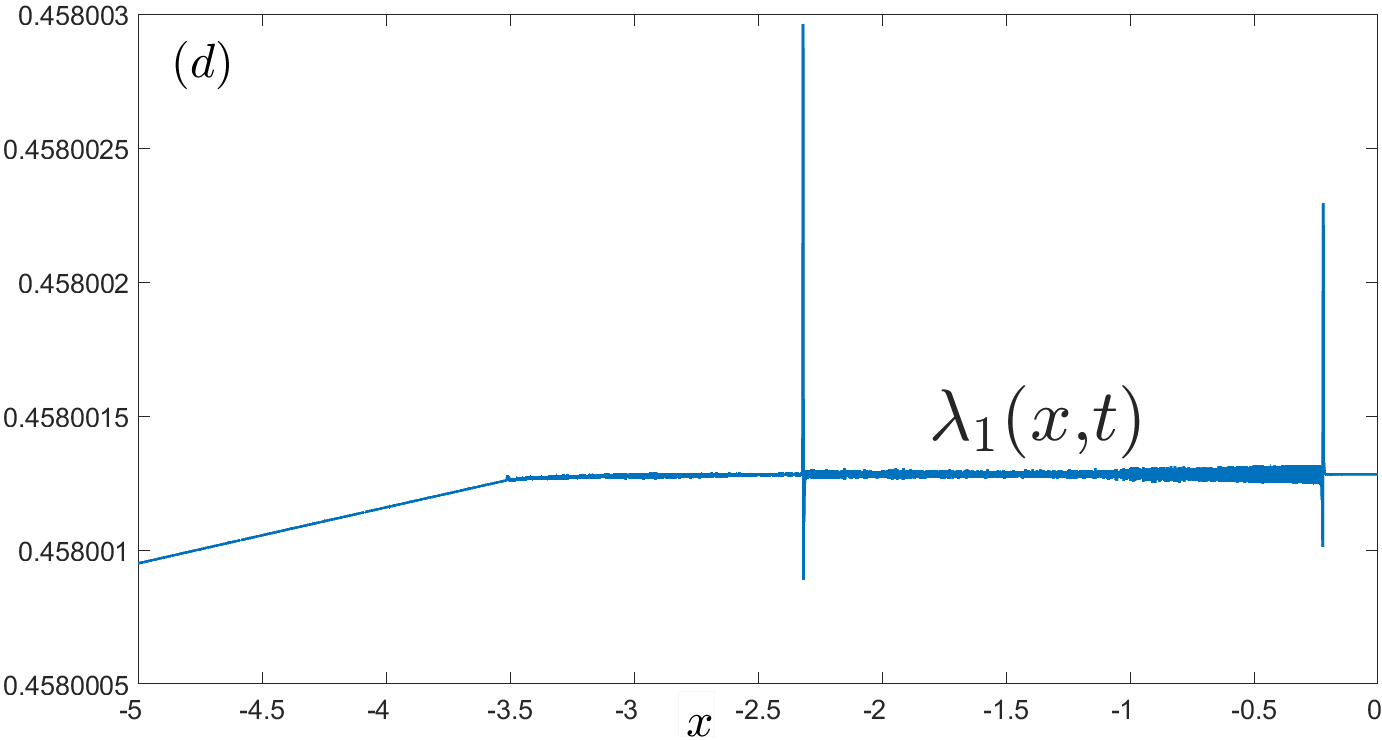}  
	
	\caption{\small Numerical 2-side sessile zig-zag solution to system \rf{BS}--\rf{BC} for\hspace{2.cm}$\eps\!=\!0.001$: {\bf (a)} $\sigma\!=\!0.2,\,L\!=\!40,\,\lambda_1\!=\!0.1606,\,\lambda_2\!=\!0.1369,\,\tilde{x}_c\!=\!-25.51$; {\bf (b)} $\sigma\!=\!1.2$,\hspace{2.cm}$L\!=\!5,\,\lambda_1\!=\!0.4580,\, \lambda_2\!=\!0.8032,\,\tilde{x}_c\!=\!-1.5$; {\bf (c)} $\sigma\!=\!1.2,\,L\!=\!5,\,\lambda_1\!=\!0.4165$,\hspace{2.cm}$\lambda_2\!=\!3.1906,\,\tilde{x}_c\!=\!-3.34$ ($x$-variable inverted); {\bf (d)} numerically observed pressure profile $\lambda_1(x,\,t)\!=\!\Pi_\eps(h(x,\,t))-h_{xx}(x,\,t)-h_{1,xx}(x,\,t)$ for solution in {\bf (b)}.}	
\end{figure}
%%%%%%%%%%%%%%%%%%%%%%%%%%%%%%%%%%%%%%%%%%%%%%%%%%%%%%%%%%%%%%%%%%%%%%%%%%%%%%%%%%%%%%%%%%%%%%%%%%%%%%%%%%%%
\noindent $\overline{h}=1$ or $\overline{h}=\sigma+1$ values, respectively, in \rf{C2_expr} and \rf{C4_expr}. In these cases, profile of $h_1^0(x)$ or $h^0(x)$ in the bulk region $(-s_2,\,-s_1)$ is given by a linear segment. The corresponding limiting formulae for positions $s,\,s_1,\,s_2,\,s_3,\,x_c,\,x_{c1},\widetilde{x}_{c1}$ can be derived then from \rf{sx_expr_2}.

Note that in contrast to stationary solutions to system \rf{BS}--\rf{BC} of sections $3-4$,  {\it 2-side sessile zig-zag} is not dynamically stable one. In particular, that is indicated 
by the spikes observed in the pressures profile of Fig.10 {\bf (d)} around $x=-s_2$ and $x=-s$ contact lines. Nevertheless, numerical simulations of system \rf{BS}-\rf{BC} show that solutions of  Fig.10 {\bf (a)}--{\bf (c)} translate extremely slow to one of the interval ends $x=-L$ or $x=0$ without changing their leading order shape. Such behavior was already observed for symmetric drop solutions in one-layer lubrication equations~\cite{BGW01} with their {\it translational instability} originating from a single exponentially small eigenvalue in the spectra~\cite{KRW12} of the corresponding linearized operators. Therefore, {\it 2-side sessile zig-zag} is a {\it weakly unstable} stationary solution to system \rf{BS}--\rf{BC}.

The constraints on {\it 2-side sessile zig-zag} are deduced from the length ones 
\be
L>s_3>s_2>s_1>s>0
\lb{LC2}
\ee
as follows. From the expressions for $x_c$ and $s_2$ in \rf{sx_expr_2} it follows that the constraint $s_2>s_1$ is satisfied if
\bes
\tfrac{\sigma C_3+\sqrt{2\sigma|\phi(1)|}(1-\sqrt{\sigma+1})}{\lambda_2^0-(\sigma+1)\lambda_1^0}>0.
\ees
If $\lambda_2^0-(\sigma+1)\lambda_1^0>0$, i.e. when $\overline{h}>\sigma+1$,  the last inequality reduces to
\be
\sigma C_3>\sqrt{2\sigma|\phi(1)|}(\sqrt{\sigma+1}-1).
\lb{Ms7}
\ee
Substituting \rf{C5} and subsequently expressions \rf{C4_expr} and \rf{lambda_grzz} into \rf{Ms7} gives
\bes
-\sigma+\sqrt{\sigma^2\overline{h}^2+\sigma(\overline{h}-1)^2}>(\overline{h}-1)(\sqrt{\sigma+1}-1),
\ees
which indeed holds if $\overline{h}>\sigma+1$. In turn, if $\lambda_2^0-(\sigma+1)\lambda_1^0<0$, the inequality in \rf{Ms7} is reversed and its satisfaction for $\overline{h}<\sigma+1$ follows by a similar argument. Therefore, we conclude that $s_2>s_1$  holds for all positive $\overline{h}$. 

Next, using expressions in \rf{sx_expr_2} the constraints $s_1>s$  and $s_3>s_2$ can be rewritten as
\bes
C_3+2\sqrt{|\phi(1)|}>0\quad\text{and}\quad \sqrt{2(\sigma+1)|\phi(1)|}-(\sigma+1)C_5>0,
\ees
respectively. Analogously, one shows that these two inequalities hold for all $\overline{h}>0$ using relations \rf{C5} or \rf{C3} with subsequent substitution of  expressions \rf{C4_expr} or \rf{C2_expr}, respectively, along with \rf{lambda_grzz}.

Finally, the first and the last inequalities in \rf{LC2} result in the constraints on the minimal interval length $L$ and the maximal free shift parameter $\widetilde{x}_c$, respectively. Namely, the two conditions imposed on {\it 2-side sessile zig-zag} can be written as
\be
L>\tfrac{\sqrt{2(\sigma+1)|\phi(1)|}}{\lambda_2^0}-\widetilde{x}_{c1}\quad\text{and}\quad \widetilde{x}_c<-\tfrac{\sqrt{2|\phi(1)|}}{\lambda_1^0}.
\lb{2sd_cons}
\ee

Additionally, we note that the maxima of $h^0(x)$ and $h_1^0(x)$ in $(-L,\,0)$ are given by values $h^m$ and $h_1^m$, respectively, only when $\widetilde{x}_{c1}\in[-s_3,\,s_2]$ and $\widetilde{x}_c\in[-s_1,\,s]$ hold.
Using expressions \rf{sx_expr_2} and \rf{C3}, \rf{C5} one checks that the latter two conditions are equivalent to 
$\lambda_2^0\le(\sigma+1)\lambda_1^0$ and $\lambda_2^0\ge\lambda_1^0$ ones, respectively, which in turn by \rf{lambda_grzz} reduce to $\overline{h}\le(\sigma+1)$ and $\overline{h}\ge1$ (as before $\overline{h}=h^m/h_1^m$).
Therefore, the following combined formula holds
\be
\displaystyle\max_{x\in(-L,\,0)}(h_1^0(x),\,h^0(x))=\left\{\begin{array}{ll}
	(h_1^m,\,C_2)& \text{if}\ \ \overline{h}\le1\ \text{(cf. Fig.10 {\bf (a)})},\\[1.5ex]
	(h_1^m,\,h^m)& \text{if}\ \ 1\le\overline{h}\le\sigma+1\ \text{(cf. Fig.10 {\bf (b)})},\\[1.5ex]
	(C_4,\,h^m)& \text{if}\ \ \overline{h}\ge\sigma+1\ \text{(cf. Fig.10 {\bf (c)})},\end{array}\right.
\lb{hm_h1m_2sgzz}
\ee 
with $C_4$ and $C_2$ given by \rf{C2_expr} and \rf{C4_expr}, respectively. 

In order to invert formula \rf{hm_h1m_2sgzz}, i.e. to find out how values $h_1^m$ and $h^m$ depend on the given maxima of  $h_1^0(x)$ and $h^0(x)$ let us denote
\be
\displaystyle\widetilde{h}=\max_{x\in(-L,\,0)}h^0(x)/\max_{x\in(-L,\,0)}h_1^0(x).
\lb{ht}
\ee
In the case $\overline{h}\le1$, inverting formula \rf{C4_expr} using \rf{hm_h1m_2sgzz} yields:
\bes
\overline{h}=\tfrac{\widetilde{h}}{\widetilde{h}-2\sigma\pm\sqrt{4\sigma(\sigma+1-\widetilde{h})}}\le 1.
\ees
Analyzing the last inequality one observes that only $+$ sign is possible there leading, subsequently, to condition $\widetilde{h}\le 1$. Arguing similarly for the case $\overline{h}\ge\sigma+1$ by inverting now formula \rf{C2_expr} one obtains an equivalent condition $\widetilde{h}\ge\sigma+1$. Combining all together yields the following inverse formula to \rf{hm_h1m_2sgzz}:
\be
\displaystyle (h_1^m,\,h_m)=\left\{\begin{array}{ll}
	\left(\max h_1^0(x),\,\tfrac{\max h^0(x)}{\widetilde{h}-2\sigma+\sqrt{4\sigma(\sigma+1-\widetilde{h})}}\right)& \text{if}\ \ \widetilde{h}\le1,\\[2.5ex]
	\left(\max h_1^0(x),\,\max h^0(x)\right)& \text{if}\ \ 1\le \widetilde{h}\le\sigma+1,\\[2.5ex]
	\left(\tfrac{\max h^0(x)}{\sigma+1-2\sigma \widetilde{h}+\sqrt{4\sigma(\sigma+1) \widetilde{h}(\widetilde{h}-1)}},\,\max h^0(x)\right)& \text{if}\ \ \widetilde{h}\ge\sigma+1.\end{array}\right.
\lb{hm_h1m_2sgzz_inv}
\ee

%%%%%%%%%%%%%%%%%%%%%%%%%%%%%%%%%%%%%%%%%%%%%%%%%%%%%%%%%%%%%%%%%
\section{Three-CL solutions}
%%%%%%%%%%%%%%%%%%%%%%%%%%%%%%%%%%%%%%%%%%%%%%%%%%%%%%%%%%%%%%%%%
In this section, we classify and derive the leading order profiles of the solutions to system \rf{SSa}--\rf{SSc} having three contact lines. Up to a possible inversion of $x$-variable we find two solution types denoted below as {\it $h_1$-sessile zig-zag} and {\it $h$-sessile zig-zag}. Due to the fact that both solution types can be viewed as opposite parts of the {\it 2-side sessile zig-zag} described in section 5 derivation of their leading order profiles proceeds analogously.\\[.5ex]

\underline{\bf $h_1$-sessile zig-zag:}\hspace{.5cm}The solution is obtained by a combined matching of the following leading order profiles: general {\bf Type I} and {\bf Type II} bulk (with $x_{c1}=-L$) solutions to {\bf Type II} CL ones centered around $x=-s_2$ ($1^{\mathrm{st}}$ CL); general {\bf Type I} and {\bf Type III} bulk solutions to {\bf Type I} CL ones centered around $x=-s_1$ ($2^{\mathrm{nd}}$ CL); as well as {\bf UTF} and general {\bf Type III} bulk solutions to {\bf Type IV} CL centered around $x=-s$ ($3^{\mathrm{rd}}$ CL). Schematically this matching chain can be encrypted as $\bf{IIB-IICL-IB-ICL-IIIB-IVCL}$--$\bf{-UTF}$ once moving in $x$ from $-L$ to $0$.
Applying the matching procedure at the three contact points one obtains the following system:
\be
\left\{\begin{array}{l}
	-\frac{\lambda_2^0}{2(\sigma+1)}(L-s_2)^2+\widetilde{C}_1=C_4,\hspace{1.8cm}
	-\frac{\lambda_2^0}{\sigma+1}(L-s_2)=C_5,\\[1.5ex]
	\frac{\lambda_1^0-\lambda_2^0}{2\sigma}(s_2+x_{c1})^2+C_1=C_4,\hspace{2.cm}
	\frac{\lambda_2^0-(\sigma+1)\lambda_1^0}{2\sigma}(s_2+x_c)^2+C=0,\\[1.5ex]
	\frac{\lambda_1^0-\lambda_2^0}{\sigma}(-s_2-x_{c1})=-\sqrt{\frac{2|\phi(1)|}{\sigma(\sigma+1)}}+C_5,\quad
	\frac{\lambda_2^0-(\sigma+1)\lambda_1^0}{\sigma}(-s_2-x_c)=\sqrt{\frac{2(\sigma+1)|\phi(1)|}{\sigma}},\\[1.5ex]
	\frac{\lambda_1^0-\lambda_2^0}{2\sigma}(s_1+x_{c1})^2+C_1=0,\hspace{2.2cm}
	\frac{\lambda_2^0-(\sigma+1)\lambda_1^0}{2\sigma}(s_1+x_c)^2+C=C_2,\\[1.5ex]
	\frac{\lambda_1^0-\lambda_2^0}{\sigma}(-s_1-x_{c1})=-\sqrt{\frac{2|\phi(1)|}{\sigma}},\hspace{1.5cm}
	\frac{\lambda_2^0-(\sigma+1)\lambda_1^0}{\sigma}(-s_1-x_c)=\sqrt{\frac{2|\phi(1)|}{\sigma}}+C_3,\\[1.5ex]	 
	-\frac{\lambda_1^0}{2}(s_1+\widetilde{x}_c)^2+\widetilde{C}=C_2,\hspace{2.6cm}
	-\lambda_1^0(-s_1-\widetilde{x}_c)=C_3,\\[1.5ex]
	-\frac{\lambda_1^0}{2}(s+\widetilde{x}_c)^2+\widetilde{C}=0,\hspace{2.9cm}
	-\lambda_1^0(-s-\widetilde{x}_c)=-\sqrt{2|\phi(1)|}.    
\end{array}\right.
\lb{0m13sol}
\ee
Additionally, we fix $\widetilde{C}=h^m$ and $\widetilde{C}_1=h_1^m$. System \rf{0m13sol} has $14$ equations with $14$ unknowns: $\lambda_1^0,\,\lambda_2^0,\,s,\,s_1,\,s_2,\,x_c,\,x_{c1},\,\widetilde{x}_c$ and $C,\,C_1,\,C_4,\,C_5,\,C_2,\,C_3$.

From the last two equations in \rf{2m0m13sol} one obtains
\be
\lambda_1^0=\tfrac{|\phi(1)|}{h^m}.
\lb{lambda1_h1zz}
\ee
Next, relation \rf{lambda_rel} expanded for the considered solution in powers of $\eps$ yields after retaining only the leading order terms:
\be
\lambda_2^0h_1^m=-\phi(1)\quad\text{and, hence,}\quad\lambda_2^0=\tfrac{|\phi(1)|}{h_1^m}.
\lb{lambda2_h1zz}
\ee
The rest of solving system \rf{0m13sol} proceeds exactly as for \rf{2m0m13sol} in section 5. Namely, from the $1^{\mathrm{st}}-3^{\mathrm{rd}}$, $5^{\mathrm{th}}$, $7^{\mathrm{th}}$, and  $9^{\mathrm{th}}$  equations in \rf{0m13sol} we derive expressions \rf{C3}-\rf{C2_expr}, while from the $4^{\mathrm{th}}$, $6^{\mathrm{th}}$, $8^{\mathrm{th}}$, and $10^{\mathrm{th}}-12^{\mathrm{th}}$ ones \rf{C5}--\rf{C4_expr}. In turn, from the linear part of system \rf{0m13sol}  the following expressions for the positions are obtained inductively:
\bea
s_2&=&\tfrac{(\sigma+1)C_5}{\lambda_2^0}+L,\nonumber\\[1.ex] x_{c1}&=&\tfrac{\sigma}{\lambda_1^0-\lambda_2^0}\left(\sqrt{\tfrac{2|\phi(1)|}{\sigma(\sigma+1)}}-C_5\right)-s_2,\quad\quad
x_c=-\tfrac{\sqrt{2(\sigma+1)\sigma|\phi(1)|}}{\lambda_2^0-(\sigma+1)\lambda_1^0}-s_2,\nonumber\\[1.ex] 
s_1&=&\tfrac{\sqrt{2\sigma|\phi(1)|}}{\lambda_1^0-\lambda_2^0}-x_{c1}=-\tfrac{\sigma}{\lambda_2^0-(\sigma+1)\lambda_1^0}\left(\sqrt{2|\phi(1)|}{\sigma}+C_3\right)-x_c,\nonumber\\[1.ex]
\widetilde{x}_c&=&\tfrac{C_3}{\lambda_1^0}-s_1,\quad\quad\quad\quad\quad\quad\quad\quad\quad\quad\quad s=-\tfrac{\sqrt{2|\phi(1)|}}{\lambda_1^0}-\widetilde{x}_c.
\lb{sx_expr_3}
\eea
Note that consistency of two different expressions for $s_1$ in \rf{sx_expr_3} as well as validity  crosschecks for the quadratic relations in \rf{0m13sol} can be verified using \rf{C3}--\rf{C2} and \rf{C5}--\rf{C4_1} analogously as was done in section 5 for system \rf{2m0m13sol}.

In summary, formulae \rf{lambda1_h1zz}--\rf{lambda2_h1zz},  \rf{C3}, \rf{C2_expr}, \rf{C5}, \rf{C4_expr}, and \rf{sx_expr_3} provide a complete information about the leading order (as $\eps\go0$) profile of the derived {\it $h_1$-sessile zig-zag} solution, typical shapes of which are shown in Fig.11. They also cover one limiting case $\lambda_2^0=\lambda_1^0$ in system \rf{0m13sol}, which correspond to $\overline{h}=1$ value in \rf{C2_expr} and \rf{C4_expr}. In this case, profile of $h_1^0(x)$  in the bulk region $(-s_2,\,-s_1)$ is given by a linear segment. The corresponding limiting formulae for positions $s,\,s_1,\,s_2,\,x_c,\,x_{c1},\widetilde{x}_c$ can be derived then from \rf{sx_expr_3}. The limiting case $\lambda_2^0=(\sigma+1)\lambda_1^0$ is not possible for solutions of \rf{0m13sol}, as shown in the next paragraph it corresponds to the merge $s_2\go L$. 

The constraints on the solution are deduced from the length ones 
\be
L>s_2>s_1>s>0
\lb{LC3}
\ee
as follows. Arguing similar as for \rf{LC2} in section 5 one shows that conditions $s_2>s_1$ and $s_1>s$ hold automatically for all $\bar{h}=h^m/h_1^m>0$. Next, using \rf{sx_expr_3} condition $L>s_2$ transforms to $C_5<0$ and, subsequently, by \rf{C3} to $\lambda_2^0<(\sigma+1)\lambda_1^0$, i.e. using \rf{lambda1_h1zz}--\rf{lambda2_h1zz}  to the constraint $\bar{h}<\sigma+1$ (cf. Fig.11 {\bf (a)}). Finally, using \rf{sx_expr_3} together with \rf{C5} one shows  that condition $s>0$ imposes the following bound:
%%%%%%%%%%%%%%%%%%%%%%%%%%%%%%%%%%%%%%%%%%%%%%%%%%%%%%%%%%%%%%%%%%%%%%%%%%%%%%%%%%%%%%%%%%%%%%%%%%%%%%%%%%%%%%%%
\begin{figure}[H] 
	\centering
	\vspace{-2.9cm}
	\hspace{-.4cm}\includegraphics[width=.5\textwidth]{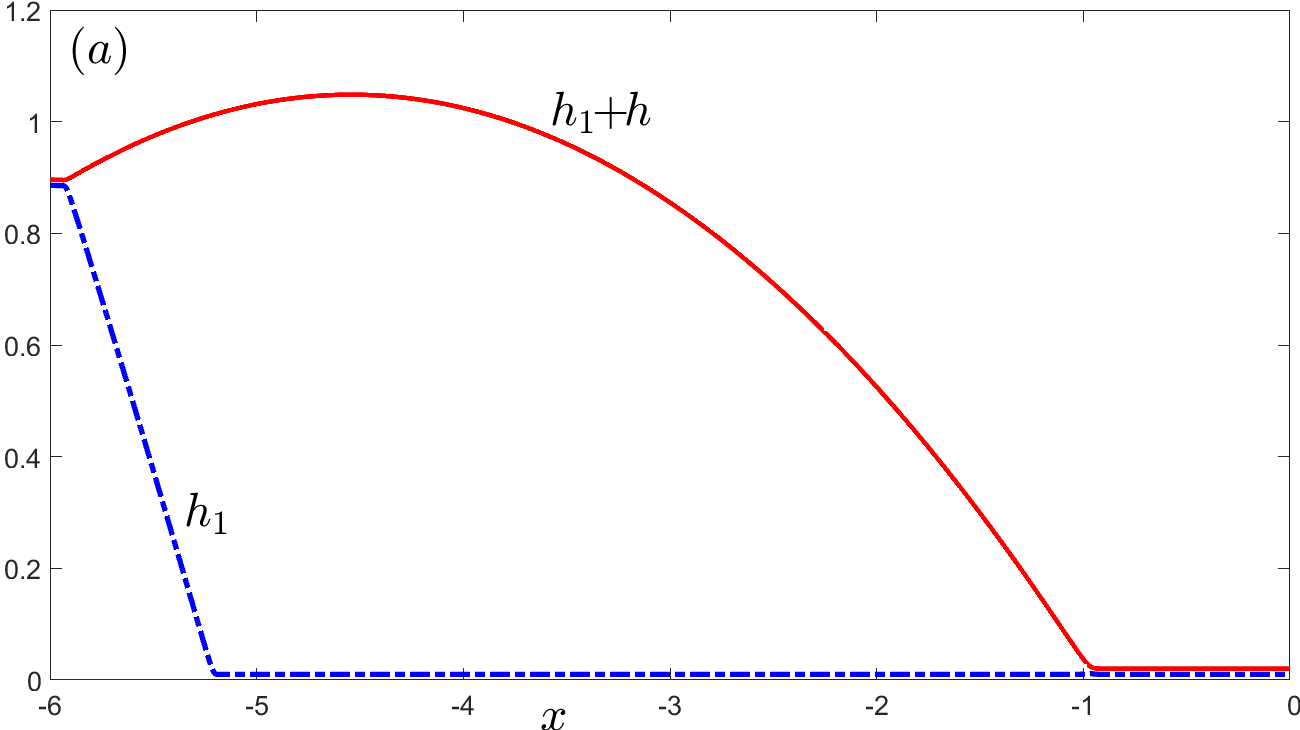}  
	\hspace{.2cm}\includegraphics[width=.5\textwidth]{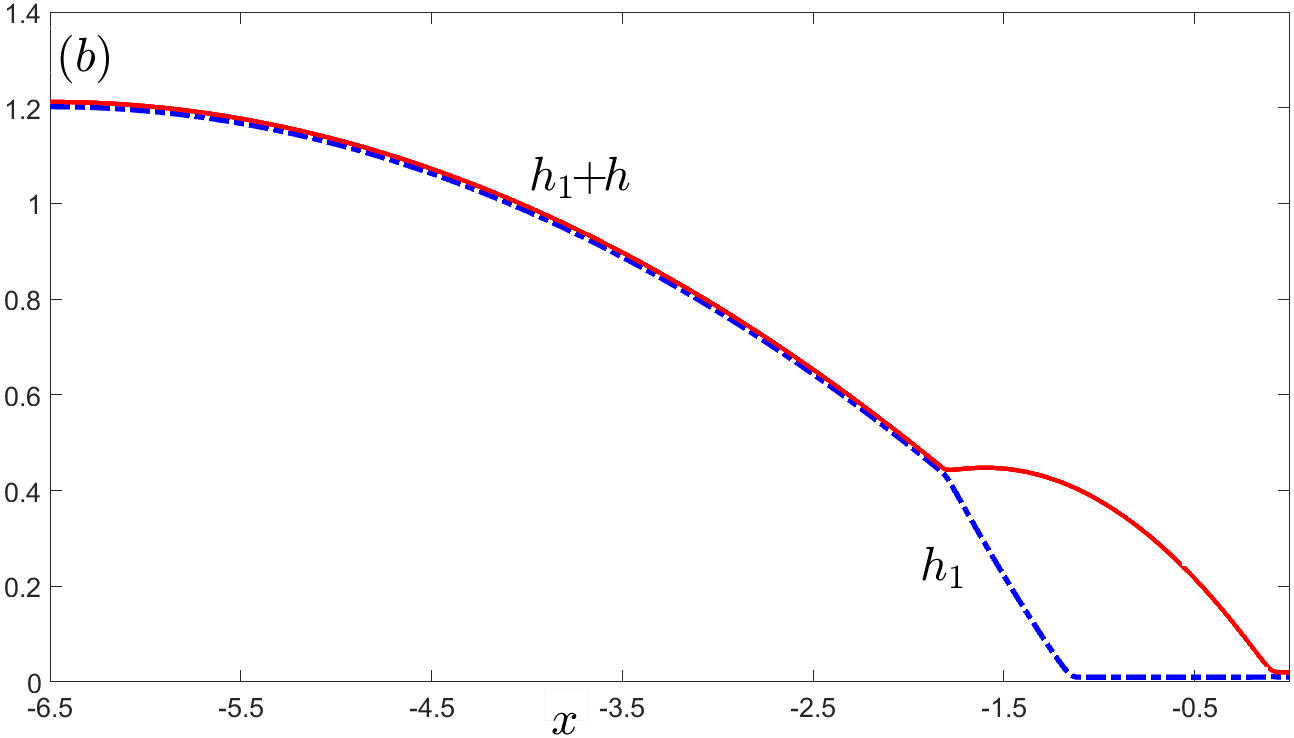}  
	
	\caption{\small Numerical $h_1$-sessile zig-zag stationary solution to system \rf{BS}--\rf{BC} for $\eps\!=\!0.01$: {\bf (a)} $\sigma\!=\!0.2,\,L\!=\!6.0$ and $\lambda_1\!=\!0.16204,\, \lambda_2\!=\!0.19024$; {\bf (b)} $\sigma\!=\!1.0,\,L\!=\!6.5$ and $\lambda_1\!=\!0.38927,\, \lambda_2\!=\!0.13978$.}	
\end{figure}
%%%%%%%%%%%%%%%%%%%%%%%%%%%%%%%%%%%%%%%%%%%%%%%%%%%%%%%%%%%%%%%%%%%%%%%%%%%%%%%%%%%%%%%%%%%%%%%%%%%%%%%%%%%%%%%%
\be
L>\tfrac{2|\phi(1)|(\sigma\lambda_1^0+\sqrt{\sigma}\lambda_2^0)+C_2(\lambda_1^0-\lambda_2^0)^2}{\sqrt{2\sigma|\phi(1)|}\lambda_1^0\lambda_2^0}>0\ \text{(cf. Fig.11 {\bf (b)})}.
\lb{h1zz_minL}
\ee
Additionally, note that $\displaystyle\max_{x\in(-L,\,0)}h^0(x)=h^m$ only when $\widetilde{x}_c\in[-s_1,\,s]$. Using \rf{sx_expr_3}, \rf{C5} one checks that the latter condition is equivalent to $\lambda_2^0\ge\lambda_1^0$ which, in turn, by \rf{lambda_grzz} reduces to $\overline{h}=h^m/h_1^m\ge1$. Therefore, the combined formula
\be
\displaystyle\max_{x\in(-L,\,0)}h^0(x)=\left\{\begin{array}{ll}
	C_2& \text{if}\ \ \overline{h}\le1\ \text{(cf. Fig.11 {\bf (b)})},\\[1.5ex]
	h^m& \text{if}\ \ \overline{h}\ge1\ \text{(cf. Fig.11 {\bf (a)})}\end{array}\right.
\lb{hm_h1gzz}
\ee
holds with $C_2$ given by \rf{C4_expr}. Further, setting $\widetilde{h}$ as in \rf{ht} and proceeding analogously to the derivation of \rf{hm_h1m_2sgzz_inv} in section 5, i.e. by 
inverting formula \rf{C4_expr} for $\overline{h}\le1$, one obtains the inverse formula to \rf{hm_h1gzz}:
\be
\displaystyle h^m=\left\{\begin{array}{ll}
	\tfrac{\max h^0(x)}{\widetilde{h}-2\sigma+\sqrt{4\sigma(\sigma+1-\widetilde{h})}}& \text{if}\ \ \widetilde{h}\le1,\\[2.5ex]
	\max h^0(x).& \text{if}\ \ \widetilde{h}\ge1.\end{array}\right.
\lb{hm_h1gzz_inv}
\ee 
\vspace{2.ex}

\underline{\bf $h$-sessile zig-zag:}\hspace{.3cm}The matching chain for this solution can be encrypted as $\bf{IIIB-ICL-IB-IICL-IIB-IIICL-UTF}$ if moving in $x$ from $-L$ to $0$. 

Applying the matching procedure at the three contact points yields:
\be
\left\{\begin{array}{l}
    -\frac{\lambda_1^0}{2}(L-s_2)^2+\widetilde{C}=C_2,\hspace{2.3cm}
    -\lambda_1^0(L-s_2)=C_3,\\[1.5ex]
    \frac{\lambda_1^0-\lambda_2^0}{2\sigma}(s_2+x_{c1})^2+C_1=0,\hspace{1.9cm}
    \frac{\lambda_2^0-(\sigma+1)\lambda_1^0}{2\sigma}(s_2+x_c)^2+C=C_2,\\[1.5ex]
    \frac{\lambda_1^0-\lambda_2^0}{\sigma}(-s_2-x_{c1})=\sqrt{\frac{2|\phi(1)|}{\sigma}},\hspace{1.6cm}
    \frac{\lambda_2^0-(\sigma+1)\lambda_1^0}{\sigma}(-s_2-x_c)=-\sqrt{\frac{2|\phi(1)|}{\sigma}}+C_3,\\[1.5ex]
    \frac{\lambda_1^0-\lambda_2^0}{2\sigma}(s_1+x_{c1})^2+C_1=C_4,\hspace{1.7cm}
    \frac{\lambda_2^0-(\sigma+1)\lambda_1^0}{2\sigma}(s_1+x_c)^2+C=0,\\[1.5ex]
    \frac{\lambda_1^0-\lambda_2^0}{\sigma}(-s_1-x_{c1})=\sqrt{\frac{2|\phi(1)|}{\sigma(\sigma+1)}}+C_5,\quad
   \frac{\lambda_2^0-(\sigma+1)\lambda_1^0}{\sigma}(-s_1-x_c)=-\sqrt{\frac{2(\sigma+1)|\phi(1)|}{\sigma}},\\[1.5ex]
   -\frac{\lambda_2^0}{2(\sigma+1)}(s_1+\widetilde{x}_{c1})^2+\widetilde{C}_1=C_4,\hspace{1.3cm}
   -\frac{\lambda_2^0}{\sigma+1}(-s_1-\widetilde{x}_{c1})=C_5,\\[1.5ex]
	-\frac{\lambda_2^0}{2(\sigma+1)}(s+\widetilde{x}_{c1})^2+\widetilde{C}_1=0,\hspace{1.6cm}
	-\frac{\lambda_2^0}{\sigma+1}(-s-\widetilde{x}_{c1})=-\sqrt{\frac{2|\phi(1)|}{\sigma+1}}.\end{array}\right.
\lb{1m02sol}
\ee
%%%%%%%%%%%%%%%%%%%%%%%%%%%%%%%%%%%%%%%%%%%%%%%%%%%%%%%%%%%%%%%%%%%%%%%%%%%%%%%%%%%%%%%%%%%%%%%%%%%%%%%%%%%%%%%%
\begin{figure}[H] 
	\centering
	\vspace{-2.9cm}
	\hspace{-.4cm}\includegraphics[width=.5\textwidth]{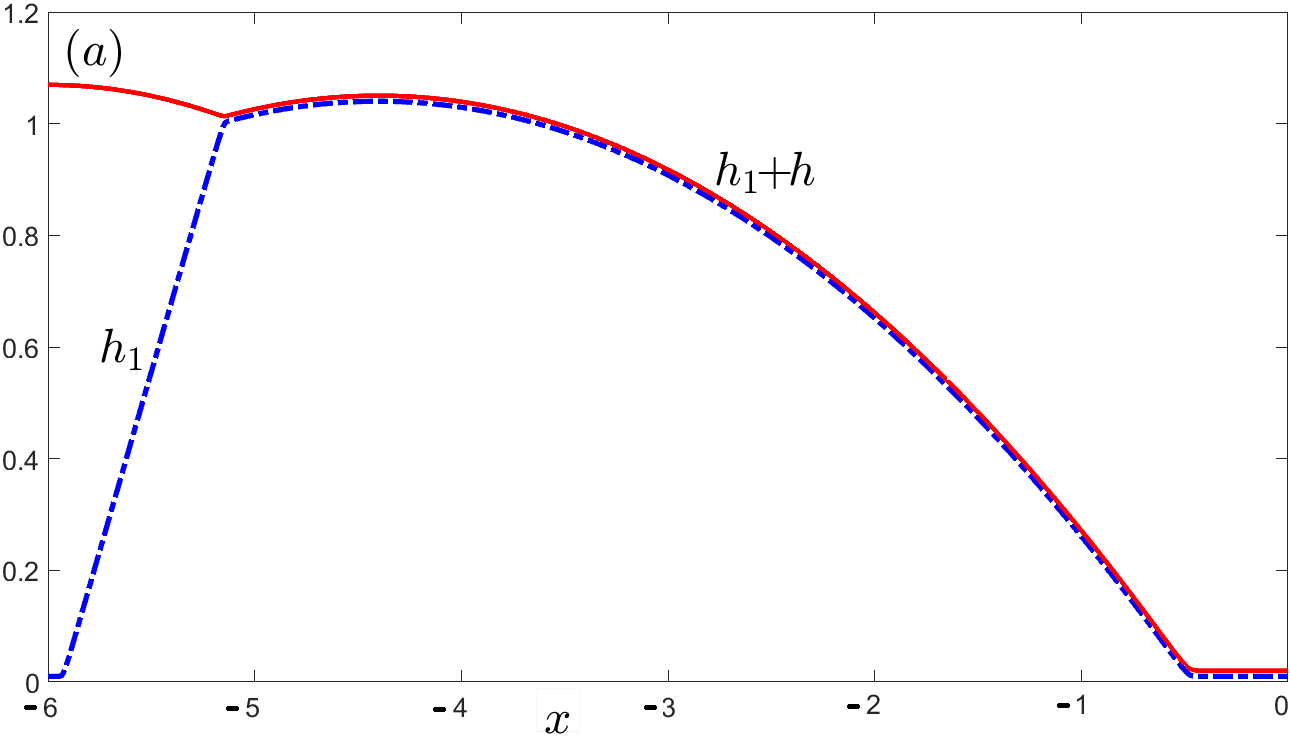}  
	\hspace{.2cm}\includegraphics[width=.5\textwidth]{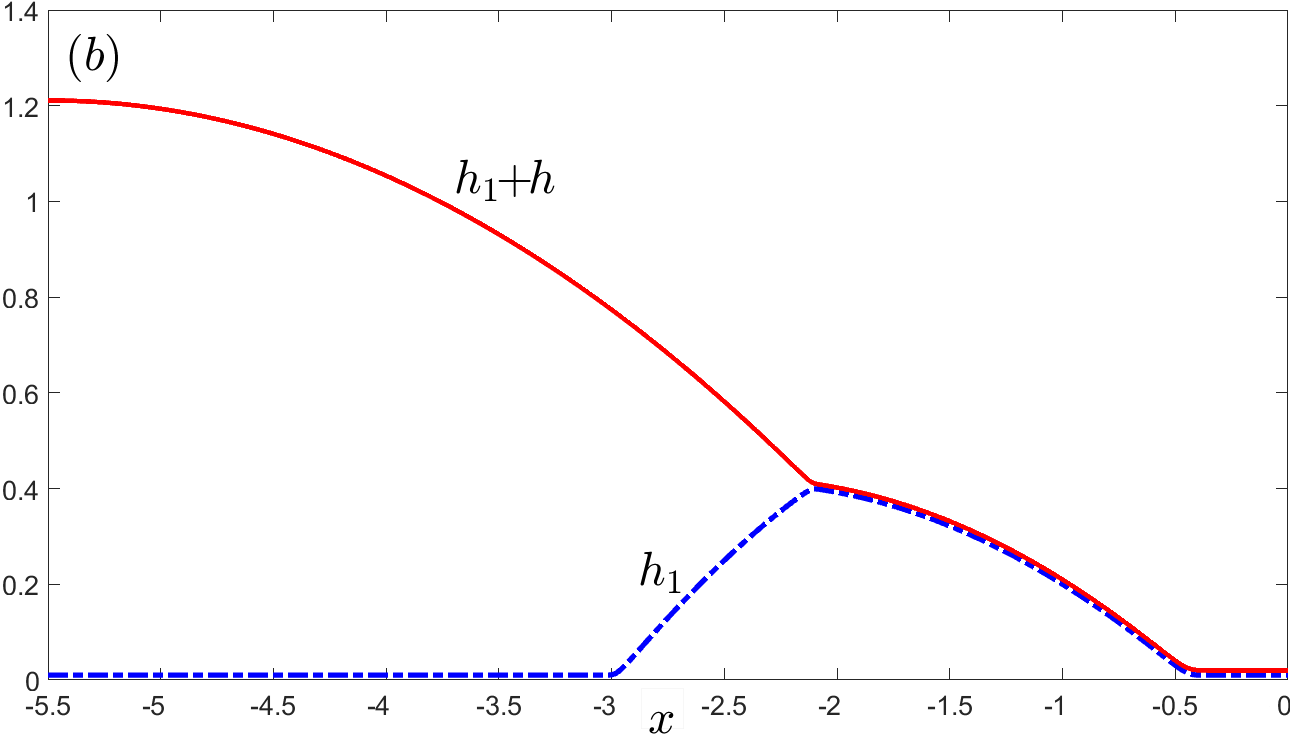}  
	
	\caption{\small Numerical $h$-sessile zig-zag stationary solution to system \rf{BS}--\rf{BC} for $\eps\!=\!0.01$: {\bf (a)} $\sigma\!=\!0.2,\,L\!=\!6.0$ and $\lambda_1\!=\!0.15877,\, \lambda_2\!=\!0.16174$; {\bf (b)} $\sigma\!=\!1.0,\,L\!=\!5.5$ and $\lambda_1\!=\!0.13991,\, \lambda_2\!=\!0.41453$.}	
\end{figure}
%%%%%%%%%%%%%%%%%%%%%%%%%%%%%%%%%%%%%%%%%%%%%%%%%%%%%%%%%%%%%%%%%%%%%%%%%%%%%%%%%%%%%%%%%%%%%%%%%%%%%%%%%%%%%%%%
Additionally, we fix $\widetilde{C}=h^m$ and $\widetilde{C}_1=h_1^m$. System \rf{1m02sol} has $14$ equations with $14$ unknowns: $\lambda_1^0,\,\lambda_2^0,\,s,\,s_1,\,s_2,\,x_c,\,x_{c1},\,\widetilde{x}_{c1}$ and $C,\,C_1,\,C_4,\,C_5,\,C_2,\,C_3$.

Solution procedure for system \rf{1m02sol} is very similar to the ones for \rf{0m13sol} and \rf{2m0m13sol} described before. Therefore, we omit the details and state just the solutions formulae. Again all expressions \rf{lambda_grzz}--\rf{C4_expr} remain true with a small exception that signs of \rf{C3} and \rf{C5} are reversed, because {\it $h$-sessile zig-zag} can be viewed as the right part of {\it the 2-side sessile zig-zag} solution (cf. Fig.10 and Fig.12). In turn, the positions expressions take the form:
\bea
s_2&=&\tfrac{C_3}{\lambda_1^0}+L,\quad
x_c=\tfrac{\sqrt{2\sigma|\phi(1)|}-\sigma C_3}{\lambda_2^0-(\sigma+1)\lambda_1^0}-s_2,\quad x_{c1}=-\tfrac{\sqrt{2\sigma|\phi(1)|}}{\lambda_1^0-\lambda_2^0}-s_2,\nonumber\\[1.ex]
s_1&=&\tfrac{\sqrt{2(\sigma+1)\sigma|\phi(1)|}}{\lambda_2^0-(\sigma+1)\lambda_1^0}-x_c=-\tfrac{\sqrt{2\sigma|\phi(1)|/(\sigma+1)}+\sigma C_5}{\lambda_1^0-\lambda_2^0}-x_{c1},\nonumber\\[1.ex]
\widetilde{x}_{c1}&=&\tfrac{(\sigma+1)C_5}{\lambda_2^0}-s_1,\quad s=-\tfrac{\sqrt{2(\sigma+1)|\phi(1)|}}{\lambda_2^0}-\widetilde{x}_{c1}.
\lb{sx_expr_4}
\eea
Typical shapes of {\it $h$-sessile zig-zag} solutions are shown in Fig.12.

The constraints on {\it $h$-sessile zig-zag} solution are deduced from \rf{LC3}. Again, conditions $s_2>s_1$ and $s_1>s_2$ hold for all $\bar{h}=h^m/h_1^m>0$, while $L>s_2$ reduces now to condition $\bar{h}>1$ (cf. Fig.12 {\bf (a)}). In turn, condition $s>0$ impose a bound on the minimal interval length $L$:
\be
L>\tfrac{2|\phi(1)|\sqrt{\sigma}((\sigma+1)\lambda_1^0+\sqrt{\sigma}\lambda_2^0)+C_4(\lambda_2^0-(\sigma+1)\lambda_1^0)^2}{\sqrt{2\sigma(\sigma+1)|\phi(1)|}\lambda_1^0\lambda_2^0}>0\ \text{(cf. Fig.12 {\bf (b)})}.
\lb{hzz_minL}
\ee

Additionally, the combined formula for the maximum of $h_1^0(x)$ takes the form
\be
\displaystyle\max_{x\in(-L,\,0)}h_1^0(x)=\left\{\begin{array}{ll}
	h_1^m& \text{if}\ \ 1<\overline{h}\le\sigma+1\ \text{(cf. Fig.12 {\bf (a)})},\\[1.5ex]
	C_4& \text{if}\ \ \overline{h}\ge\sigma+1\ \text{(cf. Fig.12 {\bf (b)})}\end{array}\right.
\lb{h1m_hgzz}
\ee 
with $C_4$ given by \rf{C2_expr}, while inverting the latter formula for $\overline{h}\ge\sigma+1$ yields:
\be
\displaystyle h_1^m=\left\{\begin{array}{ll}
	\max h_1^0(x)& \text{if}\ \ 1<\widetilde{h}\le\sigma+1,\\[1.5ex]
     	\tfrac{\max h^0(x)}{\sigma+1-2\sigma \widetilde{h}+\sqrt{4\sigma(\sigma+1)\widetilde{h}(\widetilde{h}-1)}}& \text{if}\ \ \widetilde{h}\ge\sigma+1,\end{array}\right.
\lb{h1m_hgzz_inv}
\ee 
with $\widetilde{h}$ as in \rf{ht}.

%%%%%%%%%%%%%%%%%%%%%%%%%%%%%%%%%%%%%%%%%%%%%%%%%%%%%%%%%%%%%%%%%
\section{Combined solution diagrams}
%%%%%%%%%%%%%%%%%%%%%%%%%%%%%%%%%%%%%%%%%%%%%%%%%%%%%%%%%%%%%%%%%
In sections 3-6, we described 11 types of stationary solutions to bilayer system \rf{BS}--\rf{BC} together with the leading order explicit conditions on the {\it model parameters} for which 
such solutions exist and are well defined for all sufficiently small $\eps>0$. In this section, using the latter conditions we plot and discuss the combined diagrams showing for broad ranges
of the {\it model parameters} typical shapes of the existence domains (EDs) for all found solutions together.
In Fig.13, we fix parameters $L$ and $|\phi(1)|$ and plot for three different values of $\sigma$ the solution diagrams w.r.t. two variables
$$h^{max}\!=\!\max_{x\in(-L,\,0)}h(x)\quad\text{and}\quad h_1^{max}\!=\!\max_{x\in(-L,\,0)}h_1(x)$$
with their ratio denoted below by $\tilde{h}$ similarly as in \rf{ht}. The domain boundaries for different solutions are indicated in Fig.13 by numbers as follows:
$\raisebox{.5pt}{\textcircled{\raisebox{-.9pt} {1}}}$--lens, $\raisebox{.5pt}{\textcircled{\raisebox{-.9pt} {2}}}$--internal drop, $\raisebox{.5pt}{\textcircled{\raisebox{-1.1pt} {3}}}$--$h_1$-drop,
$\raisebox{.5pt}{\textcircled{\raisebox{-.9pt} {4}}}$--$h$-drop, $\raisebox{.5pt}{\textcircled{\raisebox{-1.pt} {5}}}$--zig-zag, $\raisebox{.5pt}{\textcircled{\raisebox{-.9pt} {6}}}$--sessile lens,
$\raisebox{.5pt}{\textcircled{\raisebox{-1.pt} {7}}}$--sessile internal drop, $\raisebox{.5pt}{\textcircled{\raisebox{-1.pt} {8}}}$-$2$-drops, $\raisebox{.5pt}{\textcircled{\raisebox{-.9pt} {9}}}$--2-side sessile zig-zag, $\raisebox{.5pt}{\textcircled{\hspace{-.05cm}\raisebox{-1.1pt} {1\!0}}}$--$h_1$-sessile zig-zag, $\raisebox{.5pt}{\textcircled{\raisebox{-1.pt} {1\!1}}}$--$h$-sessile zig-zag.

From their explicit formulae we observe that all EDs are open and simply connected and their boundaries in many cases, except only for solutions $\raisebox{.5pt}{\textcircled{\raisebox{-.9pt} {6}}}$--$\raisebox{.5pt}{\textcircled{\raisebox{-1.pt} {7}}}$
and $\raisebox{.5pt}{\textcircled{\raisebox{-.9pt} {9}}}$--$\raisebox{.5pt}{\textcircled{\raisebox{-1.pt} {1\!1}}}$, are formed by intersecting straight lines. For example, EDs of lens $\raisebox{.5pt}{\textcircled{\raisebox{-1.1pt} {1}}}$ or internal drop $\raisebox{.5pt}{\textcircled{\raisebox{-.9pt} {2}}}$ solutions are semi-infinite stripes bounded by $\tilde{h}=\sigma+1$ or $\tilde{h}=1$ lines, respectively, and another vertical or horizontal ones (cf. \rf{lens_constr} and \rf{intdr_constr}); $2$-drops $\raisebox{.5pt}{\textcircled{\raisebox{-1.pt} {8}}}$ ED is given by a right triangle with its hypotenuse connecting the base points of the horizontal and vertical lines forming the ED boundaries for $h_1$-drop $\raisebox{.5pt}{\textcircled{\raisebox{-1.1pt} {3}}}$ and $h$-drop $\raisebox{.5pt}{\textcircled{\raisebox{-.9pt} {4}}}$ solutions, respectively (cf. \rf{2dr_constr} and \rf{h1dr_constr}--\rf{hdr_constr}); zig-zag solution $\raisebox{.5pt}{\textcircled{\raisebox{-1.pt} {5}}}$ lives in a pentagon domain formed by the intersection of the coordinate axes and three lines (cf. \rf{zz_constr}):
\sbea
\lb{BB_zz}
\hspace{-1.cm}(\sqrt{\sigma+1}+1)(h^{max}+\sqrt{\sigma+1}h_1^{max})&=&\tilde{L},\quad(\sigma+1)h_1^{max}-h^{max}=\tilde{L},\\
\sqrt{\sigma+1}(h^{max}-h_1^{max})&=&\tilde{L}\quad\text{with}\quad\tilde{L}=L\sqrt{\tfrac{\sigma|\phi(1)|}{2}}.
\seea
For {\it sessile lens} $\raisebox{.5pt}{\textcircled{\raisebox{-.9pt} {6}}}$ the curvelinear shape of its ED boundary qualitatively differs in two cases when $\sigma<1$ or $\sigma>1$ (cf. \rf{h1m_ls}--\rf{gl_constr}).
For $\sigma<1$ the ED is confined to a sector bounded by parabola
\bes
h_1^{max}=\tfrac{\tilde{L}^2-\sigma|h^{max}|^2}{(\sigma+1)\tilde{L}}\quad\text{with}\ \tilde{L}=L\sqrt{\tfrac{(\sigma+1)|\phi(1)|}{2}}\quad\text{and line}\  \tilde{h}=\tfrac{1+\sigma}{1-\sigma}.
\ees
The latter boundary line corresponds to singular merge $h_1^0(0)\go 0$ (cf. Fig.7 {\bf (a)}) and coincides with $ h^{max}$ axis when $\sigma\!=\!1$.
For $\sigma>1$ the ED is a sector bounded by the same parabola and two other lines:
\bes
h_1^{max}=\tfrac{\tilde{L}-h^{max}}{\sigma+1}\quad\text{and}\  \tilde{h}=\tfrac{1+\sigma}{\sqrt{\sigma}-1}
\ees
with the latter one corresponding to the singular merge shown in  Fig.7 {\bf (d)}.
%%%%%%%%%%%%%%%%%%%%%%%%%%%%%%%%%%%%%%%%%%%%%%%%%%%%%%%%%%%%%%%%%%%%%%%%%%%%%%%%%%%%%%%%%%%%%%%%%%%%%%%%%%
\begin{figure}[H] 
	\centering
	\vspace{-4.2cm}
	\hspace{-.cm}\includegraphics[width=.8\textwidth]{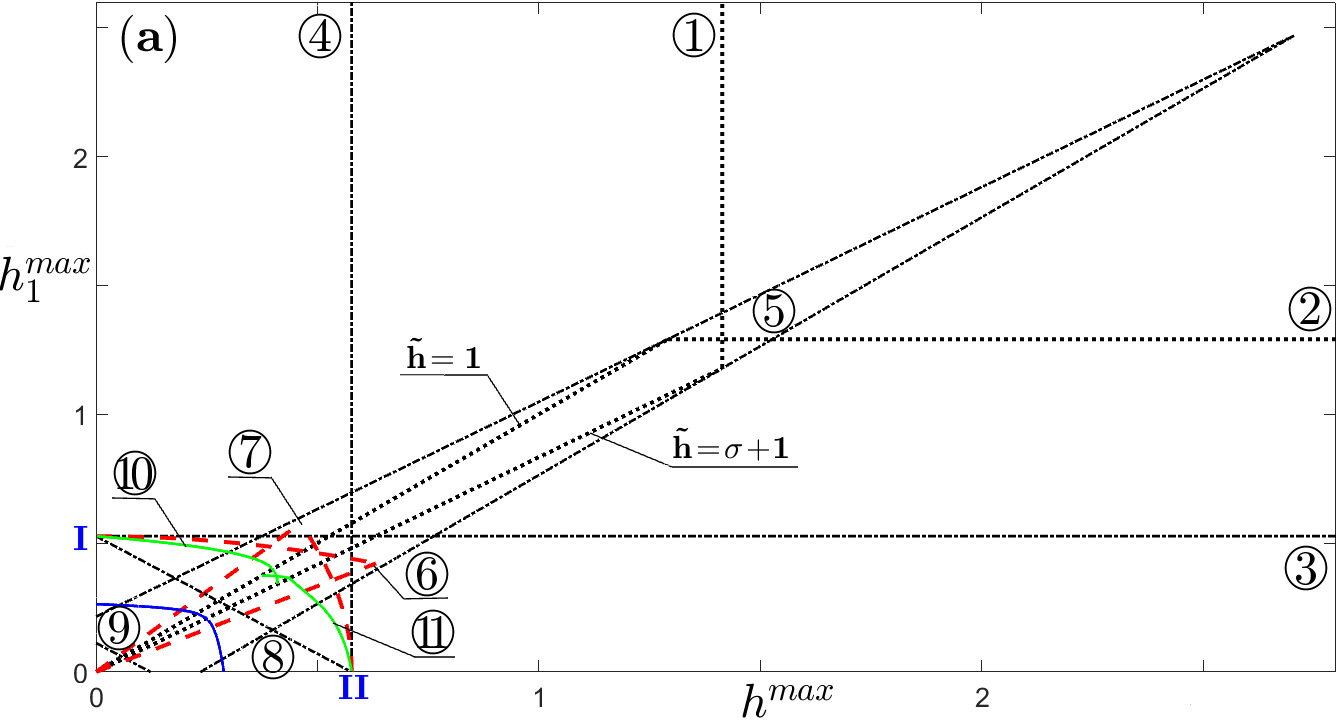}\\[1ex]   	
	\hspace{-.cm}\includegraphics[width=.78\textwidth]{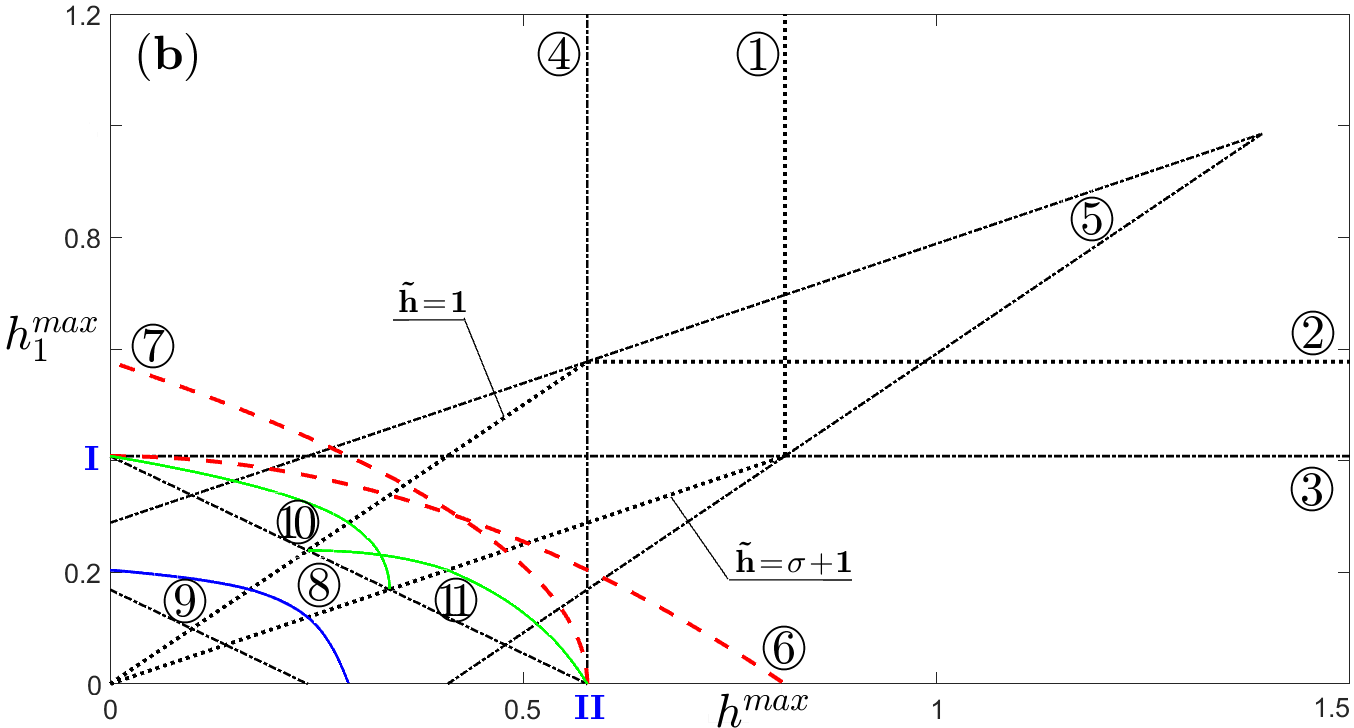}\\[1ex]   
	\hspace{-.cm}\includegraphics[width=.78\textwidth]{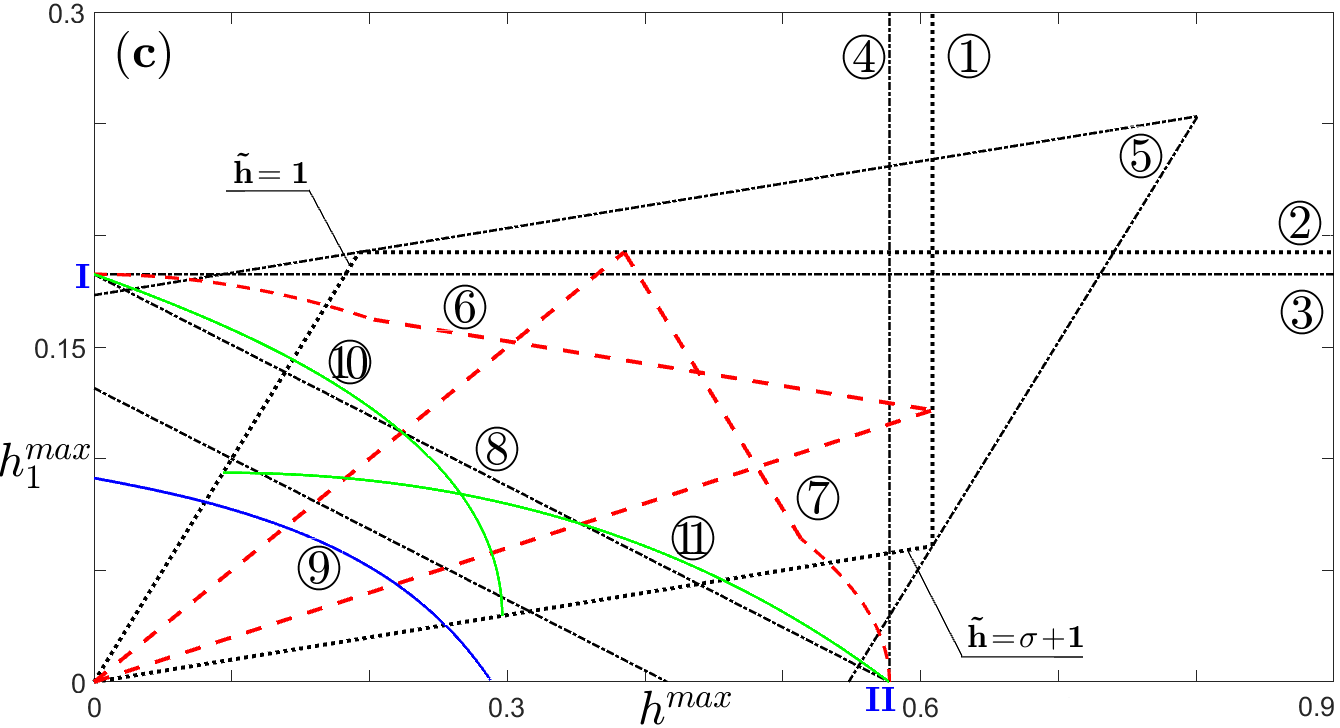}\\[1ex]  
	\hspace{-.cm}\includegraphics[width=.55\textwidth]{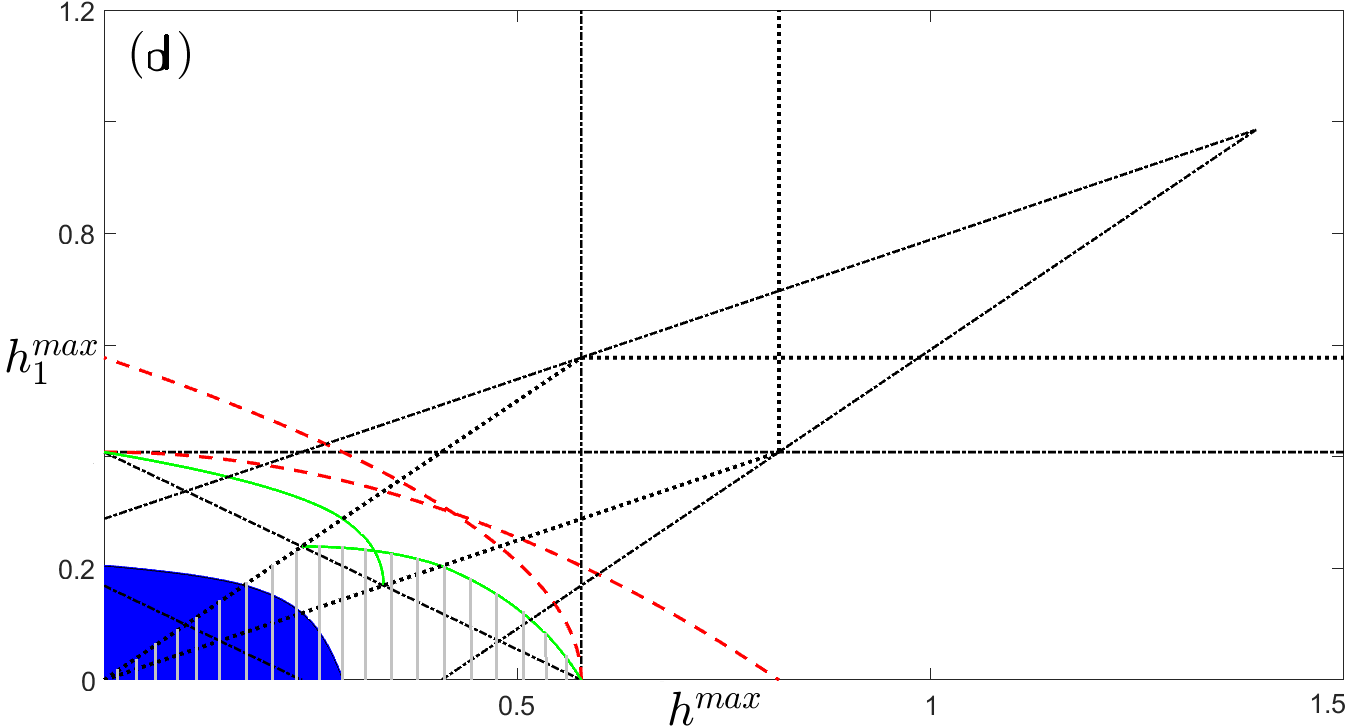}	
	\caption{\small Combined diagrams showing existence domains (EDs) for 11 solution types found in this article for $L=2$ and {\bf (a)} $\sigma=0.2$; {\bf (b)} and  {\bf (d)} $\sigma=1$; {\bf (c)} $\sigma=9$.
		ED boundaries are shown for  the $1^{\mathrm{st}}-2^{\mathrm{nd}}$ solutions as dotted, $3^{\mathrm{rd}}-5^{\mathrm{th}}$ and $8^{\mathrm{th}}$ as dash--dotted, $6^{\mathrm{th}}-7^{\mathrm{th}}$ as dashed,  $9^{\mathrm{th}}-11^{\mathrm{th}}$ as solid curves.  In {\bf (d)}, EDs of $9^{\mathrm{th}}$ and $11^{\mathrm{th}}$ solutions are filled with solid and dashed colors, respectively. {\it Symmetric points} $\mathrm{\tb\bf I}$ and $\mathrm{\tb\bf II}$ are located at $h_1^{max}=L\sqrt{\frac{|\phi(1)|}{2(\sigma+1)}}$ and $h^{max}=L\sqrt{\frac{|\phi(1)|}{2}}$, respectively. }
\end{figure}
%%%%%%%%%%%%%%%%%%%%%%%%%%%%%%%%%%%%%%%%%%%%%%%%%%%%%%%%%%%%%%%%%%%%%%%%%%%%%%%%%%%%%%%%%%%%%%%%%%%%%%%%%%%%
Similarly, for {\it sessile internal drop} $\raisebox{.5pt}{\textcircled{\raisebox{-1.pt} {7}}}$ the shape of its ED depends on the case $\sigma<1$ or $\sigma>1$ (cf. \rf{hm_ids}--\rf{gdr_constr}).
For $\sigma<1$ the ED is a sector bounded by parabola
\bes
h^{max}=\tfrac{\tilde{L}^2-\sigma|h_1^{max}|^2}{\tilde{L}}\quad\text{with}\ \tilde{L}=L\sqrt{\tfrac{|\phi(1)|}{2}}\quad\text{and line}\  \tilde{h}=1-\sigma
\ees
coinciding with $h_1^{max}$ axis when $\sigma=1$. For $\sigma>1$ the ED sector is bounded by the same parabola and two other lines:
\bes
h^{max}=\tilde{L}-h_1^{max}\quad\text{and}\  \tilde{h}=\sqrt{\sigma}-1
\ees
with the latter one corresponding to the singular merge shown in  Fig.8 {\bf (b)}.

For {\it 2-side sessile zig-zag} $\raisebox{.5pt}{\textcircled{\raisebox{-.9pt} {9}}}$ ED is a sector (cf.  Fig.13 {\bf (d)}) formed by the intersection of the coordinated axes $h^{max}$ and  $h_1^{max}$
and the curve
\bes
L\sqrt{\tfrac{|\phi(1)|}{2}}=h_1^m\left[\bar{h}+\sqrt{\sigma+1}+\tfrac{-\sqrt{\sigma+1}\sqrt{\sigma\bar{h}^2+(\bar{h}-\sigma-1)^2}+\bar{h}\sqrt{\sigma+(\bar{h}-1)^2}}{\bar{h}-\sigma-1}\right].
\ees
This formula follows from the constraints \rf{2sd_cons} after substituting into them expressions \rf{sx_expr_2} where explicit formulae  \rf{lambda_grzz}, \rf{C2_expr}, \rf{C4_expr} and \rf{C3}, \rf{C5} are used  with as before $\bar{h}=h^m/h_1^{m}$. By plotting this curve in Fig.13 we employed the transformation relations \rf{hm_h1m_2sgzz}, \rf{hm_h1m_2sgzz_inv}  between $(h_1^m,\,h^m)$ and $(h_1^{max},\,h^{max})$.

For {\it$h_1$-sessile zig-zag} $\raisebox{.5pt}{\textcircled{\hspace{-.05cm}\raisebox{-1.1pt} {1\!0}}}$ ED is a sector formed by the intersection of line $\tilde{h}=\sigma+1$
and the curve
\bes
h^m=\tfrac{\tilde{L}^2-(\sigma+1)|h_1^{max}|^2}{2(\tilde{L}-h_1^{max})}\quad\text{with}\ \tilde{L}=L\sqrt{\tfrac{|\phi(1)|}{2}}.
\ees
This formula follows from the constraint \rf{h1zz_minL} using \rf{lambda_grzz} and \rf{C4_expr}. By plotting this curve in Fig.13 we employed the transformation relations \rf{hm_h1gzz}--\rf{hm_h1gzz_inv}  between $h^m$ and $h^{max}$. Similarly, for {\it$h$-sessile zig-zag} $\raisebox{.5pt}{\textcircled{\raisebox{-1.pt} {1\!1}}}$ ED is a sector (cf.  Fig.13 {\bf (d)}) formed by the intersection of line $\tilde{h}=1$ and the curve
\bes
h_1^m=\tfrac{\tilde{L}^2-(\sigma+1)|h^{max}|^2}{2(\sigma+1)(\tilde{L}-h^{max})}\quad\text{with}\ \tilde{L}=L\sqrt{\tfrac{(\sigma+1)|\phi(1)|}{2}}.
\ees
This formula follows from the constraint \rf{hzz_minL} using \rf{lambda_grzz} and \rf{C2_expr}. By plotting this curve in Fig.13 we employed the transformation relations \rf{h1m_hgzz}--\rf{h1m_hgzz_inv}  between $h_1^m$ and $h_1^{max}$.

In summary, the EDs boundaries of the solutions to system \rf{BS}--\rf{BC} and its stationary counterpart \rf{SSa}--\rf{SSc} derived in sections 3-6 depend  linearly on the product $L\sqrt{|\phi(1)|/2}$ and, therefore, their form is simply rescaled under changes of parameters $L$ or $|\phi(1)|$, while qualitatively differs only in three cases $\sigma>1$, $\sigma=1$ or $\sigma<1$ demonstrated in Fig.13. Interestingly, all presented solution diagrams are symmetric around the special line $\tilde{h}=\sqrt{\sigma+1}$ connecting the origin with the acute peak of the ED boundary of {\it zig-zag} $\raisebox{.5pt}{\textcircled{\raisebox{-1.pt} {5}}}$. More precisely, the solution diagrams do not change if they are reflected around  $\tilde{h}=\sqrt{\sigma+1}$ line  in the transversal direction of $\tilde{h}=-\sqrt{\sigma+1}$ one. Note that the latter is the direction of the ED boundary of $2$-drops solution $\raisebox{.5pt}{\textcircled{\raisebox{-1.pt} {8}}}$  as well as of the first line in \rf{BB_zz}. In special case $\sigma=1$, the solution diagram  possesses additional symmetries (see Fig.13 {\bf (b), (d)}). For example, the ED boundaries for solutions $\raisebox{.5pt}{\textcircled{\raisebox{-.9pt} {2}}}$, $\raisebox{.5pt}{\textcircled{\raisebox{-1.pt} {8}}}$ and $\raisebox{.5pt}{\textcircled{\raisebox{-1.pt} {1\!1}}}$ intersect at a single joint point, while those for $\raisebox{.5pt}{\textcircled{\raisebox{-.9pt} {2}}}$, $\raisebox{.5pt}{\textcircled{\raisebox{-.9pt} {4}}}$ and $\raisebox{.5pt}{\textcircled{\raisebox{-1.pt} {5}}}$ at another one.

Additionally, Fig.13 reveals special {\it symmetric points} depicted there as $\mathrm{\tb\bf I}$ and $\mathrm{\tb\bf II}$ at which two quadruples of ED boundaries originate together with three out of four in them tangentially to each other. At point  $\mathrm{\tb\bf I}$, these are ED boundaries for $\raisebox{.5pt}{\textcircled{\raisebox{-1.1pt} {3}}}$,$\raisebox{.5pt}{\textcircled{\raisebox{-.9pt} {6}}}$, $\raisebox{.5pt}{\textcircled{\raisebox{-1.pt} {8}}}$ and $\raisebox{.5pt}{\textcircled{\hspace{-.05cm}\raisebox{-1.1pt} {1\!0}}}$, while at  $\mathrm{\tb\bf II}$ for $\raisebox{.5pt}{\textcircled{\raisebox{-.9pt} {4}}}$, $\raisebox{.5pt}{\textcircled{\raisebox{-1.pt} {7}}}$, $\raisebox{.5pt}{\textcircled{\raisebox{-1.pt} {8}}}$ and $\raisebox{.5pt}{\textcircled{\raisebox{-1.pt} {1\!1}}}$. Using the explicit solution formulae stated in sections $3-6$ we find that the leading order profiles of $\raisebox{.5pt}{\textcircled{\raisebox{-.9pt} {6}}}$, $\raisebox{.5pt}{\textcircled{\raisebox{-1.pt} {8}}}$ and $\raisebox{.5pt}{\textcircled{\hspace{-.05cm}\raisebox{-1.1pt} {1\!0}}}$ ones converge pointwise in $x$ variable to that one of $\raisebox{.5pt}{\textcircled{\raisebox{-1.1pt} {3}}}$, while the ones of  $\raisebox{.5pt}{\textcircled{\raisebox{-1.pt} {7}}}$, $\raisebox{.5pt}{\textcircled{\raisebox{-1.pt} {8}}}$ and $\raisebox{.5pt}{\textcircled{\raisebox{-1.pt} {1\!1}}}$ to the one of $\raisebox{.5pt}{\textcircled{\raisebox{-.9pt} {4}}}$ when points $\mathrm{\tb\bf I}$ or $\mathrm{\tb\bf II}$, respectively, are approached in Fig.13. Based on that we conjecture that $\mathrm{\tb\bf I}$ and $\mathrm{\tb\bf II}$ are bifurcation points for some of these solutions.

Finally, we note the following limiting behavior by approaching the parts of ED boundaries lying on coordinate axes $h^{max}$ or $h_1^{max}$ in Fig.13: solutions $\raisebox{.5pt}{\textcircled{\raisebox{-.9pt} {1}}}$ and  $\raisebox{.5pt}{\textcircled{\raisebox{-.9pt} {2}}}$ converge pointwise in $x$ to constant ones; all other solutions except $\raisebox{.5pt}{\textcircled{\raisebox{-1.pt} {5}}}$ and $\raisebox{.5pt}{\textcircled{\raisebox{-.9pt} {9}}}$ to $h_1$-drop $\raisebox{.5pt}{\textcircled{\raisebox{-1.1pt} {3}}}$ or $h$-drop $\raisebox{.5pt}{\textcircled{\raisebox{-.9pt} {4}}}$, respectively;  $\raisebox{.5pt}{\textcircled{\raisebox{-.9pt} {9}}}$  converges to symmetric $h_1$ or $h$-drops, respectively, while
 $\raisebox{.5pt}{\textcircled{\raisebox{-1.pt} {5}}}$ to some different single drop like profiles having not-stationary contact angles.
 
 %%%%%%%%%%%%%%%%%%%%%%%%%%%%%%%%%%%%%%%%%%%%%%%%%%%%%%%%%%%%%%%%%
 \section{Unstable coarsening solutions}
 %%%%%%%%%%%%%%%%%%%%%%%%%%%%%%%%%%%%%%%%%%%%%%%%%%%%%%%%%%%%%%%%%
 
In previous sections, we described ten types of stable stationary solutions to bilayer system \rf{BS}--\rf{IC} and a special one ({\it2-side sessile zig-zag}) which moves due to exponentially small translation instabilities (cf. Fig.10 {\bf (d)}), but retains its leading order shape and can be made stable by changing the boundary conditions \rf{BC}, e.g. to periodic ones.
 
By deriving the leading order profiles for these eleven solutions in sections 3--6 the main observation used was that they can be composed by asymptotic matching  of chains of {\it bulk, contact line and UTF} profiles introduced in section 2. As a consequence of that or alternatively, any solution having several CLs can be viewed as a composition of several one-CL solutions (serving as the building blocks), namely  {\it lens} ($\bf 0^-$), {\it internal drop} ($\bf1^-$), {\it$h_1$-drop} ($\bf2$), and {\it$h$-drop} ($\bf3$) ones, considered in section 3. The number abbreviations stated in brackets for the latter solutions introduce a useful notation for generating chain shortcuts for the composite solutions. For example, we abbreviate two-CLs {\it zig-zag} solution by $\bf(1^-0)$ and four-CLs {\it 2-side sessile zig-zag} one by $\bf(2^-0^-13)$. Note that possible $x$ variable reversion in the one-CL solutions is also taken into account in this notation: e.g. $\bf(0)$ abbreviates the {\it lens} with its center located at $x=-L$ and having the inverted in $x$ profile to that one of $\bf(0^-)$ shown in Fig.2 {\bf (a)}. Obviously, this composition approach can be further used to build and test even more complex chains with arbitrary numbers of one-CL solutions involved.   
%%%%%%%%%%%%%%%%%%%%%%%%%%%%%%%%%%%%%%%%%%%%%%%%%%%%%%%%%%%%%%%%%%%%%%%%%%%%%%%%%%%%%%%%%%%%%%%%%%%%%%%%%%%%
\begin{figure}[H] 
	\centering
	\vspace{-2.9cm}
	\hspace{-2.4cm}\includegraphics[width=.38\textwidth]{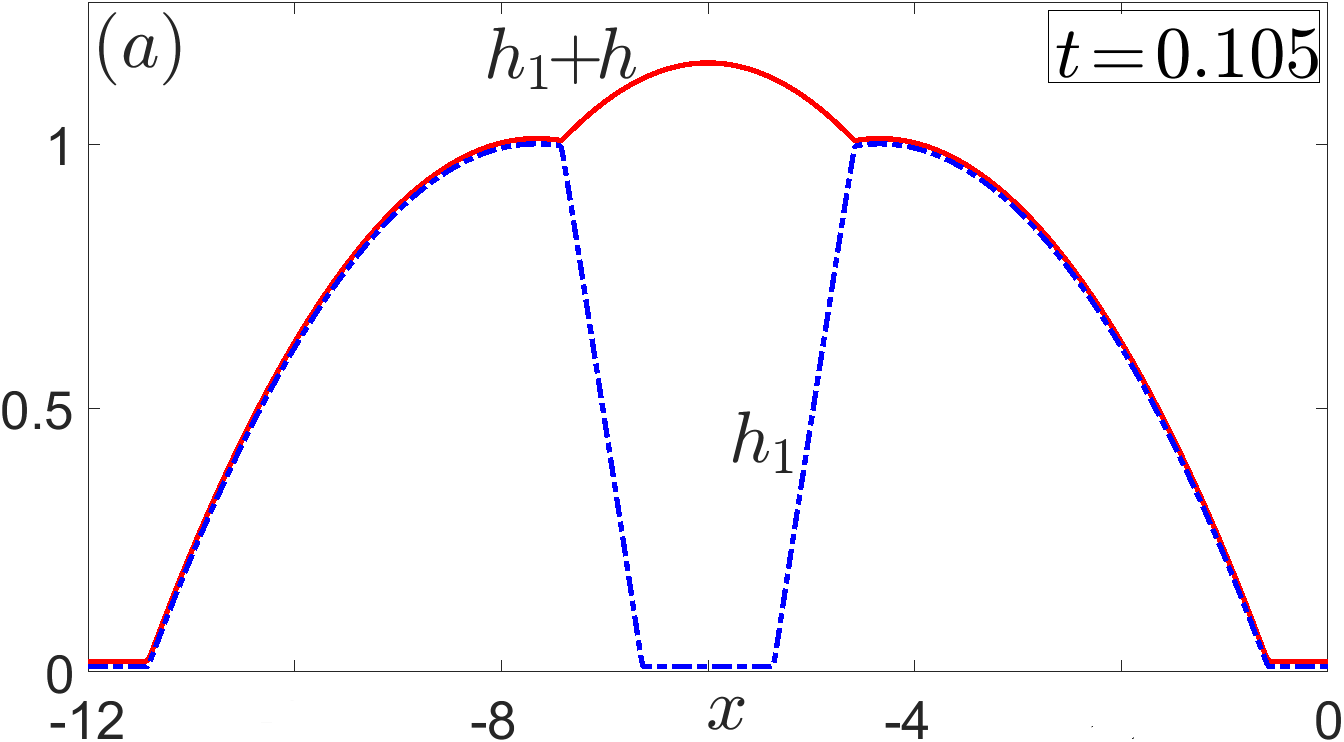}  	
	\hspace{.1cm}\includegraphics[width=.38\textwidth]{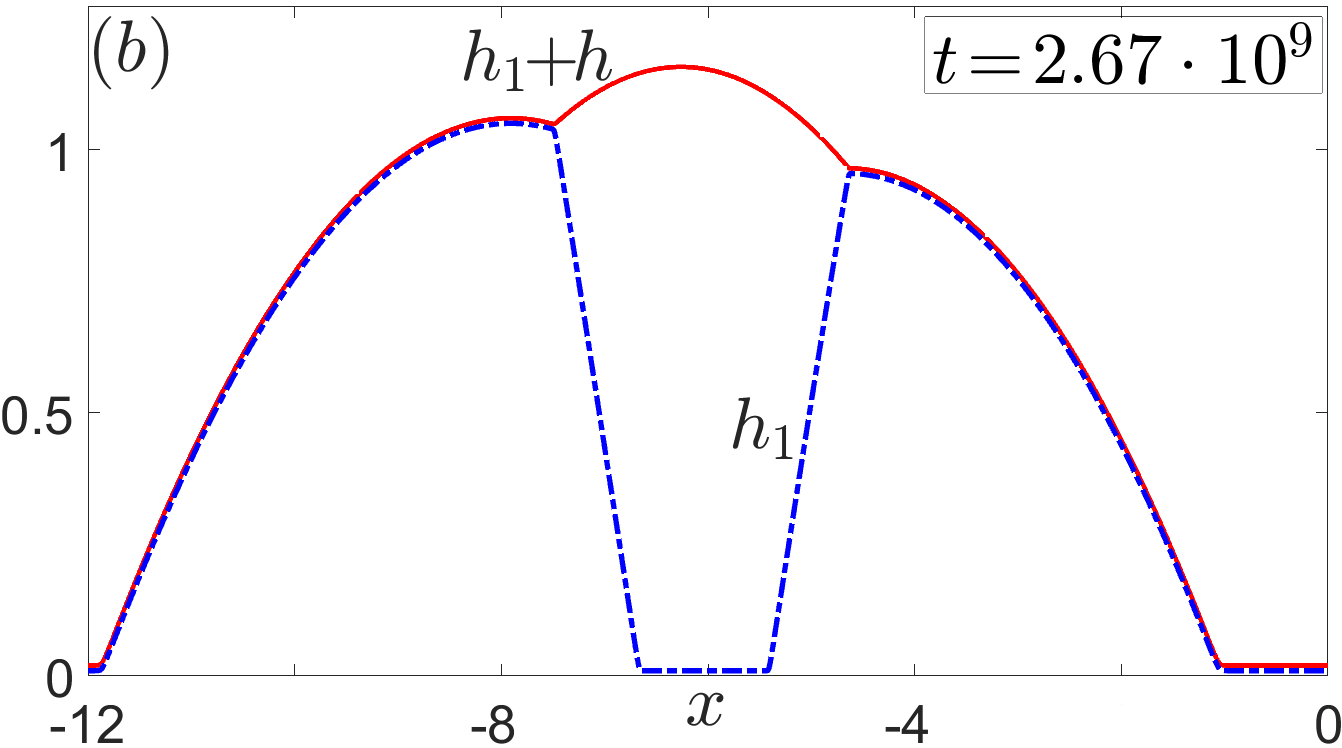}
	\hspace{.1cm}\includegraphics[width=.38\textwidth]{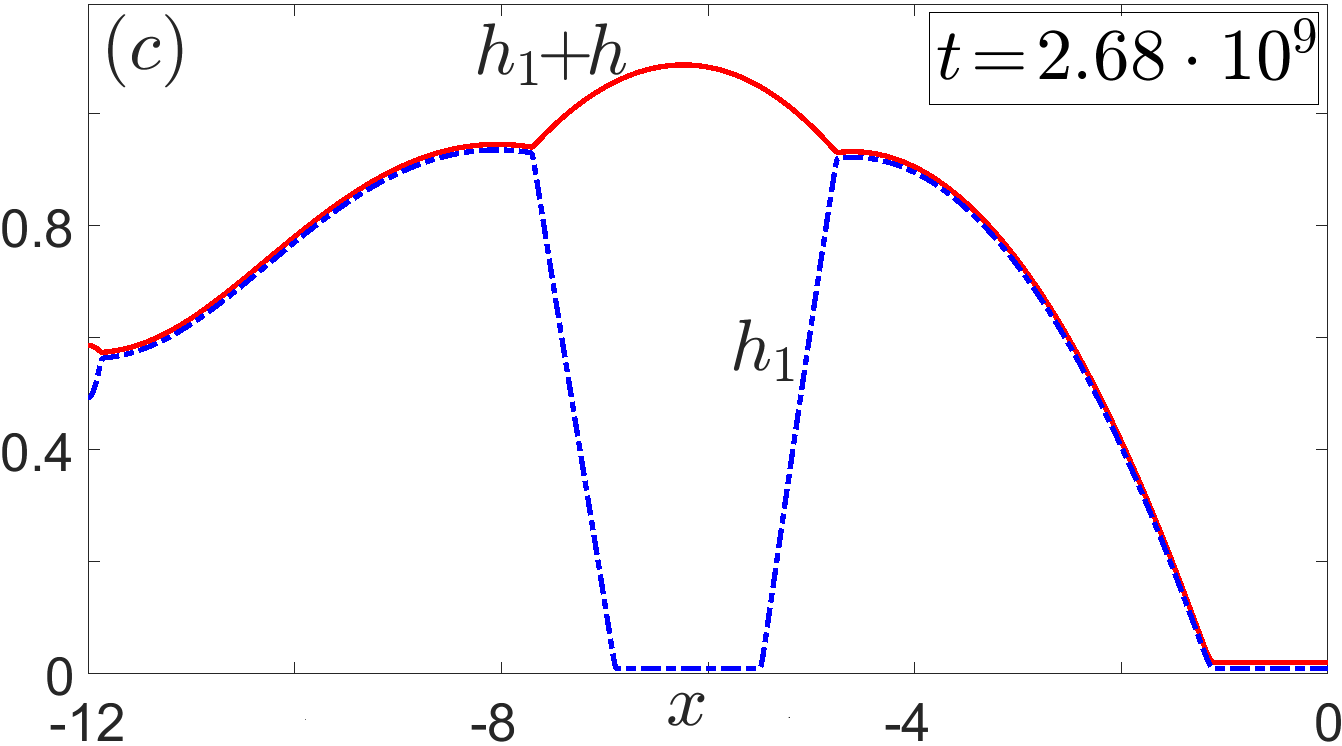}\\[1ex]	
	\hspace{-2.4cm}\includegraphics[width=.38\textwidth]{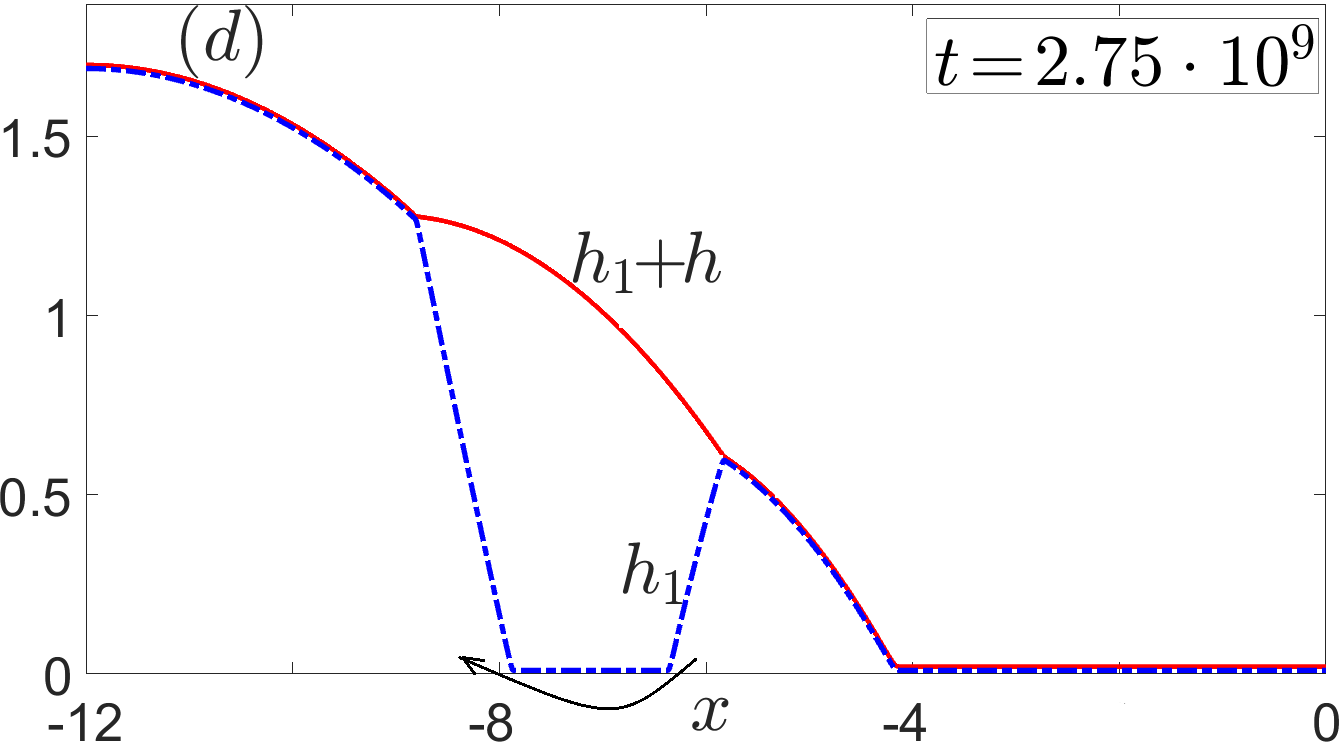}  	
	\hspace{.1cm}\includegraphics[width=.38\textwidth]{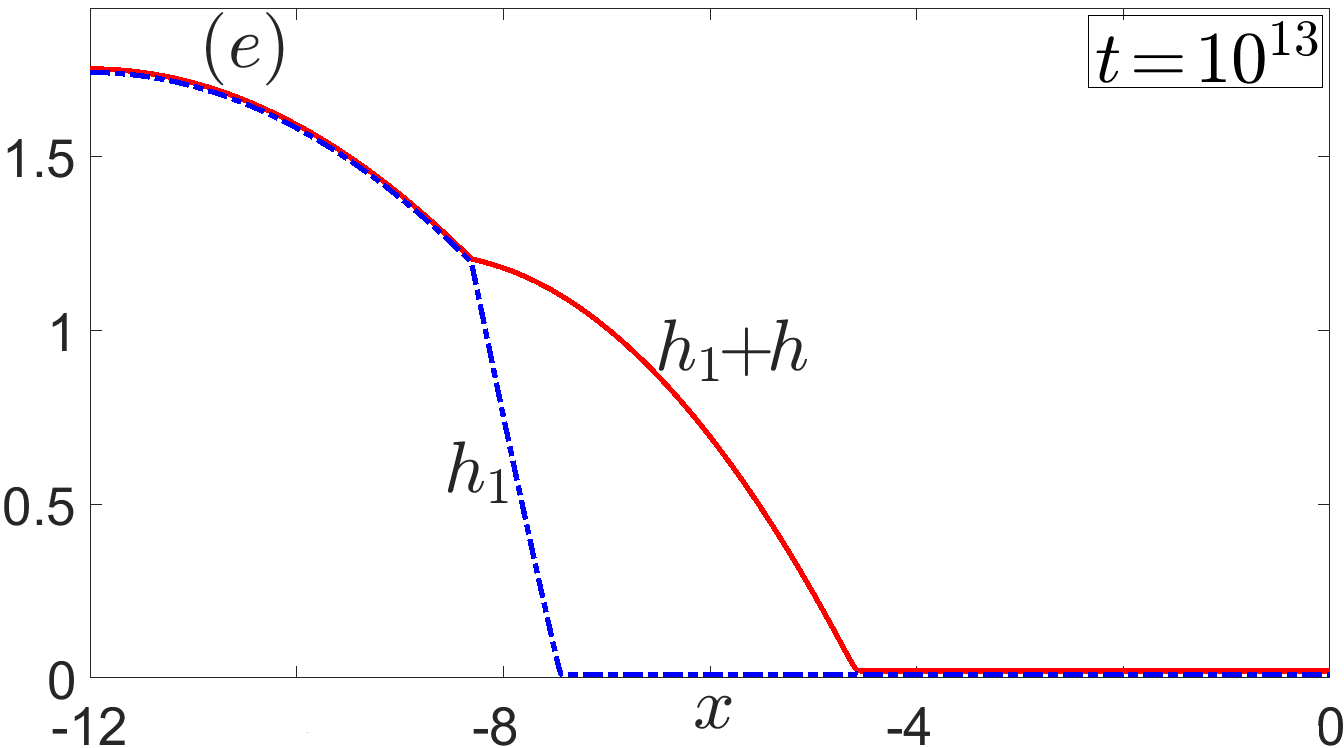}
	\hspace{.1cm}\includegraphics[width=.38\textwidth]{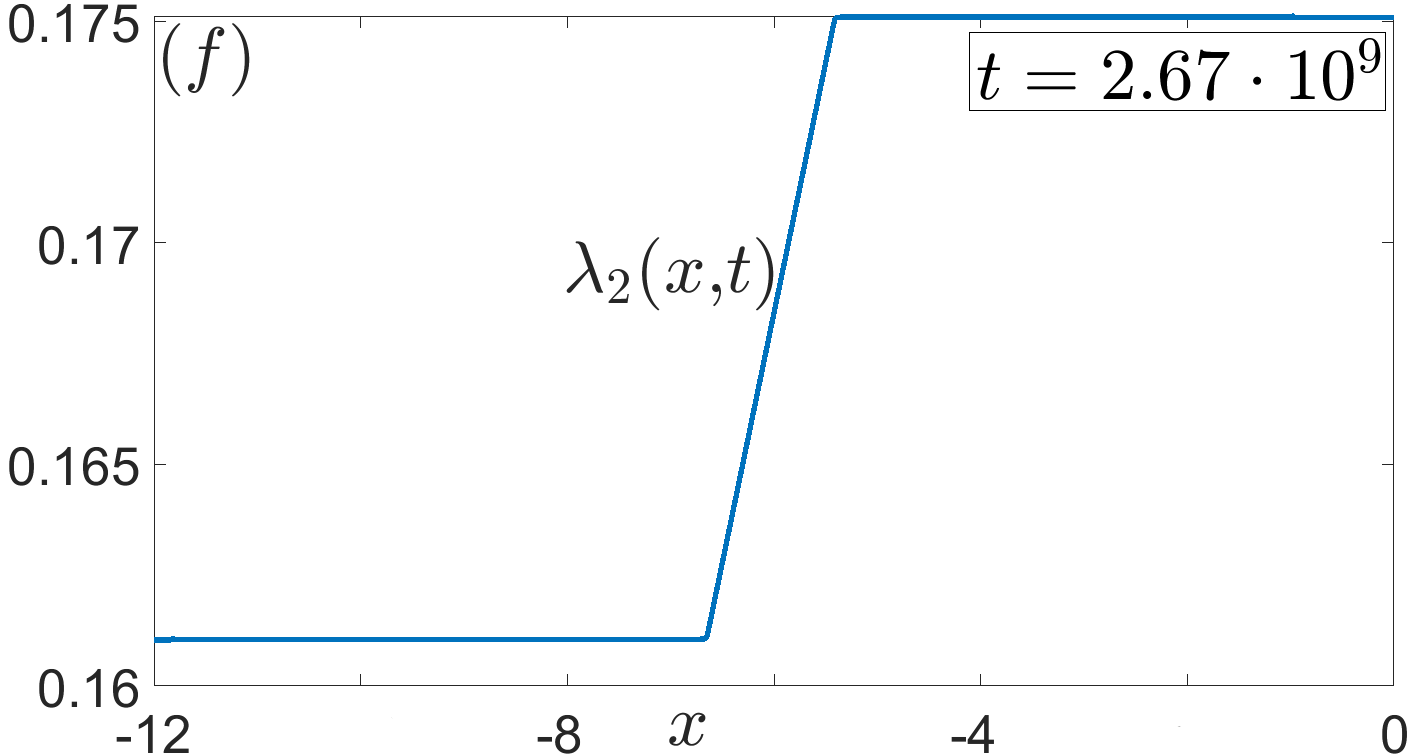}			
	
	\caption{\small {\bf(a)}--{\bf(e)}: different time snapshots for coarsening numerical solution to \rf{BS}--\rf{BC} considered with $L\!\!=\!\!12,\,\sigma\!\!=\!\!0.2,\,\eps\!\!=\!\!0.01$. Initial condition \rf{IC} is symmetric {\it $h$-sessile zig-zag} solution $(\bf2^-0^-11^-02)$ to stationary system \rf{SSa}--\rf{SSc}. {\bf(f)}: numerical pressure profile $\lambda_2(x,t)\!\!=\!\!\Pi_\eps(h_1(x,t))\!\!-\!\!h_{xx}(x,t)\!-\!(\sigma+1)h_{1,xx}(x,t)$ for {\bf(b)}.}
\end{figure}
%%%%%%%%%%%%%%%%%%%%%%%%%%%%%%%%%%%%%%%%%%%%%%%%%%%%%%%%%%%%%%%%%%%%%%%%%%%%%%%%%%%%%%%%%%%%%%%%%%%%%%%%%%%%
The main result of this section is to provide a numerical evidence that other composite solutions (beside the ones described in sections 3--4 and 6) are dynamically unstable when considered as initial profiles \rf{IC} for system \rf{BS}--\rf{BC}. For that, firstly, we demonstrate using  simulations presented in Fig.14-16 that all symmetric versions of the solutions introduced in sections 3-6, namely those that are obtained by reflection around boundaries $x=0$ or $x=-L$, turn out to be numerically unstable when set as initial profiles \rf{IC} for \rf{BS}-\rf{BC}.

In Fig.14, we take as initial condition \rf{IC} a symmetric version of {\it $h$-sessile zig-zag} solution abbreviating as ($\bf 2^-0^-11^-02$) and having six CLs, with its center located at $x=-L/2$. The corresponding numerical solution to system \rf{BS}--\rf{BC} exhibits {\it coarsening} characterized by a drift of the composite drop to $x=-L$ and profound symmetry break of its shape (cf. Fig.14 {\bf (b)}). Due to that pressure profile $\lambda_2(x,\cdot)$ (cf. Fig.14 {\bf (f)}) becomes non-constant and shows small magnitude but profound spatial variation in the middle region (where $h_1(x)\approx\eps$) of the composite drop. Such pressure variations are also typical for the {\it coarsening dynamics} observed in single layer thin liquid films~\cite{GW03,GW05}. Moreover, after the boundary touch a new small {\it lens subdrop} nucleates at $x=-L$ (cf. Fig.14 {\bf (c)}) and collapses after relatively short time while, subsequently, a new composite five-CL solution $\bf (0^-11^-02)$ forms. Finally, $\bf(0^-11^-02)$ coarsens in $h_1$ layer via the induced mass flux (depicted by an arrow in Fig.14 {\bf (d)}) from  the smaller bulk region towards the larger one and converges in the long time to stationary {\it $h_1$-sessile zig-zag} solution $\bf(0^-13)$ (cf. Fig.14 {\bf (e)}). 

In Fig.15, initial profile \rf{IC} is a symmetric version of {\it sessile lens} abbreviating as  $(\bf2^-0^-02)$ and having four CLs, with its center located at $x=-L/2$. Again the corresponding numerical solution to \rf{BS}--\rf{BC} coarsens, drifts to $x=-L$ and exhibits profound symmetry break of its {\it'lens with a collar'} like shape (cf. Fig.15 {\bf (a)-(c)}). By that pressure profile $\lambda_1(x,\cdot)$ is non-constant and exhibits two jumps at the collar lateral sides (cf. Fig.15 {\bf (i)}). Next, the collar breaks into two internal drops (cf. Fig.15 {\bf (d)}) and after subsequent boundary touch the whole compound drop quickly reforms its $h_1$ layer (cf. Fig.15 {\bf (e)}--{\bf (g)}) yielding a  composite four-CL solution $\bf (11^-03)$. 
%%%%%%%%%%%%%%%%%%%%%%%%%%%%%%%%%%%%%%%%%%%%%%%%%%%%%%%%%%%%%%%%%%%%%%%%%%%%%%%%%%%%%%%%%%%%%%%%%%%%%%%%%%%%
\begin{figure}[H] 
	\centering
	\vspace{-2.9cm}
	\hspace{-2.4cm}\includegraphics[width=.38\textwidth]{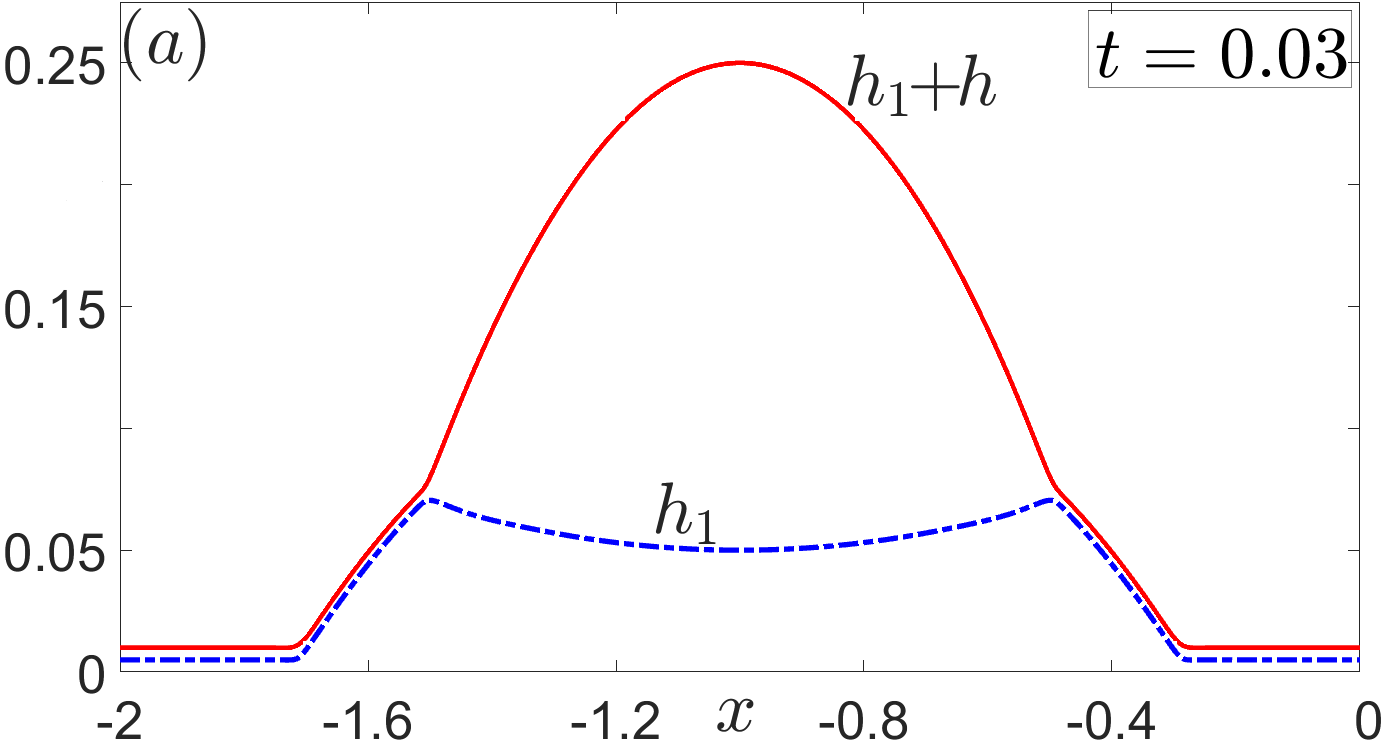}  	
	\hspace{.1cm}\includegraphics[width=.38\textwidth]{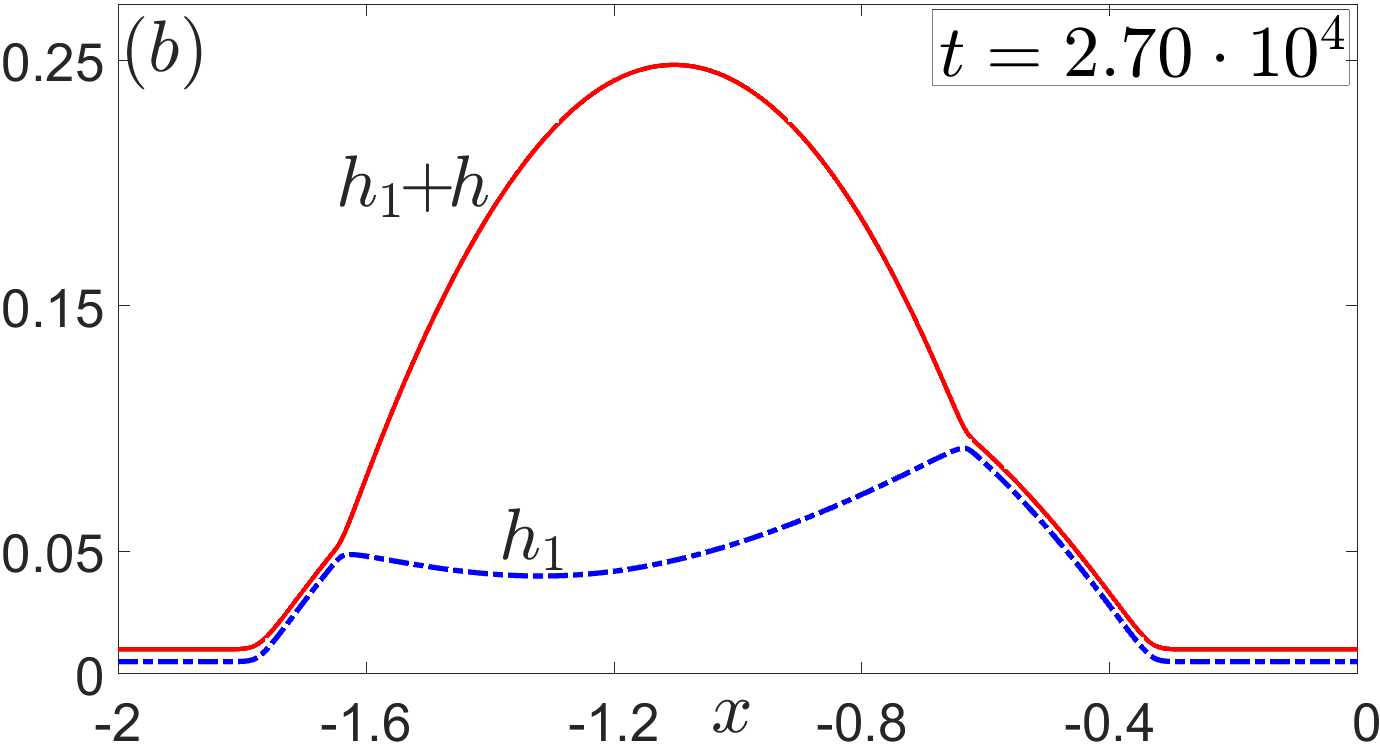}
	\hspace{.1cm}\includegraphics[width=.38\textwidth]{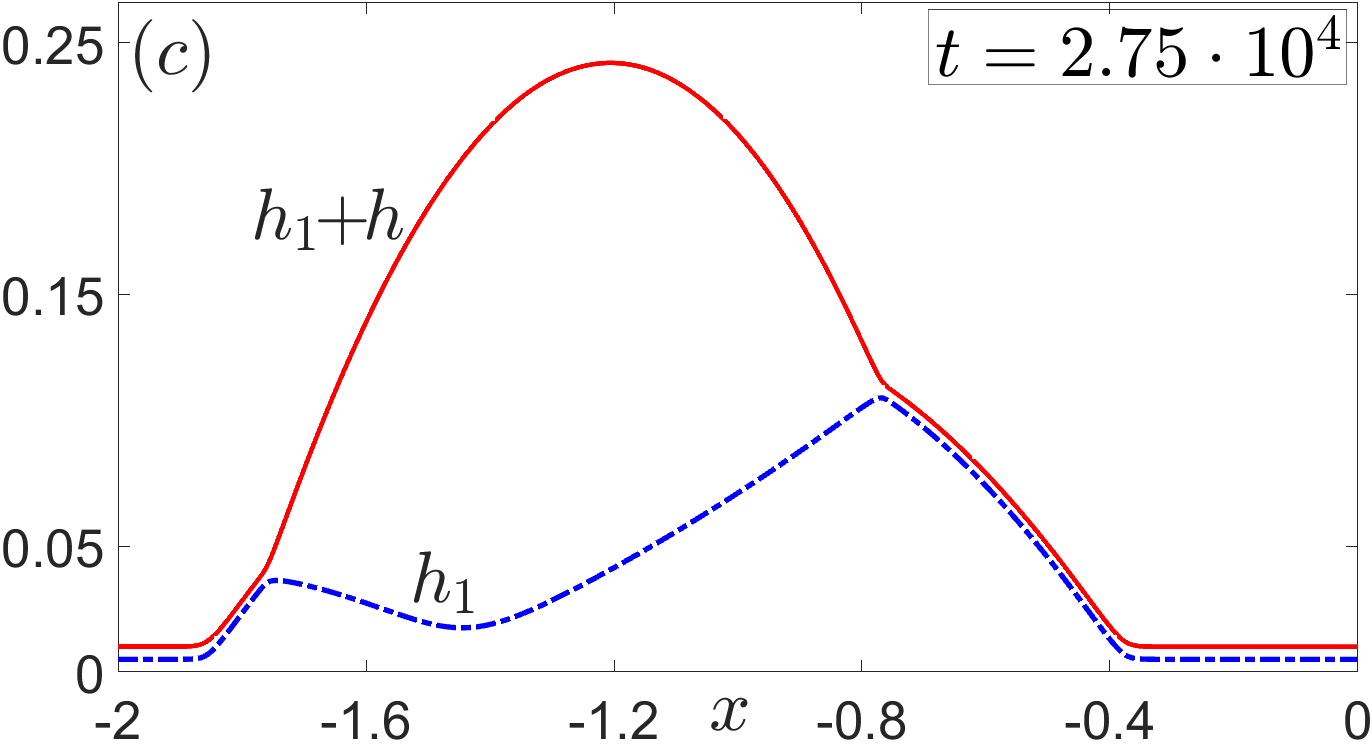}\\[1ex]	
	\hspace{-2.4cm}\includegraphics[width=.38\textwidth]{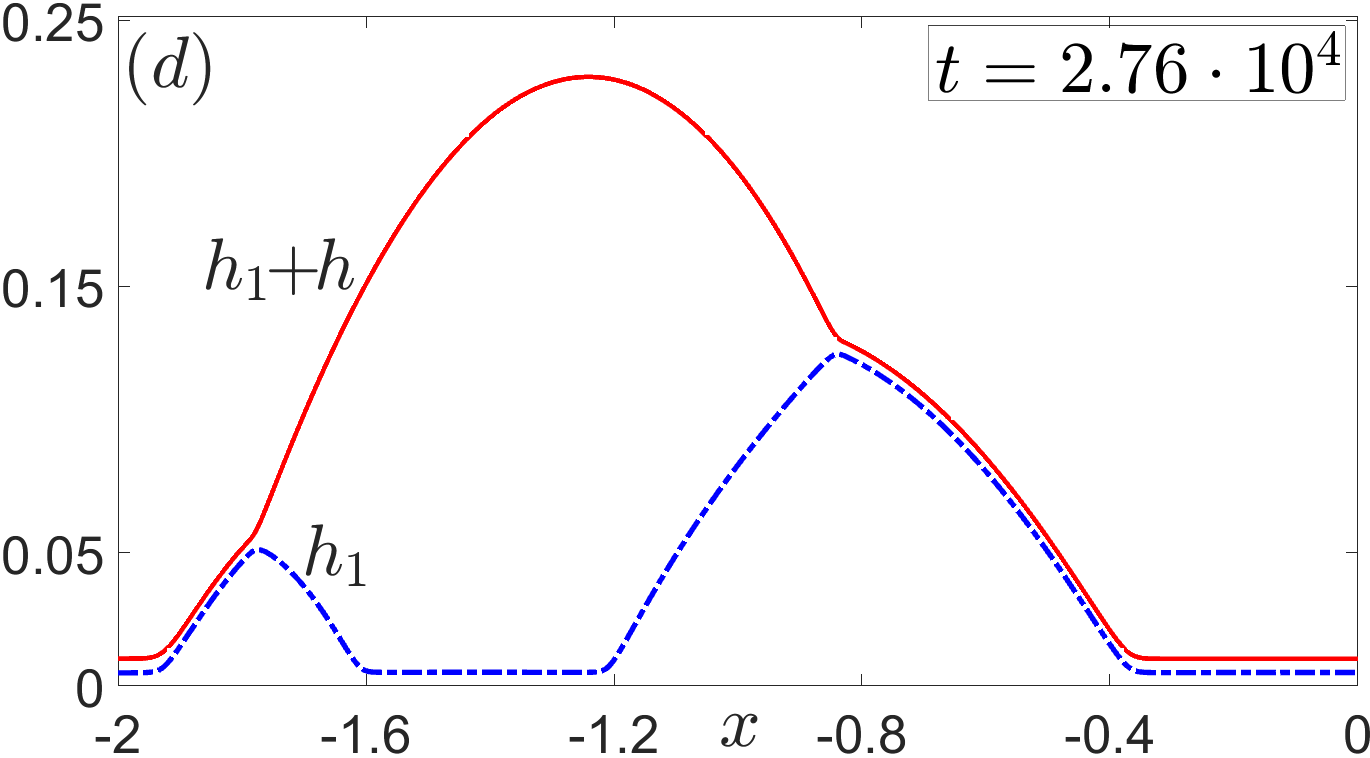}  	
	\hspace{.1cm}\includegraphics[width=.38\textwidth]{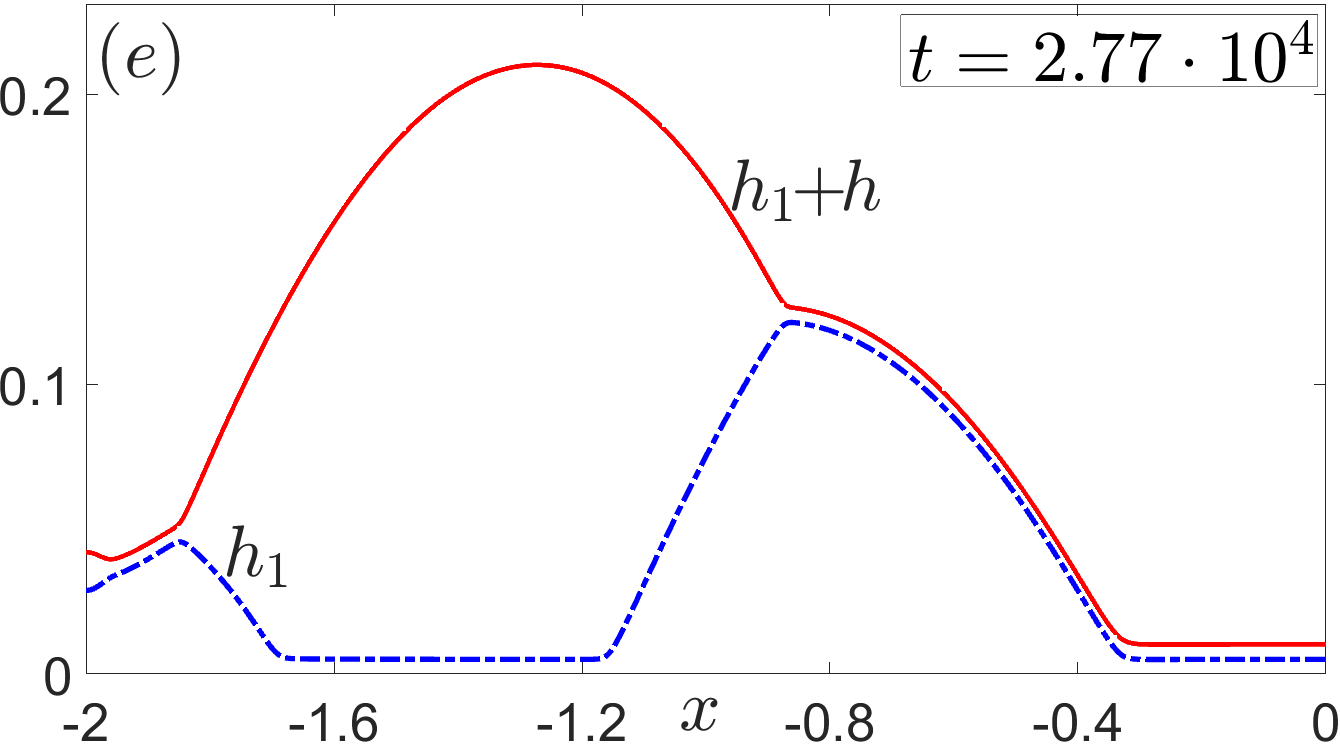}
	\hspace{.1cm}\includegraphics[width=.38\textwidth]{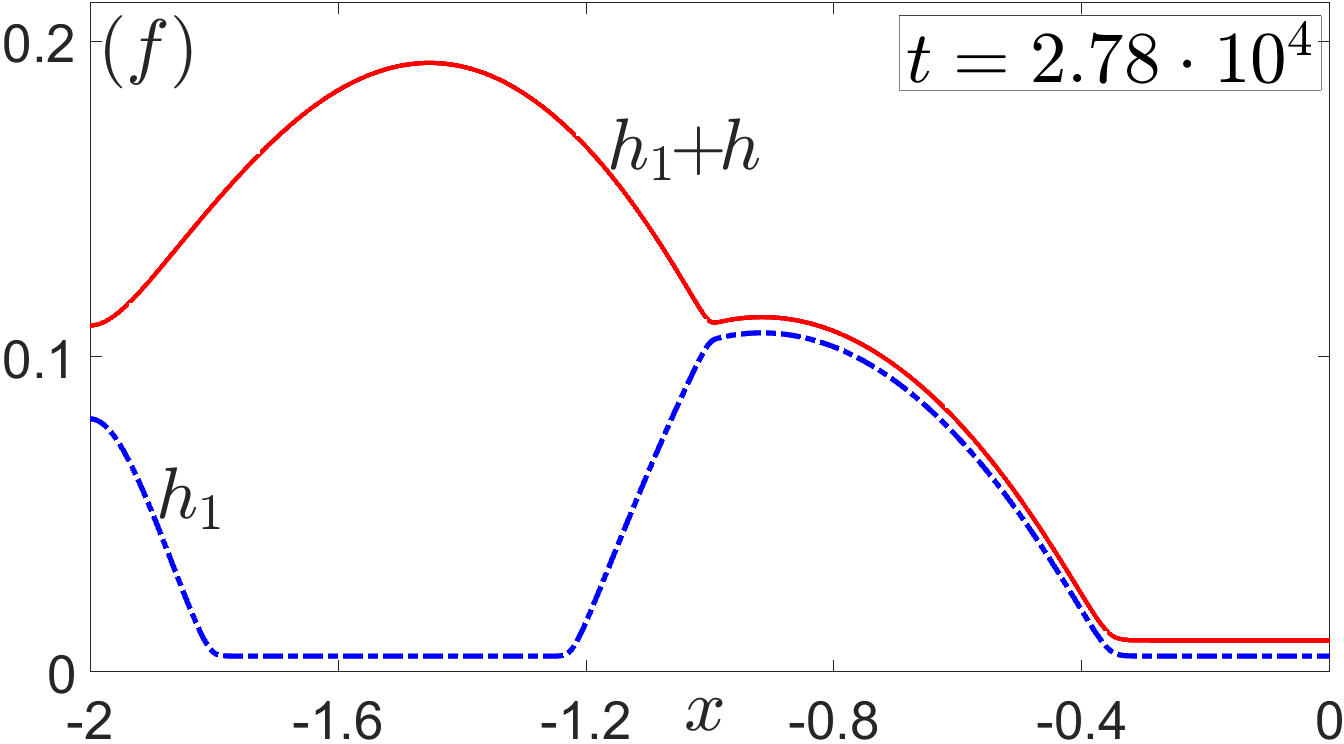}\\[1ex]	
	\hspace{-2.4cm}\includegraphics[width=.38\textwidth]{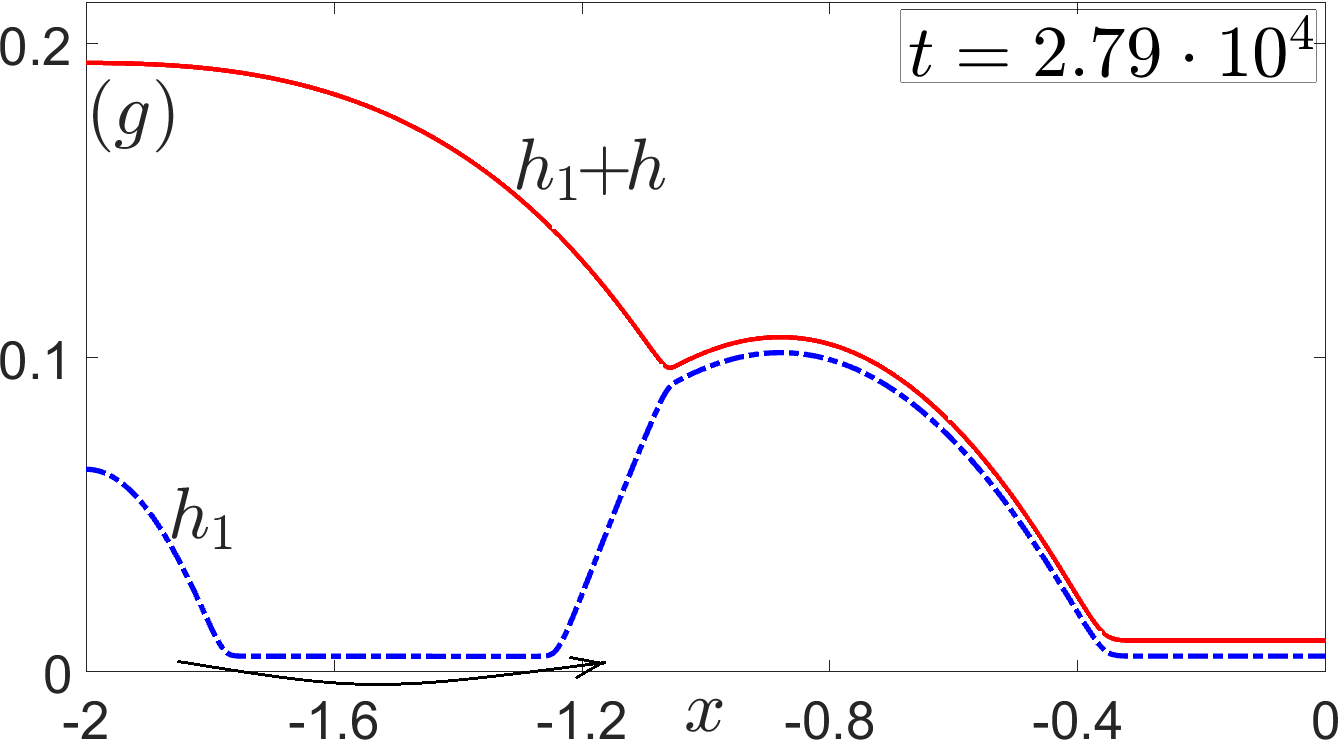}  	
	\hspace{.1cm}\includegraphics[width=.38\textwidth]{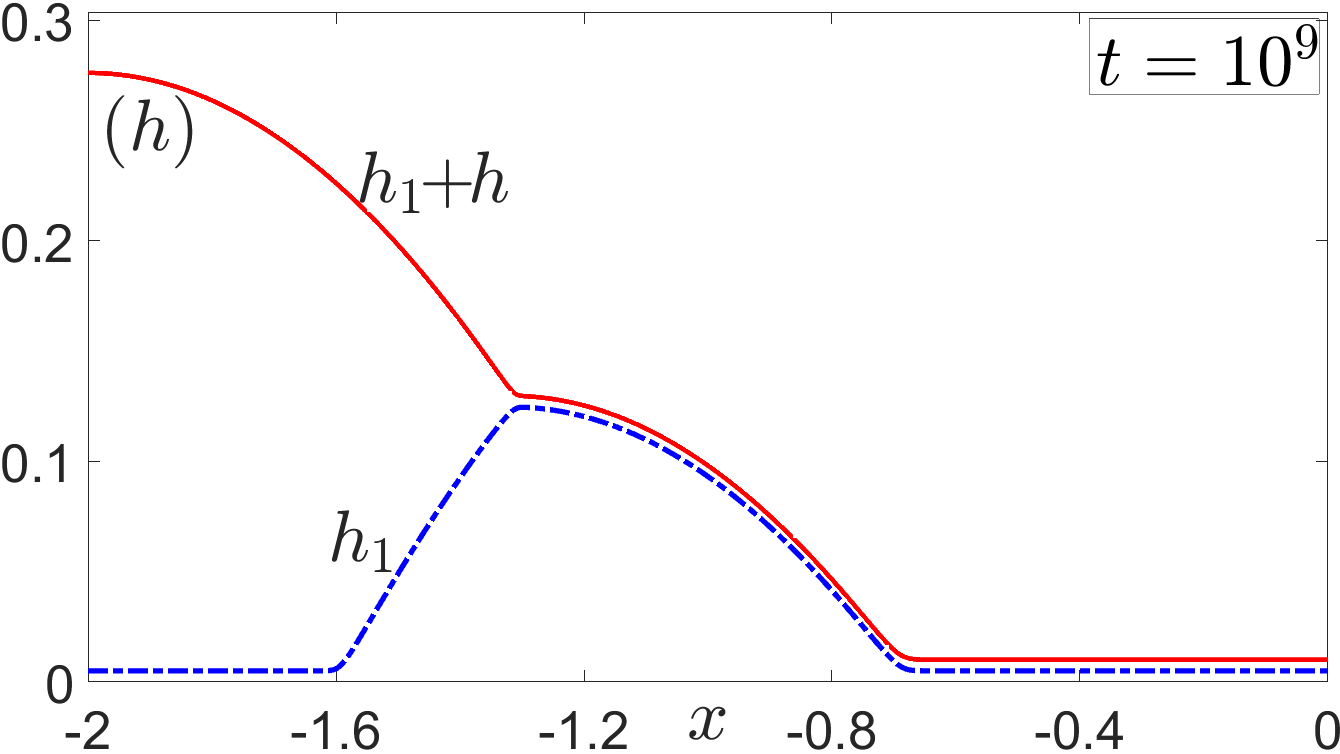}
	\hspace{.1cm}\includegraphics[width=.38\textwidth]{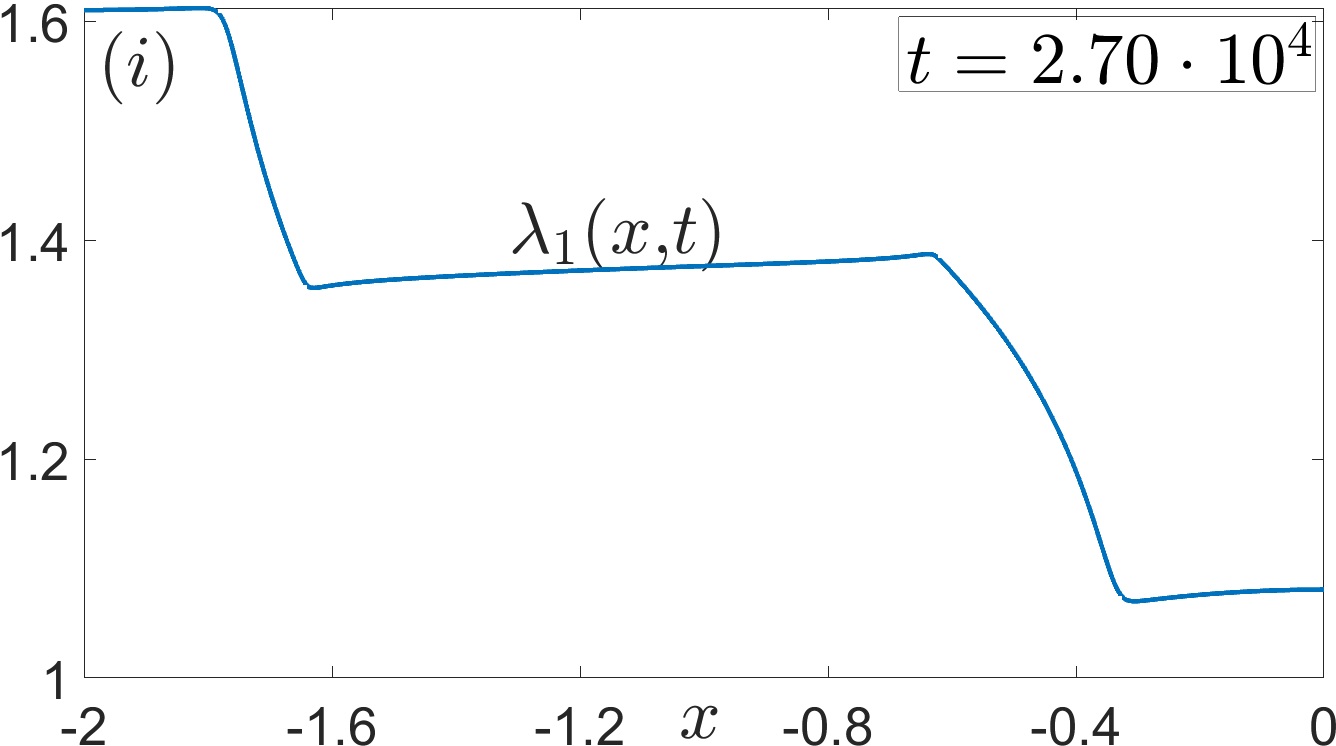}				
	
	\caption{\small {\bf(a)}--{\bf(h)}: different time snapshots of coarsening numerical solution to \rf{BS}--\rf{BC} considered with $L\!\!=\!\!2,\,\sigma\!\!=\!\!1.2,\,\eps\!\!=\!\!0.005$. Initial condition \rf{IC} is symmetric {\it sessile lens} solution $(\bf2^-0^-02)$ to stationary system \rf{SSa}--\rf{SSc}. {\bf(i)}: numerical  pressure profile $\lambda_1(x,t)\!\!=\!\!\Pi_\eps(h(x,t))\!-\!h_{xx}(x,t)\!-\!h_{1,xx}(x,t)$ for {\bf(b)}.}
\end{figure}
%%%%%%%%%%%%%%%%%%%%%%%%%%%%%%%%%%%%%%%%%%%%%%%%%%%%%%%%%%%%%%%%%%%%%%%%%%%%%%%%%%%%%%%%%%%%%%%%%%%%%%%%%%%%
\noindent Finally, $\bf (11^-03)$ coarsens via the induced mass flux in $h_1$ layer (depicted by an arrow in Fig.15 {\bf (g)}) from the smaller internal drop towards the larger one and converges in the long time to stationary {\it $h$-sessile zig-zag} $\bf(1^-02)$ (cf. Fig.15 {\bf (h)}). 

In Fig.16, initial profile \rf{IC} is  a symmetric version of {\it $h_1$-sessile zig-zag} solution ($\bf 3^-1^-00^-13$) having six CLs. The corresponding numerical solution to system \rf{BS}--\rf{BC} coarsens by drifting to $x=0$ and showing small symmetry break of its shape (cf. Fig.16 {\bf (b)}) reflected by non-constant pressure profile $\lambda_2(x,\cdot)$ with two jumps in the regions of two lateral drops (cf. Fig.16 {\bf (f)}). After the boundary touch a new {\it lens subdrop} nucleates at $x=0$ (cf. Fig.16 {\bf (c)-(d)}) and a four-CL solution $\bf (3^-1^-00^-)$ is formed. Finally, the latter coarsens via the induced mass flux  in $h$ layer (depicted by an arrow in Fig.16 {\bf (d)}) from the small lateral drop towards the lens and converges to stationary {\it $h$-sessile lens} (inverted in $x$ variable) $\bf(2^-0^-)$  (cf. Fig.16 {\bf (e)}). 

From the simulations presented in Fig.14--16 we conclude that whenever a composite  stationary solution to \rf{BS}-\rf{BC} contains as its part a reflected elementary one-CL one the former is unstable under small perturbations and coarsens via
the induced mass flux in either $h_1$ or $h$ layer having $O(\eps)$ thickness. Typical such examples are coarsening of five-CL solution  $\bf (0^-11^-02)$ (the reflected part is ${\bf(11^-)}$) in Fig.14 {\bf(d)--(e)}; four-CL one $\bf (11^-03)$  (the reflected part is $(\bf11^-)$) in Fig.15 {\bf(g)--(h)}; and four-CL one $(\bf 3^-1^-00^-)$ (the reflected part is ${\bf(00^-}$)) in Fig.16 {\bf(d)--(e)}.

Simulations shown in Fig.17 independently confirm the fact that other composite solutions (beside the eleven ones considered in sections $3-6$) are dynamically unstable. In Fig.17 {\bf (a)}, we first derived the leading order profile for
three-CL
%%%%%%%%%%%%%%%%%%%%%%%%%%%%%%%%%%%%%%%%%%%%%%%%%%%%%%%%%%%%%%%%%%%%%%%%%%%%%%%%%%%%%%%%%%%%%%%%%%%%%%%%%%%%
\begin{figure}[H] 
	\centering
	\vspace{-2.9cm}
	\hspace{-2.4cm}\includegraphics[width=.38\textwidth]{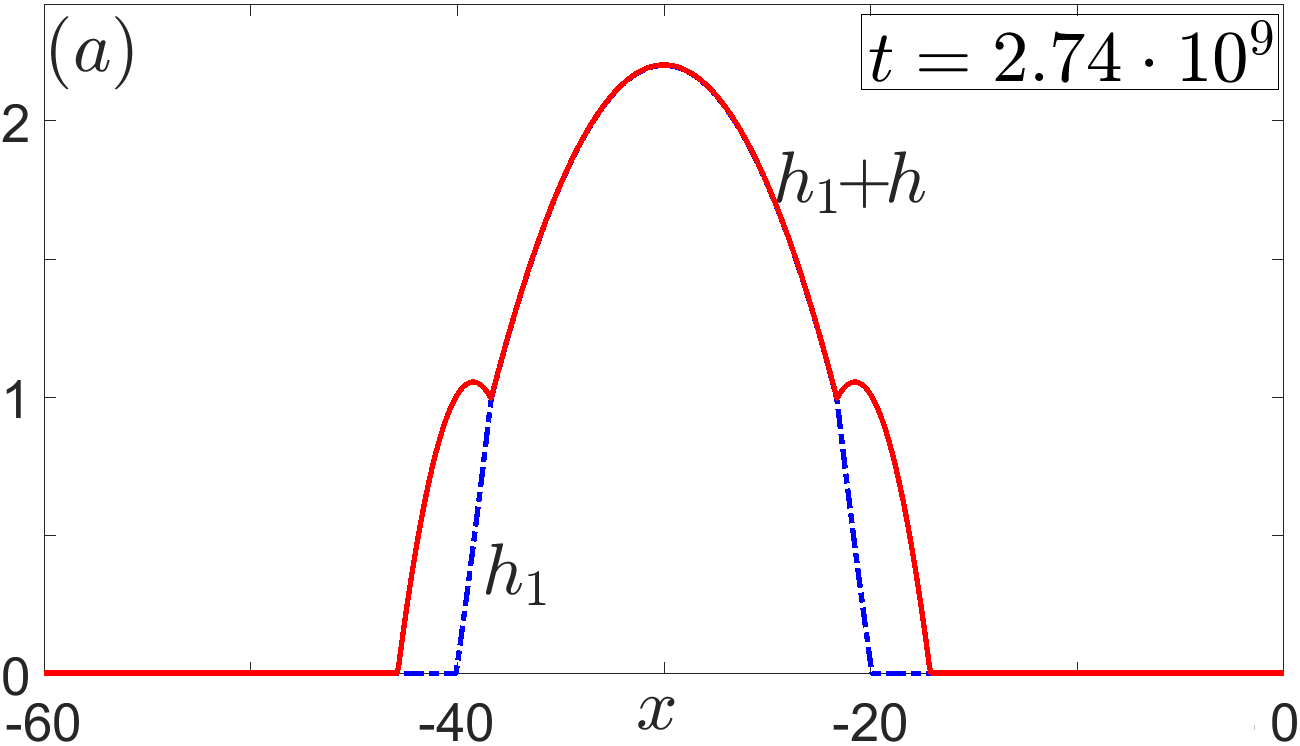}  	
	\hspace{.1cm}\includegraphics[width=.38\textwidth]{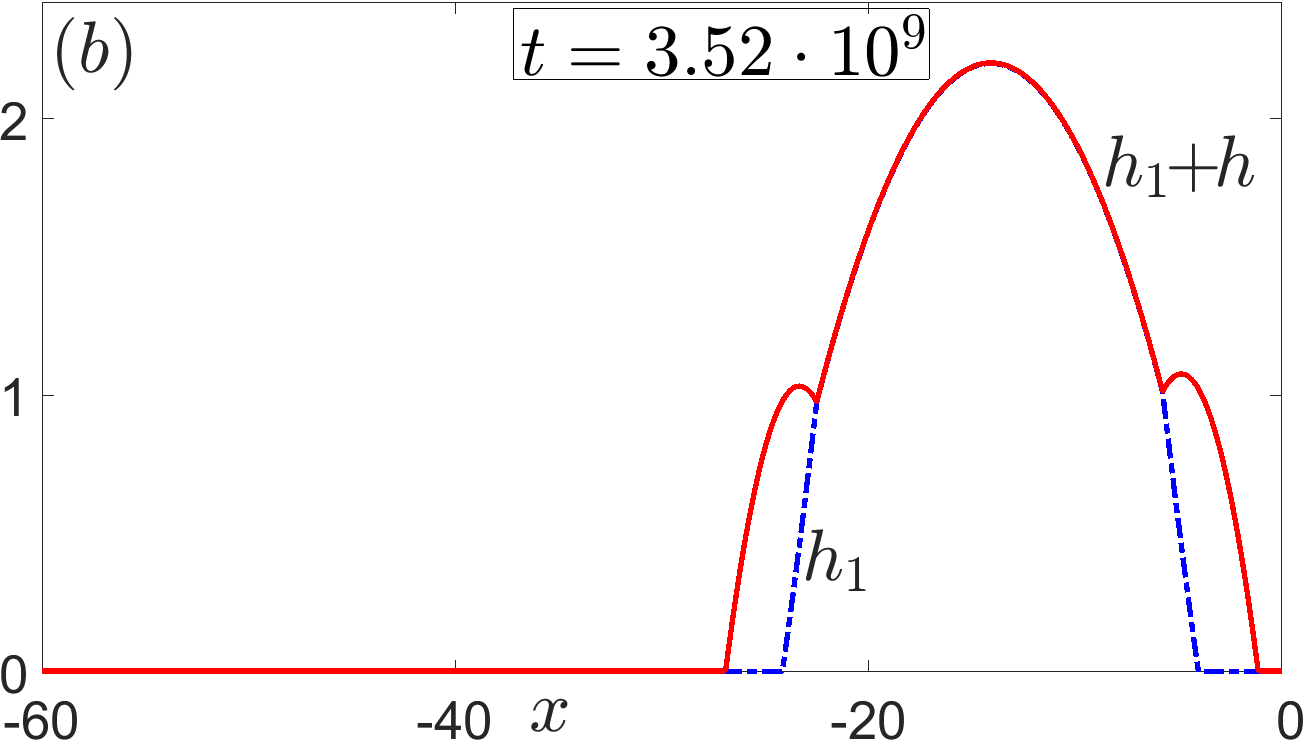}
	\hspace{.1cm}\includegraphics[width=.38\textwidth]{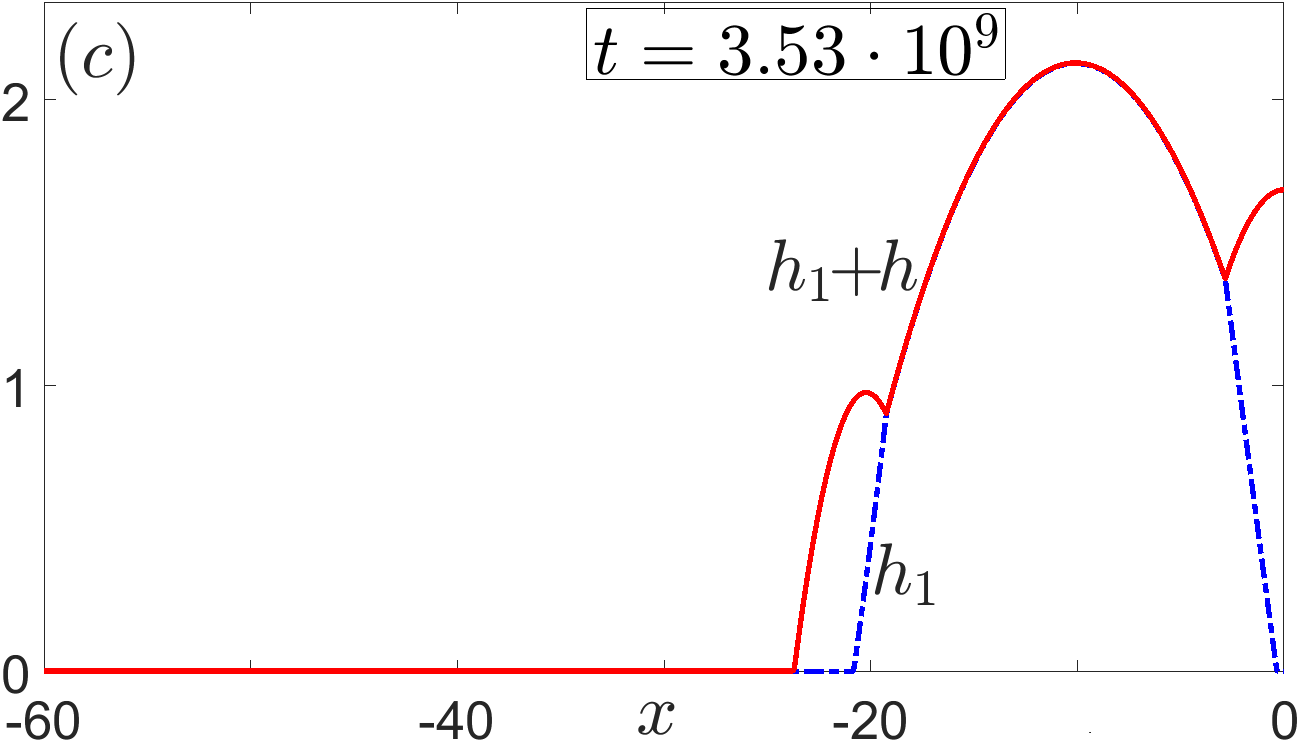}\\[1ex]	
	\hspace{-2.4cm}\includegraphics[width=.38\textwidth]{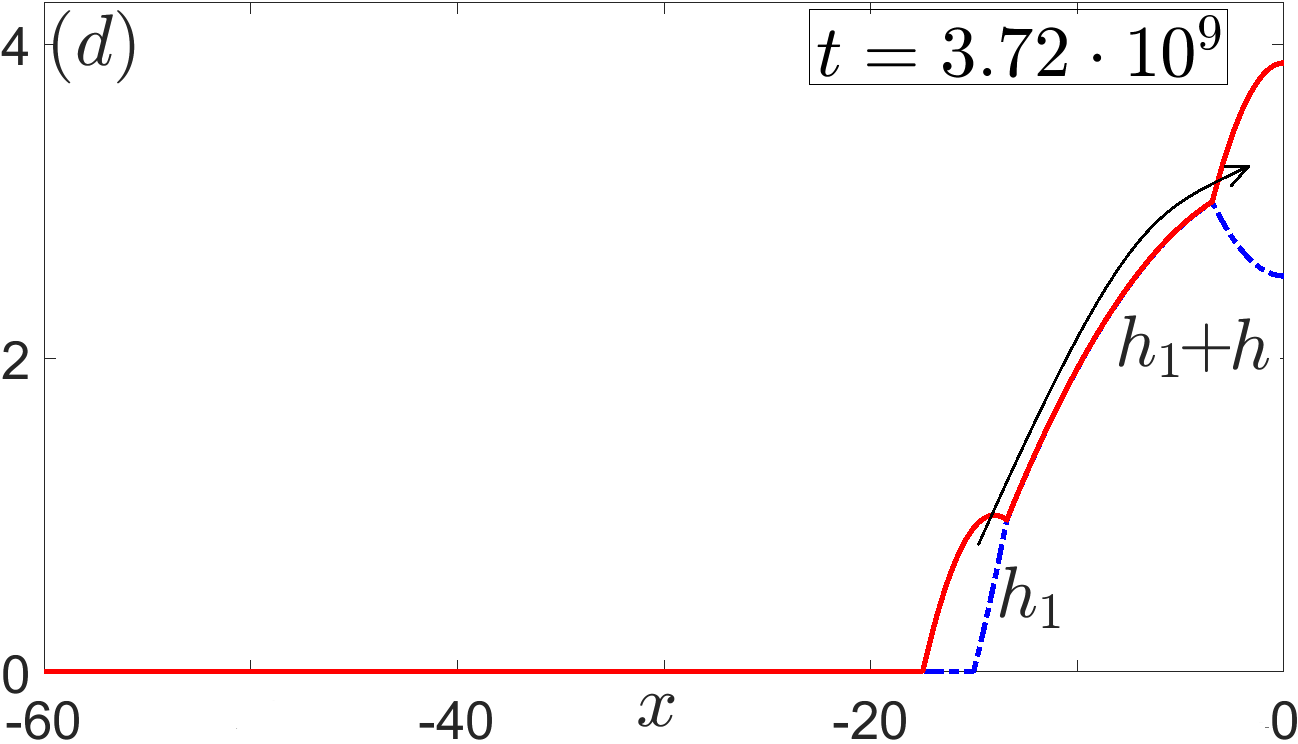}  	
	\hspace{.1cm}\includegraphics[width=.38\textwidth]{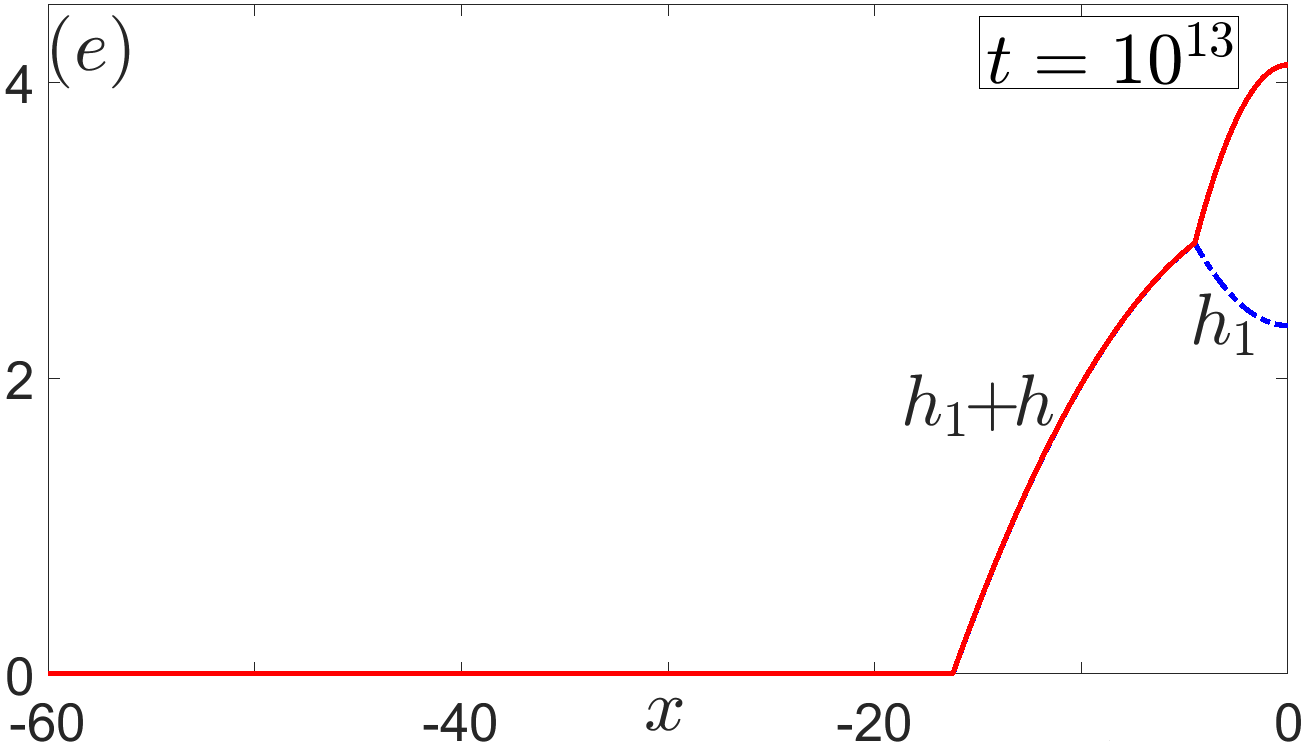}
	\hspace{.1cm}\includegraphics[width=.38\textwidth]{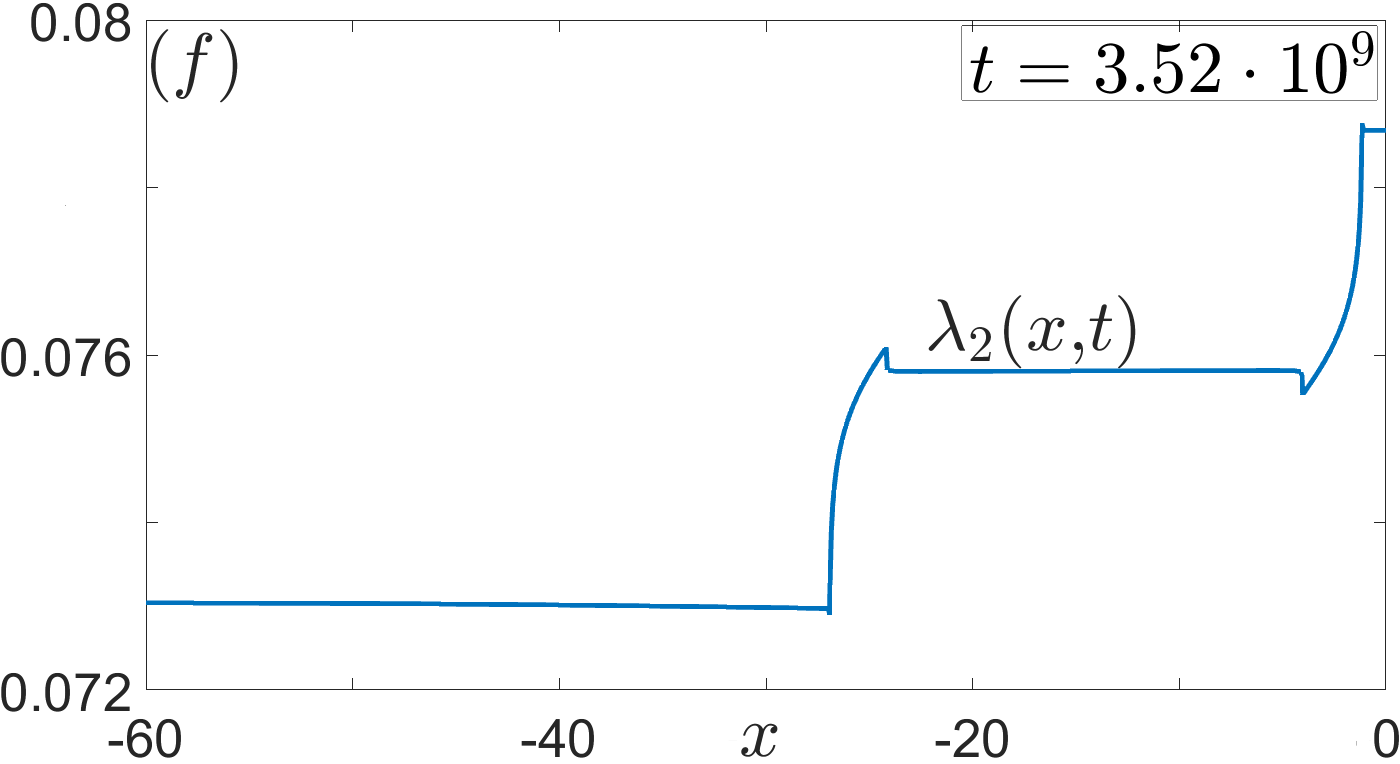}			
	
	\caption{\small{\bf(a)}--{\bf(e)}: different time snapshots of coarsening numerical solution to \rf{BS}--\rf{BC} considered with $L\!\!=\!\!60,\,\sigma\!\!=\!\!1.2,\,\eps\!\!=\!\!0.002$. Initial condition \rf{IC} is symmetric {\it $h_1$-sessile zig-zag} solution $(\bf 3^-1^-00^-13)$ to stationary system \rf{SSa}--\rf{SSc}. {\bf(f)}: numerical pressure profile $\lambda_2(x,t)\!\!=\!\!\Pi_\eps(h_1(x,t))\!-\!h_{xx}(x,t)\!-\!(\sigma+1)h_{1,xx}(x,t)$ for {\bf(b)}.}
\end{figure}
%%%%%%%%%%%%%%%%%%%%%%%%%%%%%%%%%%%%%%%%%%%%%%%%%%%%%%%%%%%%%%%%%%%%%%%%%%%%%%%%%%%%%%%%%%%%%%%%%%%%%%%%%%%%
\noindent composite solution  $(\bf1^-0^+0^-)$ to stationary system \rf{SSa}--\rf{SSc} and then set it as initial condition \rf{IC} for system \rf{BS}--\rf{BC}. Note that this composite solution can be obtain by asymptotic matching of {\it zig-zag} and {\it lens} ones at the $3^{\mathrm{rd}}$ CL $x=-s$ and the resulting system of the leading order matching conditions, in fact, decouples. Accordingly, for $(\bf1^-0^+0^-)$ pressures $\lambda_1^0$ and $\lambda_0^2$  obey formulae \rf{lambda_r}--\rf{lambda2_zz}  with  $h^m=h^0(-L)$ and $h_1^m$ being the maximum of the parabolic profile $h_1^0(x)$ in the internal region where $h(x)=O(\eps)$. Also the positions of the first two CLs for $(\bf1^-0^+0^-)$ are given by {\it zig-zag} expressions \rf{sx_expr} (up to replacement $(s_1,\,s)\go(s_2,\,s_1)$ there) and, additionally, the following three parameters describe the {\it lens} part in Fig.17 {\bf (a)}:
\bes
s=-\tfrac{\sqrt{2(\sigma+1)\sigma|\phi(1)|}}{\lambda_2^0-(\sigma+1)\lambda_1^0},\quad h^0(0)=-\tfrac{(\sigma+1)|\phi(1)|}{\lambda_2^0-(\sigma+1)\lambda_1^0},\quad h_1^0(0)=h_1^m+\tfrac{|\phi(1)|}{\lambda_2^0-(\sigma+1)\lambda_1^0}.
\ees
The corresponding numerical solution to \rf{BS}--\rf{BC} with \rf{IC} given by thus described leading order profile of $(\bf1^-0^+0^-)$ exhibits coarsening due to the induced flux in $h$-layer (depicted by an arrow in Fig.17 {\bf (b)}) reflected by the pressure variations (cf. Fig.17 {\bf (f)}). By that the {\it zig-zag} part is transformed into a small {\it lens}, while at the same time the large one touches and spreads over the substrate with a small internal drop being formed at $x=0$ (cf. Fig.17 {\bf (b)--(d)}). In the long time both the small {\it lens} and internal drop are absorbed into the bulks and the solution converges to stable stationary {\it zig-zag} (inverted in $x$ variable)  one $\bf(10^-)$ (cf. Fig.17 {\bf (e)}). 

%%%%%%%%%%%%%%%%%%%%%%%%%%%%%%%%%%%%%%%%%%%%%%%%%%%%%%%%%%%%%%%%%
\section{Discussion}
%%%%%%%%%%%%%%%%%%%%%%%%%%%%%%%%%%%%%%%%%%%%%%%%%%%%%%%%%%%%%%%%%
The results of this study lay a mathematical background and provide insight for understanding of the shapes and structure of solutions to the liquid bilayer models considered in~\cite{NTGDC12,BSTA21,MAP01,PBMT04,PBMT05} and potentially to some recently developed ones~\cite{ATTGA24,HDGT24,PM23,POB21}. Here we derived explicit formulae for the leading order profiles of eleven types of stationary solutions to one-dimensional bilayer system \rf{BS}--\rf{IC} 
%%%%%%%%%%%%%%%%%%%%%%%%%%%%%%%%%%%%%%%%%%%%%%%%%%%%%%%%%%%%%%%%%%%%%%%%%%%%%%%%%%%%%%%%%%%%%%%%%%%%%%%%%%%%
\begin{figure}[H] 
	\centering
	\vspace{-2.9cm}
	\hspace{-2.4cm}\includegraphics[width=.38\textwidth]{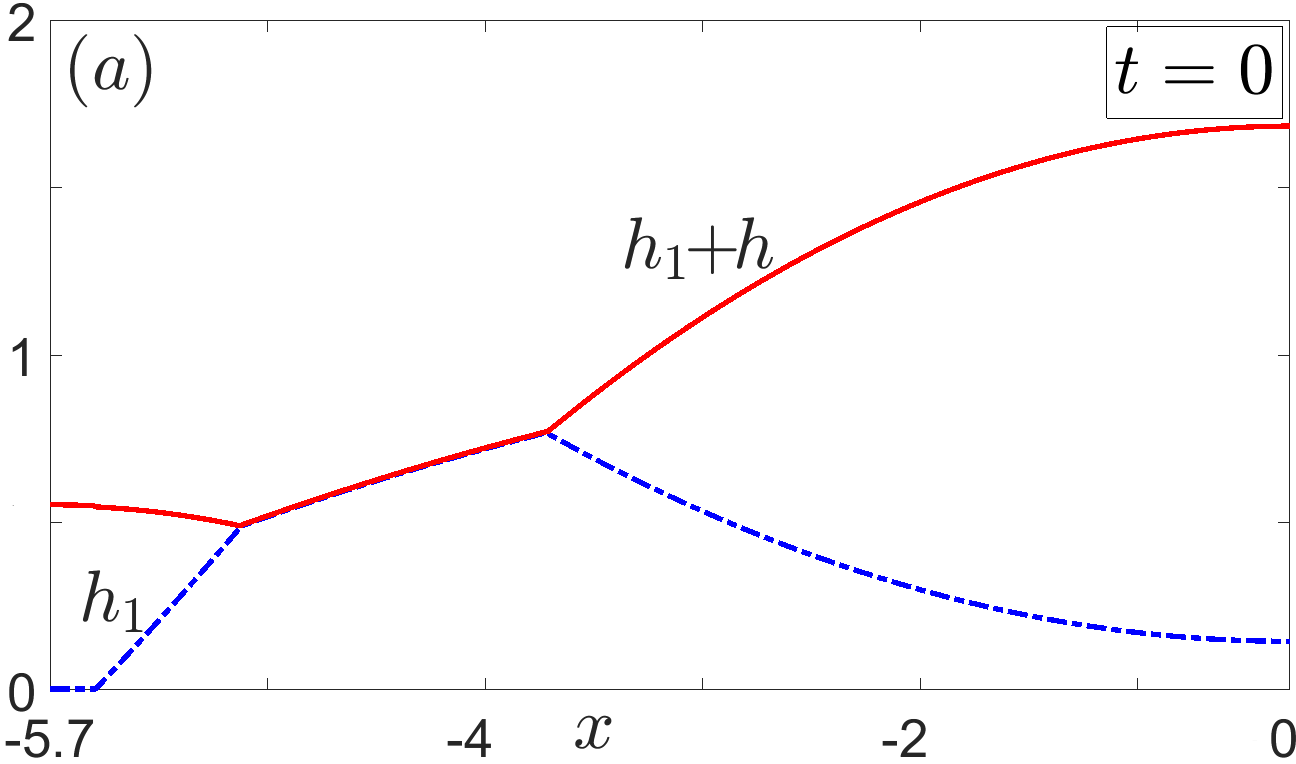}  	
	\hspace{.1cm}\includegraphics[width=.38\textwidth]{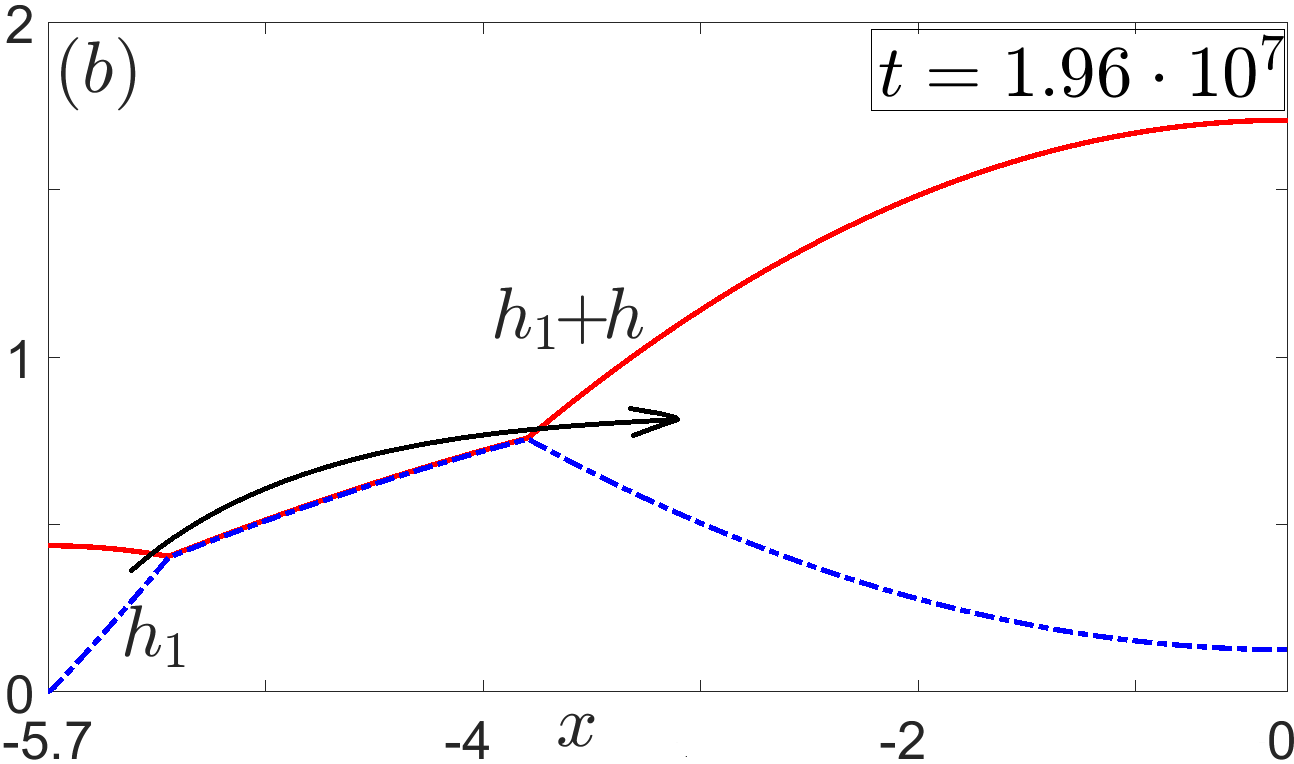}
	\hspace{.1cm}\includegraphics[width=.38\textwidth]{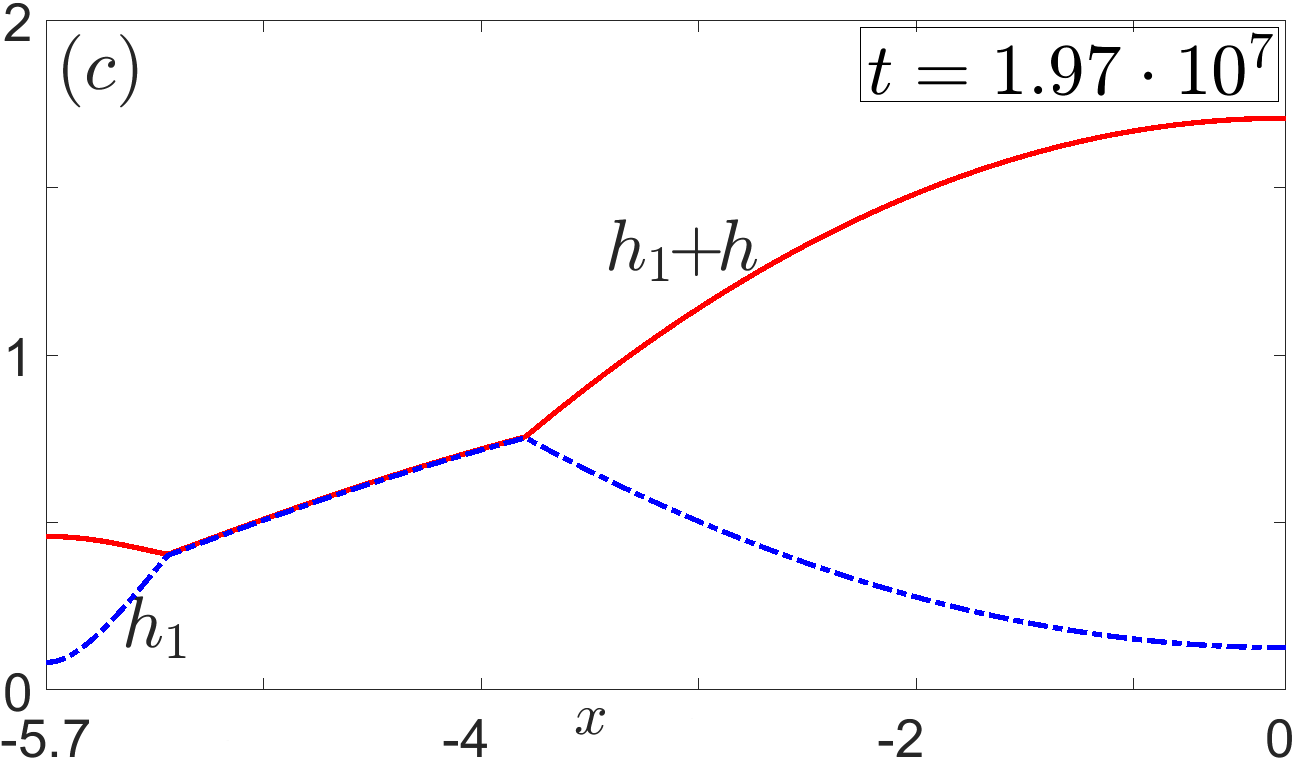}\\[1ex]	
	\hspace{-2.4cm}\includegraphics[width=.38\textwidth]{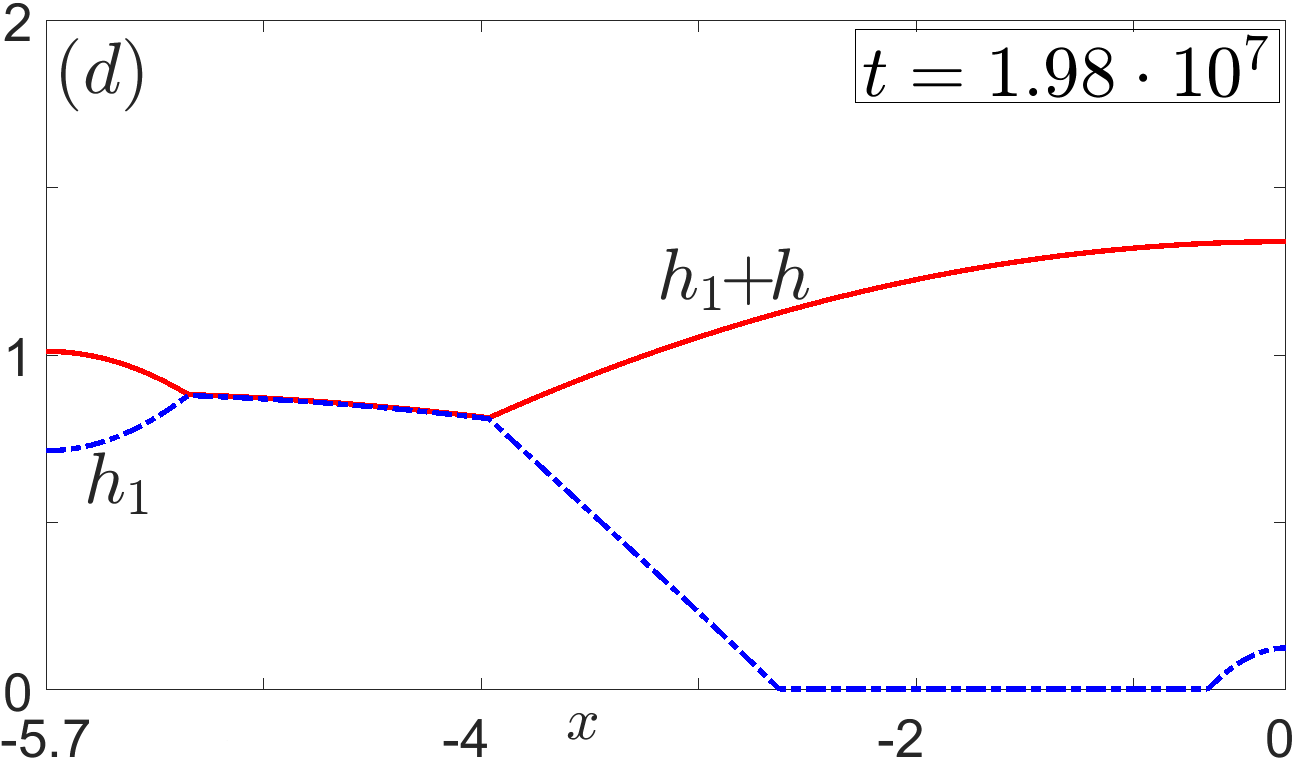}  	
	\hspace{.1cm}\includegraphics[width=.38\textwidth]{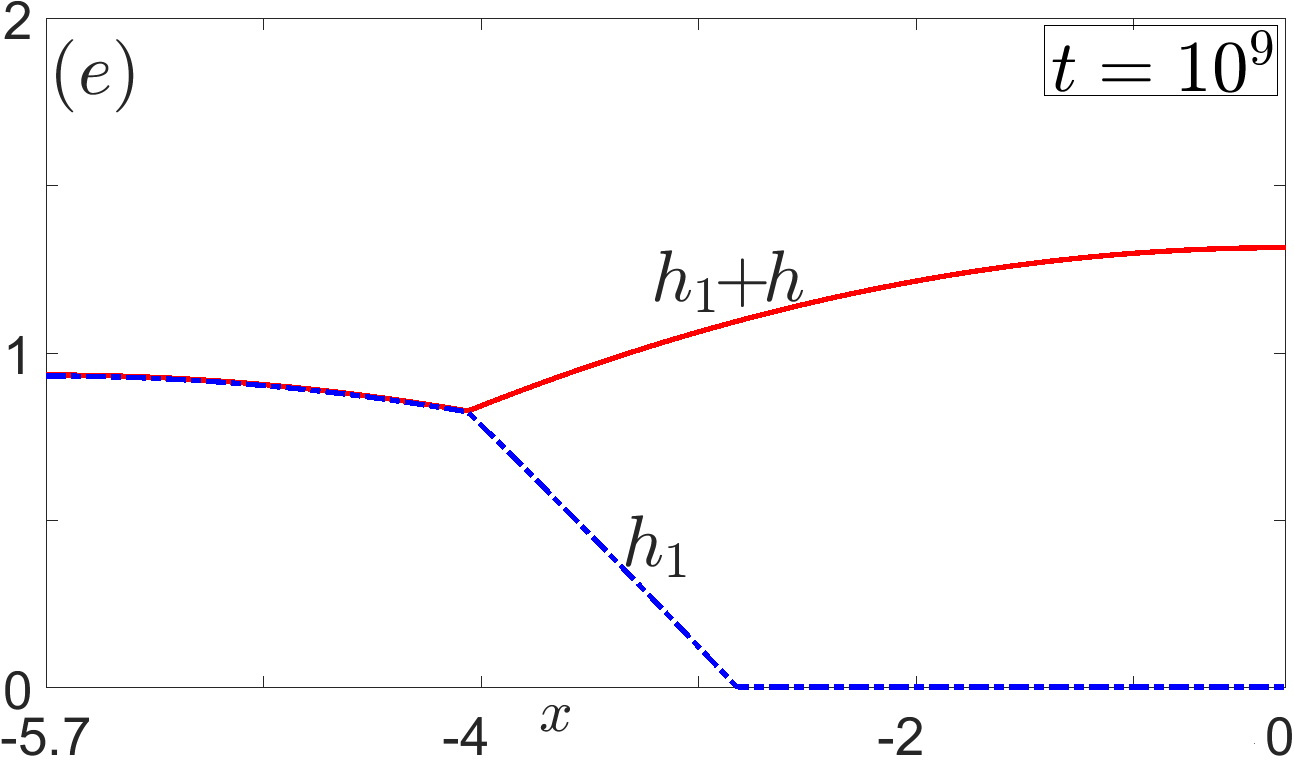}
	\hspace{.1cm}\includegraphics[width=.38\textwidth]{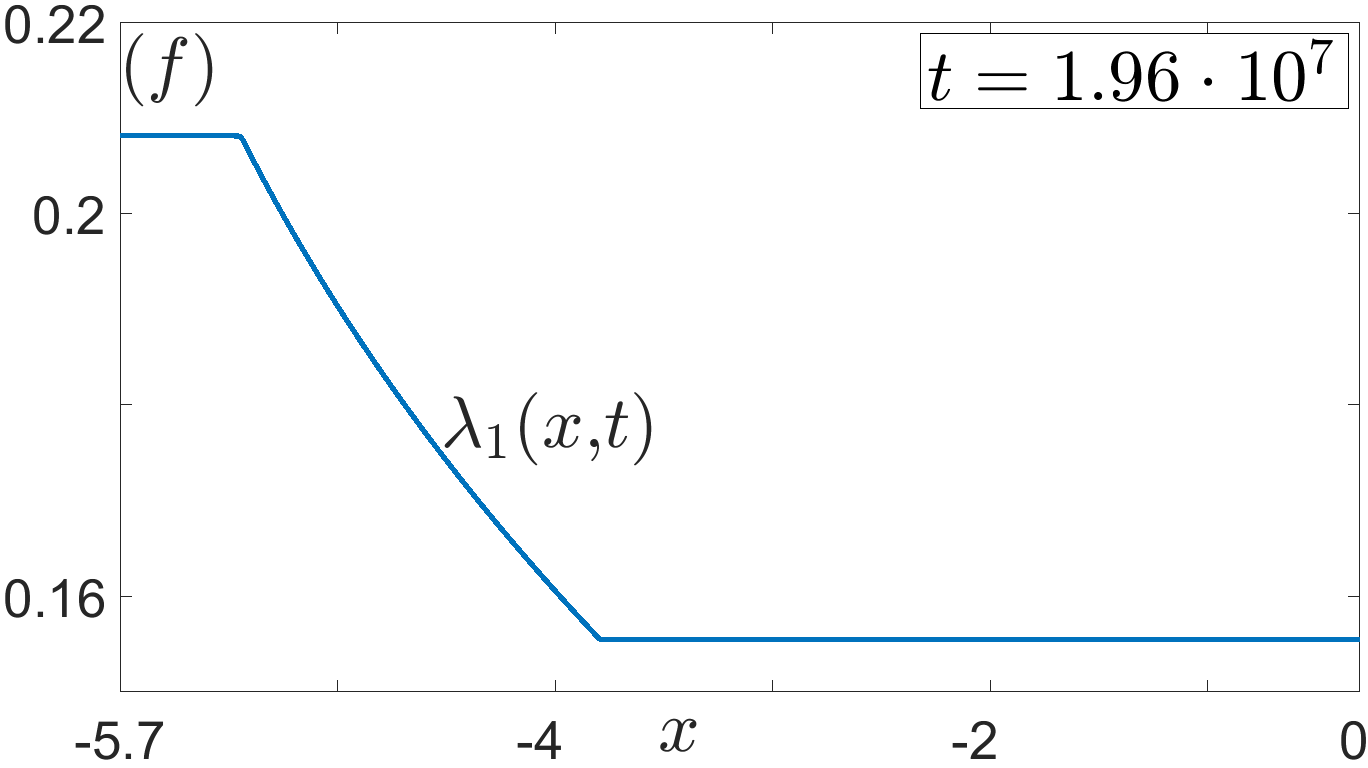}			
	
	\caption{\small{\bf(a)}--{\bf(e)}: different time snapshots of coarsening numerical solution to \rf{BS}--\rf{BC} considered with $L\!\!=\!\!5.7,\,\sigma\!\!=\!\!0.7,\,\eps\!\!=\!\!0.003$. Initial condition \rf{IC} is composite solution $(\bf1^-0^+0^-)$ to stationary system \rf{SSa}--\rf{SSc}. {\bf(f)}: numerical pressure profile $\lambda_1(x,t)\!\!=\!\!\Pi_\eps(h(x,t))\!-\!h_{xx}(x,t)\!-\!h_{1,xx}(x,t)$ in {\bf(b)}.}
\end{figure}
%%%%%%%%%%%%%%%%%%%%%%%%%%%%%%%%%%%%%%%%%%%%%%%%%%%%%%%%%%%%%%%%%%%%%%%%%%%%%%%%%%%%%%%%%%%%%%%%%%%%%%%%%%%%
\noindent  considered with intermolecular potential \rf{phi_c} depending on both layer heights. These solutions were asymptotically matched as compositions of several elementary ones with triple phase contact lines. Numerically solving  \rf{BS}--\rf{IC} we showed that most of these solutions are dynamically stable or  weakly translationally unstable. Still other composite stationary solutions to \rf{BS}--\rf{IC} turn out to be numerically unstable and experience slow time coarsening accompanied by the complex morphological transformations (cf. Fig.14--17). 

Interestingly, our results stay in a good correspondence with the observations collected by modeling and experimental studies of thin bilayer systems~\cite{NTGDC12,BSTA21,MAP01,PBMT04,PBMT05}. In~\cite{NTGDC12}, authors produced several experiments with deposition of two immiscible fluid (typically mercury and water) on a hydrophilic glass substrate and observed four types of {\it compound sessile drops}, which they termed as {\it encapsulated, lens, collar and Janus ones} (cf. Fig.1 of \cite{NTGDC12}). Some of the found solutions to \rf{BS}--\rf{IC} partially resemble the two-dimensional compound drop shapes shown in~\cite{NTGDC12,MAP01} despite the fact that the models used in the latter articles are sharp interface ones explicitly not accounting for possible intermolecular interactions in the experimental setups. Indeed, our internal drop solutions are similar to the encapsulated ones of~\cite{MAP01}, lens solutions look the same as in~\cite{NTGDC12}, our three type of sessile zig-zags resemble well the {\it Janus drops shapes}, and a collar type solution is observed in Fig.15. Mathematically seen we prefer to distinguish between sessile and bulk type solutions, because solving the leading order systems of matching conditions for them proceeds differently (e.g. compare solution algorithms to systems \rf{1m0sol} and \rf{2m0m13sol} for {\it zig-zag} and {\it 2-side sessile zigzag} in sections $4-5$).

Next, by comparing the coarsening patterns of section $8$ with the ones presented and analyzed for one-dimensional bilayer systems in~\cite{PBMT04,PBMT05} we confirm the presence of two driving {\it coarsening modes} observed there: varicose or lens (cf. Fig. 16) and zig-zag (cf. Fig. 14) ones, as well as possible transitions between them (cf. Fig. 17). Besides, we observed a new coarsening mode (not reported in~\cite{PBMT04,PBMT05}) associated with formation and interaction of internal drops (cf. Fig.15 {\bf (f)-(h)} and Fig.17 {\bf (d)-(e)}). Appearance of this new coarsening mode might be connected with our choice of algebraic decay rates in the repulsive terms of intermolecular potential \rf{phi}.  

We conclude the article by stating several interesting open questions:
\begin{itemize}
	\item[{\bf (a)}] For $h_1$ and $h$-drop stationary solutions described in section 3 we found that one of two {\it hydrodynamic pressures} $\lambda_1^0$ or  $\lambda_2^0$ can not be determined from their leading order systems of matching conditions, while in numerical simulations of system \rf{BS}--\rf{IC} they are observed being negative (cf. Fig.4-5). This unconventional observation should be justified from both physical and mathematical sides.
	
	\item[{\bf (b)}] The method for asymptotic matching of stationary solutions to \rf{BS}--\rf{IC} presented in sections $2-6$ can be used for constructing the composite ones with any arbitrary number $N$ of triple contact lines. Beside the stability one, there exists a question whether the solution existence domains (EDs) in the space of {\it model parameters} (similar to those presented in Fig.13 of section $7$) shrink to empty sets at some finite $N$ or not.
	
	\item[{\bf (c)}] Bifurcation analysis of the stationary solutions to \rf{BS}--\rf{IC} in the spirit of studies~\cite{BGW01,Zh09} is yet lacking, though the first step towards it is done in diagrams of Fig.13. In particular, our numerical simulations do not reveal solutions with four-phase merging contact lines predicted in~\cite{MAP01}, which still may belong to unstable bifurcation paths. 
	
	\item[{\bf (d)}] Our results lay mathematical background for further analysis of the coarsening dynamics in bilayer system following similar approaches to those ones developed for one-layer thin film equations~\cite{GW03,GW05,GORS09,KW10,KRW11,Ki14}.		
\end{itemize}
%%%%%%%%%%%%%%%%%%%%%%%%%%%%%%%%%%%%%%%%%%%%%%%%%%%%%%%%%%%%%%%%%
\section{Acknowledgments}
%%%%%%%%%%%%%%%%%%%%%%%%%%%%%%%%%%%%%%%%%%%%%%%%%%%%%%%%%%%%%%%%%
The author would like to thank Dirk Peschka and Barbara Wagner for our stimulating discussion of possible forms for intermolecular potential \rf{phi_c}.

\end{document}